\documentclass[11pt]{amsart}

\usepackage[margin=3cm]{geometry}
\usepackage{tikz}
\usepackage{tikz-cd}
\usepackage{amsfonts}
\usepackage{amsmath}
\usepackage{mathtools}
\usepackage{amsthm}
\usepackage{amssymb}
\usepackage{amscd}
\usepackage{graphicx}
\usepackage{subcaption}
\usepackage{faktor}
\usepackage[UKenglish]{babel}
\usepackage{indentfirst}
\usepackage[T2A]{fontenc}
\usepackage{chngcntr}
\usepackage{xcolor}
\usepackage{bm}
\usepackage{todonotes}
\usepackage{comment}
\usepackage[shortlabels]{enumitem}
\usepackage{yfonts}
\usepackage{mathabx} %for big asterisk
\usepackage{relsize} %for resizing
\usepackage{bbm}

\usepackage{array}

\usepackage{breqn}

\usepackage{geometry} 

\usepackage{tabu}

\usepackage{tikz-cd}
\usepackage{multicol}

\usetikzlibrary{matrix,fit}

\usepackage{extarrows}

\usepackage{color}\definecolor{darkblue}{rgb}{0,0.1,.5}
\usepackage[colorlinks=true,linkcolor=darkblue, urlcolor=darkblue, citecolor=darkblue]{hyperref}
\usepackage{cleveref}

\newcommand{\CP}{\mathbb{C}P^{\infty}}

\newcommand{\cK}{\mathcal{K}}
\newcommand{\cL}{\mathcal{L}}

\DeclareMathOperator*{\colim}{colim}
\DeclareMathOperator*{\hocolim}{hocolim}

\DeclareMathAlphabet\mathbfcal{OMS}{cmsy}{b}{n}

\let\epsilon\varepsilon

\newcommand{\CPi}{\mathbb{C} P ^{\infty}} %Complex projective space ^infinity
\newcommand{\cT}{\mathcal{T}} %Topology on a set$
\newcommand{\cS}{\mathcal{S}} %A set%
\newcommand{\DJK}{DJ_{\cK}} %DJ_K
\newcommand{\biggast}{\mathlarger{\mathlarger{\mathlarger{\bigast}}}} %Big asterisk for join

\newcommand{\joinspc}{\underset{i=1}{\overset{m}{\biggast}}}

\newcommand{\SX}{\Sigma X} %Suspension of a space X

\newcommand{\hooklongrightarrow}{\lhook\joinrel\longrightarrow}

\DeclareMathOperator{\sk}{sk}

\newcommand{\hr}[2][]{\hyperref[#2]{#1~\ref{#2}}}

\newtheorem{theorem}{Theorem}[section]
\newtheorem{proposition}[theorem]{Proposition}
\newtheorem{lemma}[theorem]{Lemma}
\newtheorem{corollary}[theorem]{Corollary}

\theoremstyle{remark}

\theoremstyle{definition}
\newtheorem{example}[theorem]{Example}

\newtheorem{definition}[theorem]{Definition}

\counterwithout{figure}{section}

\title{Relations among higher Whitehead maps}
\begin{document}

%\begin{center} \today \end{center}

\author{Jelena Grbi\'{c}}
\address{School of Mathematical Sciences, University of Southampton, Southampton, UK }
\email{j.grbic@soton.ac.uk}

\author{George Simmons}
\address{School of Mathematical Sciences, University of Southampton, Southampton, UK }
\email{g.j.h.simmons@soton.ac.uk}

\author{Matthew Staniforth}
\address{School of Mathematical Sciences, University of Southampton, Southampton, UK}
\email{m.staniforth@soton.ac.uk}



\maketitle
\begin{abstract}
    We study relations among generalised higher Whitehead maps. Departing from the classical approach of studying spherical Whitehead maps via the Hurewicz homomorphism, we instead define generalised higher Whitehead maps between polyhedral products. By investigating the interplay between the homotopy-theoretic properties of polyhedral products and the combinatorial properties of simplicial complexes, we describe new families of relations among these maps, while recovering and generalising known identities among Whitehead products.
\end{abstract}

\tableofcontents

\section{Introduction}

In this paper we describe new relations in $[\Sigma X, W]$, the group of homotopy classes of maps from a suspension space $\Sigma X$ to a topological space $W$, among higher Whitehead maps.  
%we study groups of homotopy classes of maps. identifying new relations among higher Whitehead maps. 
We obtain families of relations among $n$-ary higher Whitehead products for $n \geqslant 2$, recovering the Jacobi identity and Hardie's identity in homotopy groups as a special case. 

Relations among Whitehead products and higher analogues in homotopy groups have been studied since the early 1950s. Given maps $f_i \in \pi_{d_i}(Y)$ for $i = 1,2$ with $d_i \geqslant 1$, the Whitehead product $[f_1,f_2] \in \pi_{d_1 + d_2 - 1}(Y)$ is the homotopy class of the map
\[[f_1,f_2] \colon S^{d_1-1} * S^{d_2-1} \xlongrightarrow{\rho} S^{d_1} \vee S^{d_2} \xlongrightarrow{f_1 \vee f_2} Y \vee Y \xlongrightarrow{\nabla} Y \]
where $\rho$ is the attaching map of the top cell in $S^{d_1}\times S^{d_2}$ and $\nabla$ is the fold map.
The Whitehead product equips the homotopy groups of $Y$ with the structure of a graded quasi-Lie algebra. Motivated by the Jacobi identity for Lie algebras, Blakers and Massey \cite{BlakersMasseyJacobiconjecture}, and Samelson \cite{Samelson1953ACB} independently conjectured the Jacobi identity as a relation among Whitehead products. Let $f_i \in \pi_{d_i}(Y)$ for $i=1,2,3$, where $d_i \geqslant 2$. Then
\begin{equation} \label{eq:JacobiIdentity}
    (-1)^{d_1 d_3}[f_1,[f_2,f_3]] + (-1)^{d_2d_3}[f_3,[f_1,f_2]] + (-1)^{d_1d_2}[f_2,[f_3,f_1]] = 0.
\end{equation}
This relation was established by several authors in the mid 1950s, using a variety of methods. Notably, the proof of Massey and Uehara \cite{Ma-Ue57} is algebraic in nature, and is one of the first applications of the triple Massey product, while the proof of Nakaoka and Toda \cite{Na-To54} uses a geometric approach considering relative homotopy theory. 

Hardie~\cite{Hardie61} defined the exterior Whitehead product as an operation on homotopy groups, generalising the Whitehead product. For $k \geqslant 2$, let $f_i \in \pi_{d_i}(Y_i)$ for $i = 1,\dots,k$, where $d_i \geqslant 1$. The exterior Whitehead product $\{f_1,\dots,f_k\} \in \pi_{d_1 + \cdots + d_k - 1} (FW(Y_1,\dots,Y_k))$ is the homotopy class of the map
\begin{equation*}
    \{f_1,\dots,f_k\} \colon \overset{k}{\underset{i=1}{\biggast}} S^{d_i-1} \xlongrightarrow{\rho} FW(S^{d_1},\dots,S^{d_k}) \xlongrightarrow{FW(f_1,\dots,f_k)} FW(Y_1,\dots,Y_k)
\end{equation*}
where $FW(Y_1,\dots,Y_k) = \bigcup_{i=1}^k Y_1 \times \dots \times Y_{i-1} \times \ast \times Y_{i+1} \times \cdots \times Y_k$ is the fat wedge of the spaces $Y_1,\dots,Y_k$, the map $\rho$ is the attaching map of the top cell in $S^{d_1} \times \cdots \times S^{d_k}$, see~\eqref{eq:restofprodofquot}, and $FW(f_1,\dots,f_k)$ is the map induced on the fat wedge by $f_1,\dots,f_k$. For $i=1,\dots,k$, when $d_i \geqslant 2$, Hardie showed that exterior Whitehead products satisfy
\begin{equation} \label{eq:HardieIdentity}
\sum_{i=1}^k (-1)^{\eta(i)} \{f_i,\{f_1,\dots,f_{i-1},f_{i+1},\dots,f_k\}\}=0 
\end{equation}
in $\pi_{d_1+\cdots+d_k-2} (Z)$, where $Z = \bigcup_{i=1}^k Y_i \vee FW(Y_1,\dots,Y_{i-1},Y_{i+1},\dots,Y_k)$.

%There is a gap in knowledge about relations among higher Whitehead products. The first results of this nature appear in the work of Hardie~\cite{Hardie61}, who studied relations among higher Whitehead products of spherical maps.

%\todo[inline]{Comment about these being the only known relations}

%\todo[inline,color=green]{Added a sentence to the below paragraph, George and Jelena please check it}
%Apart from Hardie's identity no systematic approach to studying relations among higher Whitehead products has been presented. 
Arkowitz~\cite[Definition~2.2]{Arkowitz1962paper} defined the generalised Whitehead product. Let $f_i \in [\Sigma X_i,Y]$ for $i=1,2$. The generalised Whitehead product $[f_1,f_2] \in [\Sigma X_1 \wedge X_2,Y]$ is the homotopy class of the composite
\[[f_1,f_2] \colon \Sigma X_1 \wedge X_2 \xlongrightarrow{\rho} \Sigma X_1 \vee \Sigma X_2 \xlongrightarrow{f_1 \vee f_2} Y \vee Y \xlongrightarrow{\nabla} Y \] 
where $\rho$ is the map whose homotopy cofibre is $\Sigma X_1 \times \Sigma X_2$.

The only known relation among generalised Whitehead products is due to Arkowitz \cite{Arkowitz1962paper} who extended the Jacobi identity~\eqref{eq:JacobiIdentity} to generalised Whitehead products.
%Although various relations among Whitehead products are known, 
This motivates us to ask whether analogues of Hardie's identity~\eqref{eq:HardieIdentity} exist in the group of homotopy classes of maps $[\Sigma X,W]$ from a suspension $\Sigma X$ to an arbitrary space $W$.

In this paper, we consider generalisations of relation~\eqref{eq:HardieIdentity} to higher Whitehead maps, in the context of polyhedral products. A polyhedral product $(\underline X, \underline A)^{\cK}$ is a topological space built out of a simplicial complex $\cK$ on $m$ vertices and an $m$-tuple of $CW$-pairs $(\underline X, \underline A)$, see Definition~\ref{def:polyprod}. It can be considered as a functor from $CW$-pairs and simplicial complexes with respect to continuous maps of $CW$-pairs and simplicial inclusions. The homotopy-theoretic properties of polyhedral products tend to be encoded by underlying combinatorial structures, and by combining geometric and combinatorial techniques, we generalise exterior Whitehead products to homotopy classes of maps between polyhedral products, and study relations among them.

We begin by describing the exterior Whitehead product as the map of polyhedral products
\[\{f_1,\dots,f_m\} \colon (\underline{D},\underline{S})^{\partial \Delta^{m-1}} \longrightarrow (\underline{Y},\underline{\ast})^{\partial \Delta^{m-1}}\]
induced by the maps of pairs $(D^{d_i},S^{d_i-1}) \longrightarrow (Y_i,\ast)$ representing $f_i \in \pi_{d_i}(Y_i)$. 
%We are motivated to ask what other relations exist. In general, homotopies between maps can be arbitrarily complicated, and therefore difficult to track algebraically. This has lead to there being a lack of algebraic models available. In order to address these complications, we realise the spaces in question as polyhedral products, following the work of Abramyan-Panov \cite{AbramyanPanov}, and Zhuravleva \cite{Zhur21} \todo{Move these later to general maps when algebraic methods don't work}. By combining geometric and combinatorial techniques, we are able to keep track of complex homotopy-theoretic phenomena.
Recognising that the boundary $\partial \Delta^{m-1}$ of a simplex governs the existence of the exterior Whitehead product, we define the higher Whitehead map $h_w(f_1,\dots,f_m) \in [(\underline{CX},\underline{X})^{\partial \Delta^{m-1}},(\underline{Y},\underline{\ast})^{\partial \Delta^{m-1}}]$ associated to maps of pairs $f_i \colon (CX_i,X_i) \longrightarrow (Y_i,\ast)$, $i=1,\dots,m$, to be the map of polyhedral products
\[
h_w(f_1,\dots,f_m) \colon (\underline{CX},\underline{X})^{\partial \Delta^{m-1}} \longrightarrow (\underline{Y},\underline{\ast})^{\partial \Delta^{m-1}}
\]
induced by the maps $f_i$. 

Higher Whitehead maps of spherical maps have been studied by Abramyan and Panov~\cite{AbramyanPanov}, and Zhuravleva~\cite{Zhur21}, using algebraic techniques such as the Hurewicz map and Adams--Hilton models. In the non-spherical case, these algebraic techniques do not generalise. We instead study higher Whitehead maps $h_w(f_1,\dots,f_m)$ and relations between them geometrically, by generalising the relative homotopy group approach of Nakaoka and Toda~\cite{Na-To54}. 

We generalise Hardie’s identity to 
maps on suspension spaces. We rewrite the maps in relation~\eqref{eq:HardieIdentity}, using polyhedral products and higher Whitehead maps, as 
\[
h_w(h_w(f_1,\dots,f_{i-1},f_{i+1},\dots,f_m),f_i) \colon S^{d_1 + \cdots + d_m - 2} \longrightarrow (\underline{Y},\underline{\ast})^{\sk^{m-3} \Delta^{m-1}}
\]
and say that the skeleton $\sk^{m-3} \Delta^{m-1}$ is the simplicial complex which carries Hardie's identity. By considering maps $f_i \colon \Sigma X_i \longrightarrow Y_i$, we study relations among higher Whitehead maps of the form
\[
h_w(h_w(f_1,\dots,f_{i-1},f_{i+1},\dots,f_m),f_i) \colon \Sigma^{m-2} X_1 \wedge \cdots \wedge X_m \longrightarrow (\underline{Y},\underline{\ast})^{\sk^{m-3} \Delta^{m-1}}.
\]
We prove, as Theorem~\ref{thm:maintheorem}, a generalisation of relation~\eqref{eq:HardieIdentity} to these maps by analysing the homotopy-theoretic properties of the polyhedral product in terms of the underlying combinatorics, and developing the notion of relative higher Whitehead maps as an analogue to the relative Whitehead product.

%We further study for which simplicial complexes $\cK$ the polyhedral product $(\underline{Y},\underline{\ast})^{\cK}$ supports a relation. . To obtain new relations, we propagate the combinatorial structure of the skeleton $\sk^{m-3} \Delta^{m-1}$ to more complex combinatorial objects.
%\todo[inline]{Via polyhedral product we bring in control by the combinatorics. Comment on the complexity of combinatorics allows us to increase the complexity of the forms of the maps/relations. Paragraph to paragraph emphasise that the combinatorial complexity is still increasing}
%\todo[inline]{Then go in and unpack this paragraph, case by case, through substitution, folding, etc.} 
%\todo[inline]{Motivate substitution, folding, relations, by saying that the combinatorics captures the `base case', eg. fat wedge, Hardie relation (skeleta), etc.)}
%We increase complexity in various ways to obtain different types of relations.  
Hardie's identity~\eqref{eq:HardieIdentity} is the only known relation involving higher Whitehead maps in homotopy groups. The summands are exterior Whitehead products of a higher Whitehead map and an arbitrary map. We further generalise relation~\eqref{eq:HardieIdentity} to new, more general relations between higher Whitehead maps, by viewing it within the combinatorial framework of the polyhedral join product. We first observe that the skeleton $\sk^{m-3} \Delta^{m-1}$ is the composition complex $\sk^{m-3} \Delta^{m-1} (\emptyset,\dots,\emptyset)$, see Definition~\ref{def:PolyhedralJoin}, and that the vertex set $[m] = \{1,\dots,m\}$ can be viewed as an $m$-partition of $[m]$ into singletons. To obtain new relations, we consider an $m$-partition $\Pi = \{P_1,\dots,P_m\}$ of a vertex set and replace the skeleton $\sk^{m-3}\Delta^{m-1}=\sk^{m-3} \Delta^{m-1} (\emptyset,\dots,\emptyset)$ with the complex obtained by exchanging each $\emptyset = \partial \Delta^0$ for $\partial \Delta[P_i]$. We define the identity complex $\cK_{\Pi}$ as the resulting complex $\sk^{m-3} \Delta^{m-1} (\partial \Delta[P_1],\dots,\partial \Delta[P_m])$ which carries a new identity.

The main result of this paper, given as Theorem~\ref{thm:maintheorem}, is that for every partition $\Pi$ and associated identity complex $\cK_{\Pi}$, there is a relation among higher Whitehead maps. Suppose that each $X_i$ is a suspension space. Then
\begin{equation} \label{eq:mainidentityinintro}
\sum_{i=1}^k h_w^{\cK_{\Pi}} \left( h_w \left( f_{j_1},\dots,f_{j_{q_i}} \right),f_{i_1},\dots,f_{i_{p_i}} \right) \circ \sigma_i = 0
\end{equation}
in $\left[ \Sigma^{m-2} X_1 \wedge \cdots \wedge X_m, (\underline{Y},\underline{\ast})^{\cK_{\Pi}}\right]$, where $P_i = \{i_1,\dots,i_{p_i}\}$, $\bigcup_{j=1}^m P_j \setminus P_i = \{j_1,\dots,j_{q_i}\}$ and
\[
\sigma_i \colon \Sigma^{m-2} X_1 \wedge \cdots \wedge X_m \longrightarrow \Sigma^{p_i} (\Sigma^{q_i-2} X_{j_1} \wedge \cdots \wedge X_{j_{q_i}}) \wedge X_{i_1} \wedge \cdots \wedge X_{i_{p_i}}
\]
is induced by the coordinate permutation
\[
X_1 \times \cdots \times X_m \longrightarrow X_{j_1} \times \cdots \times X_{j_{q_i}} \times X_{i_1} \times \cdots \times X_{i_{p_i}}.
\]

For a given identity complex $\cK_{\Pi}$, we construct a new family of simplicial complexes $\cK$ for which $(\underline{Y},\underline{\ast})^{\cK}$ has new relations. To do so, we propagate the structure of the identity complex $\cK_{\Pi}$ to more complex combinatorial objects by using the operation of  simplicial substitution. For a simplicial complex $\cK$ on $[m]$ and simplicial complexes $\cS_1,\dots,\cS_m$, we denote the substitution complex of $\cS_1,\dots,\cS_m$ into $\cK$ by $\cK \langle \cS_1,\dots,\cS_m \rangle$, see Definition~\ref{def:PolyhedralJoin}. In Corollary~\ref{cor:MainThmSubstituted}, we propagate the identity complex $\cK_{\Pi}$ to $\cK_{\Pi} \langle \cS_1,\dots,\cS_m \rangle$ which supports new relations.

%\todo[inline]{Introduce new structures, eg. torsion, nested brackets of different arity}

%\todo[inline]{Motivate folding, blah blah We enhance the folding operation by combining it with simplicial substitution in order to obtain more complex relations with repeated factors. In particular Elizaveta and torsion.}

In homotopy groups, the Lie bracket given by the Whitehead product can have repeated factors; for example $[[f_1,f_2],f_1]$, in which $f_1$ is repeated. Thus far, a higher Whitehead map captured by a simplicial complex has distinct factors, as the maps therein correspond to different vertices of the simplicial complex. This motivates the question of whether higher Whitehead maps with repeated factors can also be encoded by a simplicial complex.

We introduce the combinatorial operation of folding. A fold of a simplicial complex is the image of the simplicial map induced by identifying two or more vertices. Starting with a relation with no repeated factors corresponding to the identity complex $\cK_{\Pi}$, we fold $\cK_{\Pi}$ into a new simplicial complex, which supports relations among higher Whitehead maps with repeated factors corresponding to the identified vertices under the fold. 

We further increase the family of simplicial complexes which detect relations among higher Whitehead maps with repeated factors by considering folds of $\cK_{\Pi} \langle \cS_1,\dots,\cS_m \rangle$ induced by folds of the simplicial complexes $\cS_i$. 

%As an application of our methods, we obtain an integral relation which rationalises to the rational relation shown by Zhuravleva~\cite{Zhur21}.

An $L_{\infty}$-algebra, also known as a homotopy Lie algebra, is a higher generalisation of a Lie algebra, in which higher, multilinear brackets satisfy a strong homotopy Jacobi identity. In polyhedral products, we realise higher $L_{\infty}$ brackets as higher Whitehead maps, and using identity complexes and their folds, we realise the strong homotopy Jacobi identity.

%\todo[inline]{To summarise, we integrally realise the commutators and relations in $L_{\infty}$ (check used language), perhaps mention integral, putting $L_{\infty}$ structure on polyhedral products?}

% \todo[inline]{There is a thing as $L_{\infty}$, (multilinearity, symmetric), Jacobi comes from $L_{\infty}$ when all higher brackets vanish. In a similar way we recover relations among higher brackets described by $L_{\infty}$ structures}

In the course of the paper we also study various homotopy-theoretic properties of higher Whitehead maps, such as conditions under which they are null-homotopic.
An application of this, combined with folding simplicial complexes, is the construction of $2$-torsion elements in the homotopy groups of a wide range of polyhedral products.

\section{The higher Whitehead map}

%\todo[inline,color=green]{Change $r_i$ to $q_i$}

Given maps $f_i \in \pi_{d_i}(Y_i)$ for $i=1,\dots,k$, $k \geqslant 2, d_i \geqslant 1$, Hardie~\cite{Hardie61} defined the exterior Whitehead product $\{f_1,\dots,f_k\} \in \pi_{d_1 + \cdots + d_k - 1}(FW(Y_1,\dots,Y_k))$ as the homotopy class of the composite
\begin{equation} \label{eq:Hardiehwmap}
    \{f_1,\dots,f_k\} \colon \overset{k}{\underset{i=1}{\biggast}} S^{d_i-1} \xlongrightarrow{\rho} FW(S^{d_1},\dots,S^{d_k}) \xlongrightarrow{FW(f_1,\dots,f_k)} FW(Y_1,\dots,Y_k)
\end{equation}
where
\begin{equation} \label{eq:restofprodofquot}
\rho \colon \overset{k}{\underset{i=1}{\biggast}} S^{d_i-1} = \bigcup_{i=1}^k D^{d_1} \times \cdots \times S^{d_i-1} \times \cdots \times D^{d_k} \rightarrow FW(S^{d_1},\dots,S^{d_k})
\end{equation} 
is the restriction of the product of the quotient maps sending $\partial D^{d_j} = S^{d_j-1}$ to a point for $j = 1,\dots,k$.

%\todo[inline,color=green]{George has noted a new version of the following sentence, to add in }

By expressing the exterior Whitehead product as a map of polyhedral products, we study its generalisation to maps $f_i \in [\SX_i,Y_i]$. We start by defining polyhedral products.

\begin{definition}
\label{def:polyprod}
Let $\cK$ be a simplicial complex on the vertex set $[m]$, and let $(\underline{X},\underline{A})=\{(X_i,A_i)\}_{i=1}^m$ be an $m$-tuple of $CW$-pairs. The \textit{polyhedral product} $(\underline{X},\underline{A})^{\cK}$ is defined as
\[
(\underline{X},\underline{A})^{\cK} = \bigcup_{\sigma \in \cK} (\underline{X},\underline{A})^\sigma \subseteq \prod_{i=1}^m X_i, \quad \text{ 
where }(\underline{X},\underline{A})^\sigma = \prod_{i=1}^m Y_i, \quad  Y_i = \begin{cases} X_i & \text{for }i \in \sigma \\ A_i & \text{for }i \notin \sigma. \end{cases}
\]
\end{definition}

The polyhedral product is a covariant functor with respect to both inclusions of simplicial complexes and continuous maps of $m$-tuples of $CW$-pairs. A simplicial inclusion $\cL \longrightarrow \cK$ induces an inclusion of polyhedral products
\[
(\underline{X},\underline{A})^{\cL} \longrightarrow (\underline{X},\underline{A})^{\cK}
\]
and maps of $CW$-pairs $(f_i,g_i)\colon (X_i,A_i) \longrightarrow (Y_i,B_i)$ induce a map of polyhedral products
\[
(\underline{f},\underline{g})^{\cK} \colon (\underline{X},\underline{A})^{\cK} \longrightarrow (\underline{Y},\underline{B})^{\cK}.
\]
In the case that $B_i = *$ for all $i \in [m]$, we suppress notation and denote this induced map by $(\underline{f})^{\cK}$.

The fat-wedge and the join arise as special cases of polyhedral products. If $A_i = \ast$ for $i=1,\dots,m$, then 
\begin{equation} \label{eq:fatwedge}
(\underline{X},\underline{\ast})^{\partial \Delta^{m-1}} = \bigcup_{i=1}^m (X_1 \times \cdots \times X_{i-1} \times \ast \times X_{i+1} \times \cdots \times X_m) = FW(X_1,\dots,X_m)
\end{equation}
and if $X_i = CA_i$ for $i=1,\dots,m$, then
\begin{equation} \label{eq:join}
(\underline{CA},\underline{A})^{\partial \Delta^{m-1}} = \bigcup_{i=1}^m CA_1 \times \cdots \times CA_{i-1} \times A_i \times CA_{i+1} \times \cdots \times CA_m = \joinspc A_i.
\end{equation}

\goodbreak

Polyhedral products are the colimits of diagrams over the face category of a simplicial complex containing the empty set. Whenever all maps in such a diagram are cofibrant inclusions, the colimit is homotopy equivalent to the homotopy colimit (see \cite[Proposition~8.1.1]{buchstaber2014toric}), and therefore 
\begin{equation}\label{eq:polyprodhtpycolim}(\underline{X},\underline{A})^{\cK} \simeq \hocolim_{\sigma \in \cK}(\underline{X},\underline{A})^{\sigma}.
\end{equation} 
This enables us to study the homotopy-theoretic properties of polyhedral products combinatorially.

\begin{comment}
We make use of the following technical proposition in studying the homotopy-theoretic properties of polyhedral products.

\begin{proposition}[{\cite[Proposition~8.1.1]{buchstaber2014toric}}] \label{prop:colimequalshocolim}

If for all $i \in [m]$, the inclusion $A_i \hooklongrightarrow X_i$ is a cofibration, then there is a homotopy equivalence
\begin{equation} \label{eq:colimequalshocolim}
    \hocolim_{\sigma \in \cK}(\underline{X},\underline{A})^{\sigma} \simeq \colim_{\sigma \in \cK}(\underline{X},\underline{A})^{\sigma}.
\end{equation}
\end{proposition}
\end{comment}

Using the notion of polyhedral products, we rewrite the definition of the exterior Whitehead product. Let $(\underline{D},\underline{S}) = \{(D^{d_i},S^{d_i-1})\}_{i=1}^k$ and  $(\underline{Y},\underline{\ast}) = \{(Y_i,\ast)\}_{i=1}^k$ be $k$-tuples of $CW$-pairs and for $i=1,\dots,k$, identify $f_i \colon S^{d_i} \longrightarrow Y_i$ with the map of pairs $(D^{d_i},S^{d_i-1}) \longrightarrow (Y_i,\ast)$. Then the exterior Whitehead product~\eqref{eq:Hardiehwmap} is the map of polyhedral products
\begin{equation} \label{eq:hwmap}
\{f_1,\dots,f_k\} \colon (\underline{D},\underline{S})^{\partial \Delta^{k-1}} \longrightarrow (\underline{Y},\underline{\ast})^{\partial \Delta^{k-1}}
\end{equation}
induced by the maps of pairs $(D^{d_i},S^{d_i-1}) \longrightarrow (Y_i,\ast)$ for $i=1,\dots,k$. We generalise this construction to maps $\Sigma X_i \longrightarrow Y_i$, which can be identified with maps of pairs $(CX_i,X_i) \longrightarrow (Y_i,\ast)$.

\goodbreak 
\begin{definition} \label{def:higherwhiteheadmap}
Let $m \geqslant 2$, and let $(\underline{CX},\underline{X}) = \{(CX,X)\}_{i=1}^m$ and $(\underline{Y},\underline{\ast}) = \{(Y_i,\ast)\}_{i=1}^m$ be $m$-tuples of $CW$-pairs. Let $f_i \colon (CX_i,X_i) \longrightarrow (Y_i,*)$ be maps of pairs for $i = 1,\dots,m$.  The \textit{higher Whitehead map} $h_w(f_1,\dots,f_m)$ of the maps $f_1,\dots,f_m$ is the map of polyhedral products
\begin{equation}
    h_w(f_1,\dots,f_m) \colon (\underline{CX},\underline{X})^{\partial \Delta^{m-1}} \longrightarrow (\underline{Y},\underline{\ast})^{\partial \Delta^{m-1}}
\end{equation}
induced by the maps $ f_i $.
\end{definition}

Since the polyhedral product is a homotopy colimit (see \eqref{eq:polyprodhtpycolim}), the homotopy class of $h_w(f_1,\dots,f_m)$ depends only on the homotopy classes of the maps $f_1, \dots, f_m$.

%\todo[inline,color=green]{George has noted a new version of the sentence below.}

The higher Whitehead map exhibits multilinearity, naturality, and symmetry with respect to the maps $f_i$.
\begin{proposition}
 Let $f_i \colon (CX_i,X_i) \longrightarrow (Y_i,*)$ be maps of pairs for $i = 1,\dots,m$, and for some $j \in \{1,\dots,m\}$, let $f_j' \colon (CX_j,X_j) \longrightarrow (Y_j,*)$. If $X_j$ is a suspension,
\begin{equation} \label{eq:multilinearityofhwp}
h_w(f_1,\dots,f_j + f_j',\dots,f_m) = h_w(f_1,\dots,f_j,\dots,f_m) + h_w(f_1,\dots,f_j',\dots,f_m).
\end{equation}
\end{proposition}

%\[\begin{tikzcd}
%(\mathbf{CX},\mathbf{X})^{\partial \Delta^{m-1}} \arrow[r,"\mathbf{f}"] \arrow[d,"F"] & (\mathbf{CY},\mathbf{Y})^{\partial \Delta^{m-1}} \arrow[d,"F"]\\ \bigwedge_{i=1}^{m-1} S^1 \wedge \bigwedge_{i=1}^m X_i \arrow[r,"\mathbf{f}"] & \bigwedge_{i=1}^{m-1} S^1 \wedge \bigwedge_{i=1}^m Y_i 
%\end{tikzcd}\]

%\todo[inline]{In the above diagram, let $\mathbf{f}$ be a map induced by $f_i\colon X_i \rightarrow Y_i$, and $F$ is the homotopy equivalence going from one version of the join (poly prod) to the other (smash prod). Normally we'd have to prove that this diagram commutes (up to homotopy). However Porter already proved it in the appendix cited article. }

\begin{proof}
    By~\cite[Corollary~p. 135]{Porter1965HigherOW}, there is a homotopy equivalence \[(\underline{CX},\underline{X})^{\partial \Delta^{m-1}} = \bigcup_{i=1}^m CX_1 \times \cdots \times X_i \times \cdots \times CX_m \simeq \bigwedge_{i=1}^{m-1} S^1 \wedge \bigwedge_{i=1}^m X_i,\] which is natural with respect to continuous maps on $(\underline{CX},\underline{X})$ induced by maps on $\underline{X}$. Then, given the homeomorphism
\[\left(\bigwedge_{i=1}^{m-1} S^1 \right) \wedge X_1 \wedge \cdots \wedge \left( X_j \vee X_j \right) \wedge \cdots \wedge X_m \cong \left(\bigwedge_{i=1}^{m-1} S^1 \wedge \bigwedge_{i=1}^m X_i  \right) \vee \left(\bigwedge_{i=1}^{m-1} S^1  \wedge \bigwedge_{i=1}^m X_i \right)\] 
and up to homotopy unique co-multiplication induced by any of the factors, the result follows.
\end{proof}

\begin{proposition}
\label{prop:naturality}
For $i = 1,\dots,m$, let $f_i \colon (CX_i,X_i) \longrightarrow (Y_i,*)$ be maps of pairs, and let $g_i \colon Y_i \longrightarrow Z_i$ and $l_i \colon W_i \longrightarrow X_i$ be maps. Denote by $g \colon (\underline{Y},\underline{*})^{\partial \Delta^{m-1}} \longrightarrow (\underline{Z},\underline{*})^{\partial \Delta^{m-1}}$ the map of polyhedral products induced by the maps $g_i$. Then
\begin{equation} \label{eq:naturality1}
h_w(g_1 \circ f_1, \dots, g_m \circ f_m) = g \circ h_w(f_1,\dots,f_m)
\end{equation}
and
\begin{equation} \label{eq:naturality2}
h_w(f_1 \circ \Sigma l_1,\dots,f_m \circ \Sigma l_m) = h_w(f_1,\dots,f_m) \circ \biggast_{i=1}^m l_i
\end{equation}
where $\Sigma l_i \colon (CW_i,W_i) \longrightarrow (CX_i,X_i)$ denotes the map of pairs induced by $l_i$.
\end{proposition}

\begin{proof}
Both follow by the functoriality of the polyhedral product with respect to continuous maps of pairs. %Identity \eqref{eq:naturality2} follows from the naturality of the map $\rho = (\underline{\rho})^{\partial \Delta^{m-1}} \colon (\underline{CX},\underline{X})^{\partial \Delta^{m-1}}\longrightarrow (\underline{\Sigma X},\underline{*})^{\partial \Delta^{m-1}}.$
\end{proof}

In order to prove the symmetry property of the higher Whitehead map, we express it as the composite
\begin{equation} \label{eq:hwmapdecomp}
(\underline{CX},\underline{X})^{\partial \Delta^{m-1}} \xlongrightarrow{\rho}(\underline{\SX},\underline{*})^{\partial \Delta^{m-1}} \xlongrightarrow{(\underline{f})^{\partial \Delta^{m-1}}}(\underline{Y},\underline{*})^{\partial \Delta^{m-1}}
\end{equation}
where $\rho$ denotes the map induced by the maps of pairs $(CX_i,X_i) \longrightarrow (\SX_i,*)$.

Let $\sigma \colon  (1,\dots,m)\mapsto (i_1,\dots,i_m)$ denote a permutation. For a space $Z$, consider maps $\iota \colon FW(Y_1,\dots,Y_m) \longrightarrow Z$ and $\iota' \colon FW(Y_{i_1},\dots,Y_{i_m}) \longrightarrow Z$ such that $\iota'(Y_{\sigma(i)}) = \iota(Y_i)$ for $i=1,\dots,m$. Denote by
\[h_w^Z (f_1,\dots,f_m) = \iota \circ h_w(f_1,\dots,f_m) \colon (\underline{CX},\underline{X})^{\partial \Delta^{m-1}} \longrightarrow Z\] and by \[ h_w^Z (f_{i_1},\dots,f_{i_m}) = \iota' \circ h_w(\sigma(f_1,\dots,f_m)) \colon (\underline{CX},\underline{X})_{\sigma}^{\partial \Delta^{m-1}} \longrightarrow Z \]
where $(\underline{CX},\underline{X})_{\sigma}$ denotes the permuted tuple. The permutation $\sigma$ induces the permutation on the products $\prod CX_i$ and $\prod Y_i$, which in turn induce the permutation on the join $ X_1 \ast \cdots \ast X_m$ and on the fat wedge $FW(Y_1,\dots,Y_m)$, respectively.

\begin{proposition} \label{prop:hwSymmetry}
For $i = 1,\dots,m$, let $f_i \colon (C X_i,X_i) \longrightarrow (Y_i,*)$ be maps. Then
\begin{equation} \label{eq:GenCommRule2}
   h_w^Z(f_1,\dots,f_m)= h_w^Z(\sigma(f_1,\dots,f_m))\circ\sigma 
\end{equation}
and furthermore if $X_i = S^{p_i-1}$ for $i=1,\dots,m$, then
\begin{equation} \label{eq:SpherCommRule2}
   h_w^Z(f_1,\dots,f_m)= \epsilon(\sigma) h_w^Z(\sigma(f_1,\dots,f_m))
\end{equation}
where $\epsilon(\sigma)$ is the Koszul sign of $\sigma$, i.e. the product of $(-1)^{p_i p_j}$ for every transposition $(i,j)$ of $\sigma$.
\end{proposition}

\begin{proof}
Consider the following diagram
\begin{equation*}
\begin{tikzcd}[column sep = 7em]  (\underline{CX},\underline{X})^{\partial \Delta^{m-1}} \arrow[r,"{h_w}(\sigma(f_1{,\dots,}f_m)) \circ \sigma"] \arrow[d, equal] & FW(Y_{i_1},\dots,Y_{i_m}) \arrow[dr,"\iota'"]\\
    (\underline{CX},\underline{X})^{\partial \Delta^{m-1}} \arrow[r,"{h_w}(f_1{,\dots,}f_m)"] & FW(Y_1,\dots,Y_m) \arrow[u,"\sigma"] \arrow[r,"\iota"] & Z.
\end{tikzcd}
\end{equation*}
The square commutes since the higher Whitehead map is the restriction of the product map $\prod f_i$, and therefore acts coordinate-wise. The triangle commutes by the definition of the maps $\iota$ and $\iota'$. Therefore the diagram commutes, and identity~\eqref{eq:GenCommRule2} follows. Identity~\eqref{eq:SpherCommRule2} then follows from~\eqref{eq:GenCommRule2} since in that case, the map $\sigma$ is a map $S^{p_1+\cdots+p_m-1} \longrightarrow S^{p_1+\cdots+p_m-1}$ of degree $\epsilon(\sigma)$.
\end{proof}

The higher Whitehead map is an element of the higher Whitehead product of Porter~\cite[Definition~1.4]{Porter1965HigherOW}, which is defined as follows. Consider spaces $X_1,\dots,X_m, Z$. Given maps $f_i\colon \SX_i \longrightarrow Z$ for $i = 1, \dots, m$, denote by
\[
\omega (f_1,\ldots, f_m) = \{\phi\colon FW(\SX_1,\dots,\SX_m) \longrightarrow Z \; | \; \phi \big|_{\SX_i} \simeq f_i \; \text{for } i = 1, \dots, m\}
\]
the set of extensions up to homotopy of $\bigvee f_i\colon \bigvee \SX_i \longrightarrow Z$ to $FW(\SX_1,\dots,\SX_m)$. The \emph{$k$-ary higher Whitehead product} $[f_1,\dots,f_m]$ of the maps $f_i$ is the set of homotopy classes
\begin{equation} \label{eq:hwproduct}
[f_1,...,f_m] = \left\{[\phi \circ \rho] \; | \; \phi \in \omega(f_1,\dots,f_m) \right\} \subseteq \left[\joinspc X_i, Z\right]
\end{equation}
where $\rho$ is defined as in \eqref{eq:hwmapdecomp}, noting that by \eqref{eq:join}, $\joinspc X_i = (\underline{CX},\underline{X})^{\partial \Delta^{m-1}}$, and by \eqref{eq:fatwedge}, $FW(\Sigma X_1,\ldots, \Sigma X_m) = (\underline{\Sigma X},\underline{\ast})^{\partial \Delta^{m-1}}$.
%\todo[color=green, inline]{GS: If I understood the reviewer comment correctly, then it potentially highlights confusion in how we defined $\rho$, although I don't see this confusion personally. MS: Fixed. Reviewer wanted clarity on what $\rho$ is here. Have defined $\rho$ just after equation 12, referenced that definition here, and made notation consistent.}

\begin{proposition} \label{prop:hwisHWPBasic}
For $i = 1,\dots,m$, let $f_i \colon (C X_i,X_i) \longrightarrow (Y_i,*)$ be maps. Denote by $\iota_j \colon Y_j \longrightarrow FW(Y_1,\dots,Y_m)$ the inclusion of $Y_j$ into the $j$th coordinate. Then 
\[
h_w(f_1,\dots,f_m) \in [\iota_1 \circ f_1,\dots,\iota_m \circ f_m ] \subseteq \left[\joinspc X_i, FW(Y_1,\dots,Y_m)\right].
\]
\end{proposition}
\begin{proof}
  Notice that the higher Whitehead map is the composite \[\joinspc X_i=(\underline{CX},\underline{X})^{\partial \Delta^{m-1}} \xlongrightarrow{\rho} (\underline{\SX},\underline{*})^{\partial \Delta^{m-1}} \xlongrightarrow{(\underline{f})^{\partial \Delta^{m-1
}}} (\underline{Y},\underline{*})^{\partial \Delta^{m-1}}=FW(Y_1,\ldots,Y_m).\]
For the higher Whitehead product, all the maps need to have the same codomain. Hence we include the images of the maps $f_i$ into $FW(Y_1,\ldots,Y_m)$ by coordinate inclusions. This proves the claim.
\end{proof}

%\todo[color=green, inline]{GS: For the reviewer, $h_w(f_1,\dots,f_m)$ lands in $FW(Y_1,\dots,Y_m)$, so $Z = FW(Y_1,\dots,Y_m)$. So we just need to add a line ``taking $Z$ to be ...'' in the above proof somewhere. MS: Fixed. We previously had $(\underline{f})^{\partial \Delta^{m-1}}\in \omega(f_1,\dots,f_m)$, which was wrong, because the $f_i$ had different codomains.}

The higher Whitehead product is trivial when it contains a null-homotopic map. This motivates an investigation of the conditions under which the higher Whitehead map, an element of the higher Whitehead product, is null-homotopic. The following result expresses this property in terms of the adjoints of the maps $f_i$.

\begin{proposition}
    \label{prop:hwtrivialifflift}
    For $i=1,\dots,m$ let $f_i:(CX_i,X_i) \longrightarrow (Y_i, \ast)$ be a map of pairs and let $\hat{f}_i \colon X_i \longrightarrow \Omega Y_i$ be its adjoint. Then the higher Whitehead map $h_w(f_1,\dots,f_m)$ is null-homotopic if and only if the join $\hat{f}_1 \ast \cdots \ast \hat{f}_m$ is null-homotopic.
\end{proposition}

\begin{proof}
    In the following diagram
\begin{equation} \label{eq:MainTrivialityDiagram}
    \begin{tikzcd}
    \prod_{i=1}^m \Omega Y_i \ar[r,"{\ast}"] & \Omega Y_1 \ast \cdots \ast \Omega Y_m \ar[r] & FW(Y_1,\dots,Y_m) \ar[r,"\iota"] & \prod_{i=1}^m Y_i \\
    & & X_1 \ast \cdots \ast X_m \ar[r] \ar[u,"{h_w(f_1,\dots,f_m)}"'] \ar[ul,"{\hat{f}_1 \ast \cdots \ast \hat{f}_m}",dashed] & \prod_{i=1}^m CX_i \ar[u]
    \end{tikzcd}
\end{equation}
the top row is a homotopy fibration sequence, and the right square commutes up to homotopy. The composite $\iota \circ h_w(f_1,\dots,f_m)$ is null-homotopic, and a lift to the homotopy fibre of $\iota$ is given by $\hat{f}_1 \ast \cdots \ast \hat{f}_m$. Therefore if $\hat{f}_1 \ast \cdots \ast \hat{f}_m$ is null-homotopic, so is $h_w(f_1,\dots,f_m)$. On the other hand, if $h_w(f_1,\dots,f_m)$ is null-homotopic, then there is a further lift $X_1 \ast \cdots \ast X_m \longrightarrow \Omega Y_1 \times \cdots \times \Omega Y_m$ of the map $\hat{f}_1 \ast \cdots \ast \hat{f}_m$. Since $\Omega Y_1 \times \cdots \times \Omega Y_m  \longrightarrow  \Omega Y_1 \ast \cdots  \ast \Omega Y_m$ is null-homotopic, it follows that $\hat{f}_1 \ast \cdots \ast \hat{f}_m$ is itself null-homotopic.
\end{proof}

Proposition~\ref{prop:hwtrivialifflift} shows that the higher Whitehead map $h_w(f_1,\dots,f_m)$ is null-homotopic if any factor $f_i$ is null-homotopic. In general, fully characterising the triviality of the join $\hat{f}_1 \ast \cdots \ast \hat{f}_m$, and hence the higher Whitehead map, is difficult. The following example demonstrates that for non null-homotopic maps $\hat{f}_i$, the join can be trivial, due to topological properties of the join.

\begin{example} \label{ex:hwAdjointTriv}
Let $f_1 \colon \Sigma M(\mathbb{Z}_2,1) \longrightarrow Y_1$ and $f_2 \colon \Sigma M(\mathbb{Z}_3,1) \longrightarrow Y_2$ be non null-homotopic, where $M(G,n)$ denotes the Moore space with reduced homology $G$ concentrated in degree $n$. The space $M(\mathbb{Z}_2,1) \wedge M(\mathbb{Z}_3,1)$ is simply-connected, and by the K{\"u}nneth theorem has trivial homology in all positive degrees, and therefore is contractible. Thus the join $\hat{f}_1 \ast \hat{f}_2$ has contractible domain, and hence the higher Whitehead map
\[
h_w(f_1,f_2) \colon \Sigma M(\mathbb{Z}_2,1) \wedge M(\mathbb{Z}_3,1) \longrightarrow Y_1 \vee Y_2 
\]
is null-homotopic.
\end{example} 

Non-triviality of the join $\hat{f}_1 \ast \cdots \ast \hat{f}_m$ can be detected in certain cases, in homology. If each map $\hat{f_i}$ induces a non-trivial map in homology in degree $n_i$ over a common field $\Bbbk$, then by the K{\"u}nneth theorem the smash product $\hat{f}_1 \wedge \cdots \wedge \hat{f}_m$ induces a non-trivial map in homology in degree $n_1 + \cdots + n_m$ over $\Bbbk$. In such a case, $\hat{f}_1 \ast \cdots \ast \hat{f}_m$ is not null-homotopic. 

\begin{example} \label{ex:HopfHigherWhitehead}
    Consider the Hopf map $\eta \in \pi_3(S^2)$ and let $f_2,\dots,f_k$ be identity maps $S^2 \longrightarrow S^2$. For any field $\Bbbk$, the adjoints $\hat{f}_i \colon S^1 \longrightarrow \Omega S^2$ induce non-trivial maps in degree $1$ homology for $i=2,\dots,k$, and the adjoint $\hat{\eta}: S^2 \longrightarrow \Omega S^2$ induces a non-trivial map $H_2(S^2; \Bbbk) \longrightarrow H_2(\Omega S^2;\Bbbk)$. Therefore, the map induced by $\hat{f}_1 \wedge \cdots \wedge \hat{f}_k$ is non-trivial in degree $k+1$. Therefore applying Proposition~\ref{prop:hwtrivialifflift}, the higher Whitehead map $h_w(\eta, f_2,\dots,f_m)$ is not null-homotopic.
\end{example}

Examples~\ref{ex:hwAdjointTriv} and \ref{ex:HopfHigherWhitehead} show that topological properties of the factors $f_i$ can determine the triviality of the higher Whitehead map. We now analyse combinatorial elements of the construction which imply triviality. By~\cite[Theorem~2.3]{Porter1965HigherOW}, the mapping cone of $\rho \colon (\underline{CX},\underline{X})^{\partial \Delta^{m-1}} \longrightarrow (\underline{\SX},\underline{*})^{\partial \Delta^{m-1}}$ is homotopy equivalent to the product $(\underline{\SX},\underline{*})^{ \Delta^{m-1}}$, so that the top row of the diagram 
\begin{equation} \label{eq:hwcofibdiag}
    \begin{tikzcd}
    (\underline{CX},\underline{X})^{\partial \Delta^{m-1}} \ar[r,"\rho"] \ar[dr,swap,"{h_w(f_1,\dots,f_m)}"] & (\underline{\SX},\underline{*})^{\partial \Delta^{m-1}} \ar[r] \ar[d]& (\underline{\SX},\underline{*})^{ \Delta^{m-1}} \ar[dl,dashed,"\phi"]\\
    & (\underline{Y},\underline{*})^{\partial \Delta^{m-1}}
    \end{tikzcd}
\end{equation}
 is a homotopy cofibration. Therefore, $h_w(f_1,\dots,f_m)$ is null-homotopic if and only if there exists an extension up to homotopy $\phi \colon (\underline{\SX},\underline{*})^{ \Delta^{m-1}} \longrightarrow (\underline{Y},\underline{*})^{\partial \Delta^{m-1}}$ of $(\underline{f})^{\partial \Delta^{m-1}}\colon (\underline{\SX},\underline{*})^{\partial \Delta^{m-1}} \longrightarrow (\underline{Y},\underline{*})^{\partial \Delta^{m-1}}$. 
 %\todo[color=green, inline]{GS: Reviewer wants us to link Porter's homotopy equivalence $T_1 \Sigma \cup_{W(i)} C\Sigma^{n-1} \wedge \longrightarrow T_0 \Sigma$ to first row of diagram being cofib. MS: Fixed. Added sentence stating that Porter's result is that the mapping cone is homotopy equiv to the product} 
 Notice that if we consider instead the composite \[\iota \circ h_w(f_1,\dots,f_m) \colon (\underline{CX},\underline{X})^{\partial \Delta^{m-1}} \xlongrightarrow{h_w(f_1,\dots,f_m)} (\underline{Y},\underline{*})^{\partial \Delta^{m-1}} \xlongrightarrow{\iota} (\underline{Y},\underline{*})^{\Delta^{m-1}}\]
where $\iota$ is the map induced by the inclusion $\partial \Delta^{m-1} \longrightarrow \Delta^{m-1}$, then an extension $\phi$ always exists, given by the product $(\underline{f})^{\Delta^{m-1}}$, and the higher Whitehead map is therefore trivial in the space $(\underline{Y},\underline{*})^{\Delta^{m-1}}$ corresponding to the full simplex $\Delta^{m-1}$. This observation motivates the question of whether there are other simplicial complexes $\cK$ such that the composite $\iota \circ h_w(f_1,\dots,f_m) \colon (\underline{CX},\underline{X})^{\partial \Delta} \longrightarrow (\underline{Y},\underline{\ast})^{\cK}$ is null-homotopic.

For a subset $J=\{j_1,\dots,j_n\} \subseteq [m]$, denote by $\Delta[J]$ and $\partial \Delta[J]$  the simplex on the vertex set $J$ and its boundary complex, respectively. For $\cK$ containing $ \partial \Delta[J]$, define
\begin{equation}
\label{eq:basichw^K}
h_w^{\cK}(f_{j_1},\dots,f_{j_n}) \colon (\underline{CX},\underline{X})^{\partial\Delta[J]}\longrightarrow (\underline{Y},\underline{\ast})^{\cK}
\end{equation} 
to be the composite
\[(\underline{CX},\underline{X})^{\partial \Delta[J]}\xlongrightarrow{h_w(f_{j_1},\dots,f_{j_n})}(\underline{Y},\underline{*})^{\partial \Delta[J]} \xlongrightarrow{\iota}(\underline{Y},\underline{\ast})^{\cK}\] where the map $\iota$ is induced by the inclusion $\iota \colon \partial \Delta[J] \longrightarrow \cK$.

\begin{proposition}
\label{lem:basichwtriviality}
Let $\cK$ be a simplicial complex which contains $\partial \Delta[J]$. Then the map \[h_w^{\cK}(f_{j_1},\dots,f_{j_n}) \colon (\underline{CX},\underline{X})^{\partial \Delta[J]} \longrightarrow (\underline{Y},\underline{\ast})^{\cK}\] is null-homotopic if either $f_{j_i}$ is null-homotopic for some $i \in \{1,\dots,m\}$, or $\Delta[J] \subseteq \cK$. 
\end{proposition}

\begin{proof}

Consider diagram \eqref{eq:hwcofibdiag}. If $f_i$ is null-homotopic for some $i \in \{1,\dots,m\}$, then $\phi$ can be chosen to be $
f_1 \times \cdots \times f_{i-1} \times \ast \times f_{i+1} \times \cdots \times f_m$.
In this case, the diagram commutes up to homotopy. If $\Delta[J] \subseteq \cK$, then $\phi$ can be chosen to be $f_1 \times \cdots \times f_m$.
\end{proof}

\section{The nested higher Whitehead map}

\begin{comment}
A class of Whitehead and exterior Whitehead products of particular interest are nested products, which satisfy the Jacobi and Hardie's identities respectively. We consider the case of a nested higher Whitehead map, where for some $i$, the map $f_i$ is of the form \eqref{eq:basichw^K}. Nested higher Whitehead maps include nested Whitehead and exterior Whitehead products, and through the lens of the higher Whitehead map, we introduce a richer combinatorial structure into their study.
\end{comment}

Higher Whitehead maps are defined as maps of polyhedral products over the boundary of a simplex. In this section we study a class of higher Whitehead maps with richer combinatorial structure, by supposing that the maps $f_i$ themselves are higher Whitehead maps. Such maps appear, for example, in the Jacobi identity~\eqref{eq:JacobiIdentity} and Hardie's identity~\eqref{eq:HardieIdentity}.

\begin{definition} \label{def:nestedhwmap}

A \textit{nested higher Whitehead map} is a higher Whitehead map $h_w(f_1,\dots,f_m)$ where for at least one $i \in \{1,\dots,m\}$, the map $f_i$ is itself a higher Whitehead map.

\end{definition}

%For a nested higher Whitehead map, for all $i \in [m]$ we denote the codomain of $f_i$ by $\left(\underline{Y}_i,\underline{*} \right)^{\partial \Delta^{l_i-1}}$ where for any $f_i$ which is not a higher Whitehead map, $l_i =1 $ and $\partial \Delta^0 = \bullet$, so that $\left(\underline{Y}_i,\underline{*} \right)^{\partial \Delta^{l_i-1}} = Y_i$.

Consider a nested higher Whitehead map $h_w(f_1,\dots,f_m)$. For $i \in \{1,\dots,m\}$, if $f_i \colon \Sigma X_i \longrightarrow (\underline{Y},\underline{\ast})^{\partial \Delta^{l_i-1}}$ is a higher Whitehead map, define $\partial \overline{\Delta}^{l_i-1} = \partial \Delta^{l_i-1}$, and otherwise let $\partial \overline{\Delta}^{l_i-1} = \bullet$, the simplex on one vertex. The codomain of the nested higher Whitehead map $h_w(f_1,\dots,f_m)$ is then the polyhedral product
\begin{equation} \label{eq:nestedhwcodom}
\left( \underline{\left(\underline{Y}_i,\underline{*} \right)}^{\partial \overline{\Delta}^{l_i-1}},\underline{*} \right)^{\partial \Delta^{m-1}}.
\end{equation}

%Let $m \geqslant 2$. For $i = 1,\dots,m$, let $l_i \geqslant 1$, and suppose that for $j = 1,\dots,l_i$, $f_{i_j} \colon \SX_{i_j} \longrightarrow Y_{i_j}$ are maps. Denote 
%\begin{equation} \label{eq:innerhwmapfornested}
%f_i = h_w(f_{i_1},\dots,f_{i_{l_i}}) \colon \left(\underline{CX}_i,\underline{X}_i \right)^{\partial \Delta^{l_i}-1} \longrightarrow \left(\underline{Y}_i,\underline{*}\right)^{\partial \Delta^{l_i}-1}
%\end{equation}
%where if $l_i =1$, then we set $f_i = h_w(f_{i_1}) = f_{i_1}$. 

%\todo[inline]{Put all of nested in Definition, and adjust references throughout}

%The nested higher Whitehead map of $f_1,\dots,f_m$ is defined as the higher Whitehead map
%\begin{equation} \label{eq:nestedhwmap}
%h_w(f_1,\dots,f_m) \colon \left(\underline{CX},\underline{X}\right)^{\partial \Delta^{m-1}} \longrightarrow \left(\underline{Y},\underline{*}\right)^{\partial \Delta^{m-1}}
%\end{equation}
%where $\underline{Y} = \left\{\left(\underline{Y}_i,\underline{*} \right)^{\partial \Delta^{l_i-1}}\right\}_{i=1}^m$, as a nested higher Whitehead map.\

 %   A nested higher Whitehead map is a higher Whitehead map $h_w(f_1,\dots,f_m)$ where at least one of the maps $f_i$ is a higher Whitehead map.
 %   \[
%    (\underline{Y},\underline{\ast}) = ((Y_1,\ast),\dots,(Y_{i-1},\ast),(\underline{Y}_i,\underline{\ast})^{\partial \Delta^{l_i-1}},(Y_{i+1},\ast),\dots,(Y_m,\ast))
%    \]
%    where
 %   \[(\underline{Y}_i,\underline{\ast}) = ((Y_{i_1},\ast),\dots,(Y_{i_{l_i}},\ast))
%    \]
%\todo[inline]{The point is to align notation with the polyhedral join}

To express the rich combinatorial structure arising within this polyhedral product, we introduce the polyhedral join product, defined by Vidaurre \cite{Vidaurre}. 

A simplicial pair $(\cS,\cT)$ consists of simplicial complexes $\cS$ and $\cT$, both with vertex set $[l]$, such that $\cT$ is a subcomplex of $\cS$.

\begin{definition} \label{def:PolyhedralJoin}
Let $\cK$ be a simplicial complex on $[m]$, and let $(\cS_i,\cT_i)$ be a simplicial pair on $[l_i]$ for $i=1,\dots,m$. Let $(\underline{\cS},\underline{\cT})=\{(\cS_i,\cT_i)\}_{i=1}^m$ be an $m$-tuple of simplicial pairs. The \textit{polyhedral join product} is the simplicial complex on vertex set $[l_1] \sqcup \dots \sqcup [l_m]$ defined by
\begin{equation*} \label{eq:polyjoinprod}
(\underline{\cS},\underline{\cT})^{*\cK} = \bigcup_{\sigma \in \cK} (\underline{\cS},\underline{\cT})^{*\sigma} \subseteq \biggast_{i=1}^m \cS_i, \quad \text{ where }
(\underline{\cS},\underline{\cT})^{*\sigma} = \biggast_{i=1}^m \mathcal{Y}_i,  \quad \mathcal{Y}_i = \begin{cases} \cS_i & \text{for } i \in \sigma \\ \cT_i & \text{for } i \notin \sigma.\end{cases} \end{equation*}
\end{definition}

Special cases of polyhedral join products have been studied. If $\cT_i = \{ \emptyset \}$ for all $i$, then we denote by
\begin{equation} \label{eq:subsaspolyjoin}
\cK\langle \cS_1,...,\cS_m\rangle=(\underline{\cS},\{\underline{\emptyset}\})^{*\cK}
\end{equation} 
the \emph{substitution} of $\cS_1,...,\cS_m$ into $\cK$ (see~\cite{AbramyanPanov}). If $\cS_i = \Delta^{l_i-1}$ for all $i$, then we denote by \begin{equation} \label{eq:compaspolyjoin}
\cK(\cT_1,\dots,\cT_m)=(\underline{\Delta}^{l_i-1},\underline{\cT})^{*\cK}
\end{equation}
the \emph{composition} of $\cK$ with $\cT_1,...,\cT_m$ (see~\cite{anton2013composition}). In the case that $\cK=\Delta^{m-1}$ or ${\cK = \partial \Delta^{m-1}}$, we denote by
\begin{align} \label{eq:subandcompoverfullsimplex}
    \Delta \langle \cK_1,\dots,\cK_m \rangle=\Delta^{m-1} \langle \cK_1,\dots,\cK_m \rangle
    \end{align}
and
\begin{align} \label{eq:subandcompoverboundaryofsimplex}
    \partial \Delta\langle \cK_1,\dots,\cK_m \rangle=\partial \Delta^{m-1} \langle \cK_1,\dots,\cK_m\rangle
\end{align}
the substitution complexes over the full simplex and the boundary, respectively, and adopt the same abbreviation for composition complexes $\Delta \left( \cK_1,\dots,\cK_m \right)$ and $\partial \Delta \left( \cK_1,\dots,\cK_m \right)$.
%\begin{definition}
%    Let $\cK_1,\dots,\cK_m$ be simplicial complexes on vertex sets $[m_1],\dots,[m_n]$ respectively. The substitution $\Delta^{m-1} \langle \cK_1,\dots,\cK_m \rangle$ is denoted $\Delta \langle \cK_1,\dots,\cK_m \rangle$, while the substitution $\partial \Delta^{m-1} \langle \cK_1,\dots,\cK_m\rangle$ is denoted $\partial \Delta\langle \cK_1,\dots,\cK_m \rangle$. The corresponding abbreviations are used for the composition complexes $\Delta \left( \cK_1,\dots,\cK_m \right)$ and $\partial \Delta \left( \cK_1,\dots,\cK_m \right)$.
%\end{definition}

\begin{comment}denote the corresponding substitution complexes by $\Delta \langle \cS_1,\dots,\cS_m \rangle $ and $\partial \Delta \langle \cS_1,\dots,\cS_m \rangle$, respectively, and adopt the similar abbreviation for composition complexes. 
\end{comment} 
We observe that for any simplicial complex $\cK$, $\cK \langle \bullet, \dots, \bullet \rangle = \cK = \cK ( \circ, \dots, \circ)$, where $\bullet$ denotes the simplicial complex consisting of a single vertex, and $\circ$ denotes the empty complex on a single vertex.

%\todo[color=green, inline]{GS: Reviewer wanted $\Delta \langle \cS_1,\dots,\cS_m \rangle $ and $\partial \Delta \langle \cS_1,\dots,\cS_m \rangle$ to be in a definition, but we don't actually call these complexes a unique name.. Original text in a comment. MS: Put them inside an equation environment, and referred back to these when they're used later. This I believe addresses the reviewer's concerns on this one.}

In analogy with the polyhedral product, the polyhedral join product is covariantly functorial with respect to simplicial inclusions ${\cL \longrightarrow \cK}$, and simplicial maps of pairs of simplicial complexes.

Viadurre showed that polyhedral products over polyhedral join products decompose in the following way.

\begin{theorem}[{\cite[Theorem~2.9]{Vidaurre}}]
Let $(\underline{\cS},\underline{\cT})=\{(\cS_i,\cT_i)\}_{i=1}^m$ be an $m$-tuple of simplicial pairs on vertex sets $[l_1],\dots,[l_m]$ and let $(\underline{X},\underline{A})$ be a $(l_1 + \cdots + l_m)$-tuple of $CW$-pairs. Then
\begin{equation} \label{eq:polyprodandpolyjoin}
(\underline{X},\underline{A})^{(\underline{\cS},\underline{\cT})^{*\cK} }= \left({\underline{(\underline{X},\underline{A})}^{{\cS}_i}},{\underline{(\underline{X},\underline{A})}^{\cT_i}} \right)^{\cK}. 
\end{equation}
\qed
\end{theorem}

%\todo[inline,color=green]{GS: Reviewer wanted Viadurre stated as Theorem. Done.}

Using the combinatorics of the polyhedral join product, we analyse nested higher Whitehead maps.

By~\eqref{eq:subandcompoverboundaryofsimplex} and~\eqref{eq:polyprodandpolyjoin}, we can rewrite the codomain~\eqref{eq:nestedhwcodom}  of the nested higher Whitehead map as
\[\left( \underline{\left(\underline{Y}_i,\underline{*} \right)}^{\partial \overline{\Delta}^{l_i-1}},\underline{*} \right)^{\partial \Delta^{m-1}} = \left( \underline{Y},\underline{*} \right)^{(\underline{\partial \overline{\Delta}}^{l_i-1},\underline{\emptyset})^{*\partial \Delta^{m-1}}}=\left( \underline{Y},\underline{*} \right)^{\partial \Delta \langle \partial \overline{\Delta}^{l_1-1},\dots,\partial \overline{\Delta}^{l_m-1}\rangle} .\]

By Lemma~\ref{lem:basichwtriviality}, the composition of the nested higher Whitehead map with the inclusion 
\[\left( \underline{Y},\underline{*} \right)^{\partial \Delta \langle \partial \overline{\Delta}^{l_1-1},\dots,\partial \overline{\Delta}^{l_m-1}\rangle} \longrightarrow \left( \underline{Y},\underline{*} \right)^{\Delta \langle \partial \overline{\Delta}^{l_1-1},\dots,\partial \overline{\Delta}^{l_m-1}\rangle} = \left( \underline{\left(\underline{Y}_i,\underline{*} \right)}^{\partial \overline{\Delta}^{l_i-1}},\underline{*} \right)^{ \Delta^{m-1}}\]
induced by the inclusion of simplicial complexes, is null-homotopic. This motivates the question of the description of subcomplexes $\cK$ of $\Delta \langle \partial \overline{\Delta}^{l_1-1},\dots,\partial \overline{\Delta}^{l_m-1}\rangle$ for which the map
\begin{equation} \label{eq:nestedHw}
h_w^{\cK}(f_1,\dots,f_m) \colon (\underline{CX},\underline{X})^{\partial \Delta^{m-1}} \xlongrightarrow{h_w}(\underline{Y},\underline{*})^{\partial \Delta \langle \partial \overline{\Delta}^{l_1-1},\dots,\partial \overline{\Delta}^{l_m-1}\rangle} \longrightarrow (\underline{Y},\underline{*})^{\cK}
\end{equation}
is null-homotopic. The following two propositions give topological and combinatorial conditions under which this occurs.
%\todo[inline,color=green]{For george: Change $\partial \Delta^{l_i-1}$ to $\overline{\partial \Delta}^{l_i-1}$ everywhere, and also define it (same as $\partial \Delta$ except when $l_i = 1$.}
%\todo[inline, color=green]{Make sure we don't call $h_w^{\cK}$ a higher Whitehead map. Only $h_w$ is.}

\goodbreak

\begin{proposition} \label{prop:basicnestedhwtriviality}
Let $\cK$ be a simplicial complex on $[l_1] \sqcup \cdots \sqcup [l_m]$ such that \\ $\partial \Delta \langle \partial \overline{\Delta}^{l_1-1},\dots,\partial \overline{\Delta}^{l_m-1}\rangle \subseteq \cK$. Suppose that for at least one $i \in \{1,\dots,m\}$, $f_i = h_w(f_{i_1},\dots,f_{i_{l_i}}) \colon \SX_i \longrightarrow (\underline{Y}_i,\underline{*})^{\partial \Delta^{l_i-1}}$ with $l_i \geqslant 2$.
If one of the following is satisfied:
\begin{enumerate}[(i)]
    \item \label{1)} the map $f_i\colon  \SX_i \longrightarrow Y_i$ is null-homotopic for some $i = 1,\dots,m$;
    \item \label{2)} $\Delta \langle \partial \overline{\Delta}^{l_1-1},\dots,\partial \overline{\Delta}^{l_m-1}\rangle \subseteq \cK$,
    %\item \label{3)} there exists $i \in [m]$ such that $\cK$ contains $\partial \Delta \langle \partial \Delta^{l_1-1},\dots, \Delta^{l_i-1},\dots, \partial \Delta^{l_m-1}\rangle$, where $f_i$ is a higher Whitehead map,
\end{enumerate}
then $h_w^{\cK}(f_1,\dots,f_m)$ is null-homotopic.
\end{proposition}

\begin{proof}
The claim follows by Lemma~\ref{lem:basichwtriviality}, together with \eqref{eq:polyprodandpolyjoin}.
%Assume that (iii) holds. Let $f_i' = h_w^{\Delta^{l_i-1}}(f_{i_1},\dots,f_{i_{l_i}})$. Since $\partial \Delta^{l_i-1} \subseteq \Delta^{l_i-1}$ then
%\[
%h_w^{\cK} (f_1,\dots,f_i,\dots,f_m) = h_w^{\cK} (f_1,\dots,f_i',\dots,f_m).
%\]
%The statement follows since the map $h_w^{\cK} (f_1,\dots,f_i',\dots,f_m)$ is null-homotopic by Lemma~\ref{lem:basichwtriviality}.
\end{proof}

\begin{proposition} \label{prop:combnestedhwtriviality}

Let $\cK$ be a simplicial complex on $[l_1] \sqcup \cdots \sqcup [l_m]$ such that \\ $\partial \Delta \langle \partial \overline{\Delta}^{l_1-1},\dots,\partial \overline{\Delta}^{l_m-1}\rangle \subseteq \cK$. Suppose that for at least one $i \in \{1,\dots,m\}$, $f_i = h_w(f_{i_1},\dots,f_{i_{l_i}}) \colon \SX_i \longrightarrow (\underline{Y}_i,\underline{*})^{\partial \Delta^{l_i-1}}$ with $l_i \geqslant 2$. If there exists $1 \leqslant i \leqslant m$ such that $\cK$ contains $\partial \Delta \langle \partial \overline{\Delta}^{l_1-1},\dots, \Delta^{l_i-1},\dots, \partial \overline{\Delta}^{l_m-1}\rangle$, where $f_i$ is a higher Whitehead map,
then $h_w^{\cK}(f_1,\dots,f_m)$ is null-homotopic.
    
\end{proposition}

\begin{proof}
    Let $f_i' = h_w^{\Delta^{l_i-1}}(f_{i_1},\dots,f_{i_{l_i}})$. Since $\partial \Delta^{l_i-1} \subseteq \Delta^{l_i-1}$,
\[
h_w^{\cK} (f_1,\dots,f_i,\dots,f_m) = h_w^{\cK} (f_1,\dots,f_i',\dots,f_m).
\]
The statement follows since the map $h_w^{\cK} (f_1,\dots,f_i',\dots,f_m)$ is null-homotopic by Lemma~\ref{lem:basichwtriviality}.
\end{proof}

%\todo[inline, color=green]{Separate out (iii) into another Proposition. Then the Examples make more sense from the perspective of adding simplices. Example of $h_w(h_w(1,2),h_w(3,4),h_w(5,6))$ in $1$-skel of octahedron, and fill in just the $1$-simplices. Then when we fill in enough that the complex from (iii) is there, then naturality kicks in.}

We now apply the higher Whitehead map to the study of higher Whitehead products. Whether the higher Whitehead map is null-homotopic is determined by both the combinatorics of $\cK$ and homotopy theoretic properties of the maps $f_i \colon \SX_i \longrightarrow Y_i$, as shown in Example~\ref{ex:hwAdjointTriv}. In the following example, we study a family of nested higher Whitehead maps of inclusions $S^2 \longrightarrow \DJK$. Such maps were studied algebraically by Abramyan~\cite{Abramyan}. By considering a family of non null-homotopic maps, we construct non null-homotopic maps within trivial higher Whitehead products.

%\todo[inline,color=green]{For george: Rewrite the following sentence: say Abramyan's result first, and then what our example does}

\begin{example} \label{ex:nestedexample}
Let $p \geqslant 2$, $q \geqslant 3$, and $\cS_1 = \partial \Delta^{p-1}$. Define the simplicial complex $\cK$ on the vertex set $[p+q-1] = \{1_1,\dots,1_p\}\sqcup\{2,\dots,q\}$ by 
\begin{align*}
\cK &= \partial \Delta^{q-1} \langle \cS_1, p+1, \dots,p+q-1 \rangle \cup \Delta[1_1,\dots,1_p] \\ 
& = \partial \Delta^{q-1} \langle \partial \Delta[1_1,\dots,1_p],p+1,\dots,p+q-1 \rangle \cup \Delta[1_1,\dots,1_p].
\end{align*}

Let $\mu_{1_i} \colon S^2_{1_i} \longrightarrow \CPi_{1_i}$ and $\mu_j \colon S^2_{j} \longrightarrow \CPi_j$ for $1 \leqslant i \leqslant p$ and $2 \leqslant j \leqslant q$ denote the inclusions of the bottom cell. %Denote by $(\underline{Y},\underline{*})$ the tuple $\{(DJ_{\partial \Delta^{p-1}},*),(\CPi_2,*),\dots,(\CPi_q,*)\}$. 
The map
\[
h_w^{\cK}(h_w(\mu_{1_1},\dots,\mu_{1_p}),\mu_2,\dots,\mu_q) \colon S^{2p+2q-5} \longrightarrow (\CPi,*)^{\partial \Delta^{q-1} \langle \partial \Delta^{q-1},\bullet,\dots,\bullet \rangle} \longrightarrow (\CPi,\ast)^{\cK}
%=DJ_{\partial \Delta^{q-1} \langle \partial \Delta^{q-1},\bullet,\dots,\bullet \rangle} #
\] 
is not null-homotopic, see~\cite[Proposition~7.2]{Abramyan}. By Proposition~\ref{prop:hwisHWPBasic}, this map is an element of the higher Whitehead product $[\iota_1 \circ h_w(f_{1_1},\dots,f_{1_p}),\iota_2 \circ f_2,\dots,\iota_q \circ f_q]$, which is trivial since $\iota_1 \circ h_w(f_{1_1},\dots,f_{1_p})$ is null-homotopic by Proposition~\ref{prop:basicnestedhwtriviality}. We therefore obtain higher Whitehead products with indeterminacy. 

Notice that the case $p=2$ and $q=3$ recovers the result of Abramyan~\cite{Abramyan}.

%Since the map $\iota_1 \circ h_w(f_{1_1},\dots,f_{1_p})$ is null-homotopic, then also $0 \in [\iota_1 \circ h_w(f_{1_1},\dots,f_{1_p}),\iota_2 \circ f_2,\dots,\iota_q \circ f_q]$, by Proposition~\ref{prop:basicnestedhwtriviality}, and we obtain that $h_w(h_w(\mu_{1_1},\dots,\mu_{1_p}),\mu_2,\dots,\mu_q)$ is a non-trivial element of a trivial higher Whitehead product.

% \todo[inline]{For intro: Then go onto say that outside of pathological cases (topological) the same combinatorics controls the non-triviality - we say that non-triviality is controlled both combinatorially and topologically. - may go into intro}

%The case where $p=q=3$ was considered in~\cite{Abramyan} in the special case that for each $1 \leqslant j \leqslant p$ and each $2 \leqslant i \leqslant q$, the maps $f_{1_j}$ and $f_i$ respectively are the inclusion $S^2 \longrightarrow \CPi$. In particular it was observed that there is a wedge summand $S^{10}$ of $\Z_{\cK}$ such that the composite $S^{10} \longhookrightarrow \Z_{\cK} \longrightarrow DJ_{\cK}$ is not a representative of an element of a non-trivial higher Whitehead product. It follows from the above that this composite is the higher Whitehead map $h_w(h_w(f_{1_1},f_{1_2},f_{1_3}),f_2,f_3,f_4)$

\end{example}

\section{Folded higher Whitehead maps} \label{sec:FoldedHw}

%\todo[inline,color=green]{$\pi_*(\Omega X)$ has been studied intensely in homotopy theory, including the quasi Lie structure that Whitehead products endow}

The homotopy groups $\pi_*(X)$ of a space $X$ are a central object of study in homotopy theory. The Whitehead product $[\cdot,\cdot] \colon \pi_m(X) \times \pi_n(X) \longrightarrow \pi_{m+n-1}(X)$ endows the graded module $\pi_*(X) = \bigoplus_{k \geqslant 1} \pi_k(X)$ with the structure of a graded quasi-Lie algebra. Under the isomorphism $\pi_{*+1}(X) \cong \pi_{*} (\Omega X)$, Whitehead products appear as commutators in the algebra $H_*(\Omega X; \mathbb{Z})$ under the Hurewicz map $\pi_*(\Omega X) \longrightarrow H_*(\Omega X;\mathbb{Z})$. Rationalising, the algebra $H_*(\Omega X;\mathbb{Q})$ is the universal enveloping algebra of $\pi_*(\Omega X) \otimes \mathbb{Q}$, meaning that homotopy-theoretic information about $X$ can be extracted from determining the primitive elements in a model for $H_*(\Omega X;\mathbb{Q})$.

For $\cK$ a simplicial complex, the homotopy Lie algebra $\pi_* (\Omega \DJK) \otimes \mathbb{Q}$ was studied by Zhuravleva~\cite{Zhur21}. Using an Adams--Hilton model of $H_*(\Omega \DJK,\mathbb{Q})$, she described the Hurewicz image of the adjoint of higher Whitehead maps.

%\todo[inline,color=green]{Don't talk about relations yet, just focus on the first Lie bracket. Raises the question of whether the combinatorics detects these maps}

For example, let $\cK = \partial \Delta \langle \partial \Delta \langle 1,2,3 \rangle, 4,5 \rangle$ (see \eqref{eq:subandcompoverboundaryofsimplex}) and let $u_i \colon S^1 \longrightarrow \Omega \CP$ be the adjoint of the inclusion $\mu_i \colon S^2 \longrightarrow \CP$ of the bottom cell for $i=1,\dots,5$. The nested brackets $[[u_1,u_2,u_3],[u_1,u_4,u_5]]$ and $[[[u_1,u_2,u_3],u_4,u_5],u_1]$, the adjoints to $[h_w^{\cK}(\mu_1,\mu_2,\mu_3),h_w^{\cK}(\mu_1,\mu_4,\mu_5)]$ and $[h_w^{\cK}(h_w(\mu_1,\mu_2,\mu_3),\mu_4,\mu_5),\mu_1^{\cK}]$, respectively, are defined in $\pi_8 (\Omega \DJK) \otimes \mathbb{Q}$.

We study these Whitehead products from the geometric perspective of the higher Whitehead map, and construct the maps considered by Zhuravleva integrally. Our construction realises $2$-torsion elements of homotopy groups, by giving a combinatorial interpretation of higher Whitehead maps with repeated factors.

%\todo[inline,color=green]{Sentence I wrote that doesn't seem to fit anywhere: we combinatorially interpret higher Whitehead maps of maps that repeat}

%\todo[inline,color=green]{Plain english: start with a simplicial complex in which hw map is defined (non-trivial). Is the operation of identifying vertices preserved by hw map? Our main study is when this simplicial fold produces a higher Whithead map (ie. induces a map of polyhedral products)}

Let $\cK$ be a simplicial complex realising a nested higher Whitehead map
\[
h_w^{\cK}(f_1,\dots,f_m) \colon (\underline{CX},\underline{X})^{\partial \Delta} \longrightarrow (\underline{Y},\underline{\ast})^{\cK}.
\]
We study when the simplicial fold operation $\cK \longrightarrow \overline{\cK}$ identifying two or more vertices preserves the structure of a higher Whitehead map. Specifically, we ask when the simplicial fold $\cK \longrightarrow \overline{\cK}$ induces a map of polyhedral products $\nabla \colon (\underline{Y},\underline{\ast})^{\cK} \longrightarrow (\underline{Y},\underline{\ast})^{\overline{\cK}}$, in which case we will call the composite
\[
\nabla h_w(f_1,\dots,f_m) \colon (\underline{CX},\underline{X})^{\partial \Delta} \longrightarrow (\underline{Y},\underline{\ast})^{\cK} \longrightarrow  (\underline{Y},\underline{\ast})^{\overline{\cK}}
\]
a folded higher Whitehead map. 
%One situation where this is possible is when for the vertices being folded together, the associated spaces $Y_i=Y$, for $Y$ an associative $H$-space, see~\eqref{eq:BigHMap}.

To build up to the definition of the folded higher Whitehead map, we begin by defining a fold of a simplicial complex. 

\begin{definition}\label{def:foldofsimpcmplxs}
A \emph{fold} of a vertex set $[m]$ consists of disjoint subsets $I,J \subset [m]$ and a surjective map $\psi \colon I \longrightarrow J$. Given a fold $\psi \colon I \longrightarrow J$ of $[m]$, we denote the pre-images of each $j \in J$ by $I_j = \psi^{-1} (j)$, so that $I = \bigsqcup_{j \in J} I_j$.

Let $\cK$ be a simplicial complex on $[m]$. A fold $\psi \colon I \longrightarrow J$ of $[m]$ extends to a map on $\cK$ by sending a simplex $(i_1,\dots,i_k) \in \cK$ to $(\overline{\psi}(i_1),\dots,\overline{\psi}(i_k))$, where $\overline{\psi}(i) = \psi(i)$ if $i \in I$, and $\overline{\psi}(i) = i$, otherwise. We interchange freely between referring to $\psi \colon I \longrightarrow J$ as a fold of $\cK$ and of its vertex set $[m]$.

We define the folded complex $\cK_{\nabla(I,J)}$ of the fold of $\cK$ induced by $\psi$ as the image $\overline{\psi}(\cK)$. Then, $\cK_{\nabla(I,J)}$ is a simplicial complex on vertex set $[m]\setminus I$.
\end{definition}

%\todo[inline, color=green]{Folding of $\cK$ refers to the map, fold of $\cK$ refers to the image}

It follows from the definition that
\begin{equation} \label{eq:FoldedComplex}
\cK_{\nabla(I,J)} = \{\sigma \subseteq [m] \setminus I \; | \; \sigma \in \cK \text{ or } (\sigma \setminus \{j\}) \sqcup \{i\} \in \cK \text{ for some } j \in J, i \in I_j\}.
\end{equation}
When $J= \{j\}$ consists of one element we abbreviate $\cK_{\nabla(I,\{j\})}$ to $\cK_{\nabla(I,j)}$. If further $I = \{i\}$, then we abbreviate $\cK_{\nabla(\{i\},j)} = \cK_{\nabla(i,j)}$.

\begin{example} \label{ex:FoldingExample1}
Let $\cK$ be the simplicial complex shown on the left of Figure~\ref{fig:FoldingExample1}. Let $I = \{4\}$, $J = \{1\}$ and let $\psi \colon I \longrightarrow J$ be the map sending $\{4\}$ to $\{1\}$. Then $\cK_{\nabla(4,1)}$ is the simplicial complex on vertices $\{1,2,3\}$ shown on the right of Figure~\ref{fig:FoldingExample1}. 
\begin{figure}[h]
    \centering
    \begin{tikzpicture}[scale=0.7]
        \coordinate [label=right:{1}] (1) at (0,0,0);
        \coordinate [label=below left:{2}] (2) at (-3,0,1.72);
        \coordinate [label=above:{3}] (3) at (-3,0,-1.72);
        \coordinate [label=left:{4}] (4) at (-6,0,0);
        \draw (1) -- (2) -- (3) -- (1);
        \draw (4) -- (2) -- (3) -- (4);
        \foreach \point in {1,2,3,4}
            \fill [black] (\point) circle (1.5 pt);
    \end{tikzpicture}
    \hspace{0.05 \textwidth}
    \begin{tikzpicture}[scale=0.7]
        \coordinate [label=right:{1}] (1) at (0,0,0);
        \coordinate [label=left:{2}] (2) at (-3,0,1.72);
        \coordinate [label=above:{3}] (3) at (-3,0,-1.72);
        \coordinate [label=right:{4}] (4) at (0,1,0);
        \draw (1) -- (2) -- (3) -- (1);
        \draw (4) -- (2) -- (3) -- (4);
        \foreach \point in {1,2,3,4}
            \fill [black] (\point) circle (1.5 pt);
    \end{tikzpicture}
    \hspace{0.05 \textwidth}
    \begin{tikzpicture}[scale=0.7]
        \coordinate [label=right:{1}] (1) at (0,0,0);
        \coordinate [label=left:{2}] (2) at (-3,0,1.72);
        \coordinate [label=above:{3}] (3) at (-3,0,-1.72);
        \draw (1) -- (2) -- (3) -- (1);
        \foreach \point in {1,2,3}
            \fill [black] (\point) circle (1.5 pt);
    \end{tikzpicture}
    \caption{The fold $\psi \colon \{4\} \longrightarrow \{1\}$ of a simplicial complex.}
    \label{fig:FoldingExample1}
\end{figure}
\end{example}

%\todo[inline,color=green]{Also just do the completely unfolded one on the left.}

The following properties follow from Definition~\ref{def:foldofsimpcmplxs}.

\begin{proposition}
\label{prop:propsoffoldofsimpcmplx}\
Let $\cK$ be a simplicial complex on $[m]$. Then:
\begin{enumerate}[(i)]
    \item if $i,j \in [m]$ with $i \neq j$, then $\cK_{\nabla(i,j)} \cong \cK_{\nabla(j,i)}$;
    \item if $I = \{i_1,\dots,i_n\} \subseteq [m]$, and $j \in [m]$ with $j \notin I$, then
\[\cK_{\nabla(I,j)} =  {{{\cK_{\nabla(i_1,j)}}_{\nabla(i_2,j)}}_{\cdots}}_{\nabla(i_n,j)}; \]
    \item if $I_1, I_2 \subseteq [m]$ are such that $I_1 \cap I_2 = \emptyset$, and $j \notin I_1 \sqcup I_2$, then
\[\pushQED{\qed} 
{\cK_{\nabla(I_1,j)}}_{\nabla(I_2,j)} = {\cK_{\nabla(I_2,j)}}_{\nabla(I_1,j)}.  \qedhere
\popQED \]
    %\item Let $\cK_1,\cK_2$ be simplicial complexes on $[m_1], [m_2]$ respectively. Let $\cL_1 \cong \cL_2 \cong \cL$ be isomorphic full subcomplexes of $\cK_1$ and $\cK_2$ respectively, with vertex sets $\{l_{1_1},\dots,l_{n_1}\}$ and $\{l_{1_2},\dots,l_{n_2}\}$ respectively, where $l_{n_1} = l_{n_2}$. Let $\cK = \cK_1 \sqcup \cK_2$. Then
%\[\cK_1 \bigcup_{\cL} \cK_2 \cong {{{\cK_{\nabla(l_{1_1},l_{1_2})}}_{\nabla(l_{2_1},l_{2_2})}}_{...}}_{\nabla(l_{n_1},l_{n_2})}.\]
\end{enumerate}
\end{proposition}

%\todo[inline]{Which of these do we actually use?}

%It follows from Proposition~\ref{prop:propsoffoldofsimpcmplx} that if $J = \{j_1,\dots,j_k\}$ we have
%\[\cK_{\nabla(I,J)} = {\cK_{{\nabla(I_{j_1},j_1)}_{\cdots}}}_{\nabla(I_{j_k},j_k)}\]
%We now define a fold of full subcomplexes. Suppose that for $J_1 = \{{j_1}_1,\dots,{j_1}_{n_1} \}, J_2 = \{{j_2}_1,\dots,{j_2}_{n_2} \} \subseteq [m]$ with $J_1 \cap J_2 = \emptyset$, a map $\psi\colon J_1 \longrightarrow J_2$, ${j_1}_l \longmapsto {j_2}_l$ is given. Then the fold of $\cK_{J_1}$ onto $\cK_{J_2}$ by $\psi$ is defined
%\begin{equation} \label{eqref:foldoffullsubcmplxs} \cK_{\nabla(J_1,J_2)_{\psi}} = {{\cK_{\nabla({j_1}_1,\psi({j_1}_1))}}_{\dots}}_{\nabla({j_1}_n,\psi({j_1}_n))} \end{equation}
%We denote by $\nabla(J_1,J_2)_{\psi} \colon \cK \longrightarrow \cK_{{\nabla(J_1,J_2)}_{\psi}}$ the induced map of simplicial complexes. Where the map $\psi$ does not need to be emphasized, we omit it from the notation.

To describe the map of polyhedral products induced by a fold of simplicial complexes, let $\psi \colon I \longrightarrow J$ be a fold of $\cK$. Suppose that $Y_j$ is an associative $H$-space for each $j \in J$, and that $Y_i = Y_j$ for each $i \in I_j$. For $j \in J$, the $H$-multiplication map $\mu_j \colon Y_j \times Y_j \longrightarrow Y_j$ extends up to homotopy to a map
\[
\prod_{i \in I_j} Y_i \times Y_j \longrightarrow Y_j.
\]
Therefore, the fold $\psi \colon I \longrightarrow J$ induces a map $\prod_{i \in [m]} Y_i \longrightarrow \prod_{i \in ([m] - I)} Y_i$ given by
\begin{equation} \label{eq:BigHMap}
\prod_{j \in J} \prod_{i \in I_j} (Y_i \times Y_j) \times \prod_{i \in ([m] - (I \sqcup J))} Y_i \longrightarrow \prod_{j \in J} Y_j \times \prod_{i \in ([m] - (I \sqcup J))} Y_i.
\end{equation}

The fold $\psi \colon I \longrightarrow J$ of $\cK$ induces the map of polyhedral products
\begin{equation} \label{eq:InducedFoldMap}
\nabla_{(I,J)} \colon (\underline{Y},\underline{\ast})^{\cK} = \bigcup_{\sigma \in \cK} \prod_{i \in \sigma} Y_i \longrightarrow \bigcup_{\sigma \in \cK} \prod_{i \in \overline{\psi}(\sigma)} Y_i = (\underline{Y},\underline{\ast})^{\cK_{\nabla(I,J)}}
\end{equation}
defined as the restriction of~\eqref{eq:BigHMap} to $(\underline{Y},\underline{\ast})^{\cK}$.

We call $(\underline{Y},\underline{\ast})^{\cK_{\nabla(I,J)}}$ the fold of $(\underline{Y},\underline{\ast})^{\cK}$ under the map $\psi \colon I \longrightarrow J$, and the map $\nabla_{(I,J)}$ the fold map on polyhedral products induced by the map $\psi \colon I \longrightarrow J$.

%\todo[inline]{The following paragraph has been added}

Since the multiplication $\mu \colon Y_j \times Y_j \longrightarrow Y_j$ extends the fold map $Y_j \vee Y_j \longrightarrow Y_j$ up to homotopy, the fold $\psi \colon I \longrightarrow J$ of $\cK$ also induces a map of polyhedral products whenever $I_j$ is such that $\cK_{I_j\cup\, j}$ consists of disjoint vertices. In this case, the map~\eqref{eq:InducedFoldMap} is defined without requiring $Y_j$ to be an $H$-space. %We will return to this observation when considering relations among folded higher Whitehead maps in Theorem~\ref{thm:GeneralUnconnectedFold}.

%\todo[inline]{Linking sentence needed.} % This first suggestion doesn't seem to fit: The folding of simplicial complexes can be extended to incorporate the structure of the polyhedral join product, folding on $\cK$ and folding on $\cK_i$. We specialise to substitution to define folded nested higher Whitehead map}

We specialise to fold maps of polyhedral products related to higher Whitehead maps.

%The composite of a higher Whitehead map with a fold map is an element of a higher Whitehead product with repeated factors, whose homotopy-theoretic properties of these composites be analysed combinatorially. 

\begin{definition} \label{def:FoldedhwMap}
The \emph{folded higher Whitehead map} is the composite
\[
\nabla_{(I,J)} h_w^{\cK}(f_1,\dots,f_m) \colon (\underline{CX},\underline{X})^{\partial \Delta} \longrightarrow (\underline{Y},\underline{\ast})^{\cK} \longrightarrow (\underline{Y},\underline{\ast})^{\cK_{\nabla(I,J)}}.
\] 
\end{definition}

In analogy with the nested higher Whitehead map, the folded higher Whitehead map is an element of a particular higher Whitehead product.

\begin{proposition}
\label{prop:foldedhwmapishwprod}
The folded higher Whitehead map is an element of the higher Whitehead product
\[
\nabla_{(I,J)}  h_w^{\cK}(f_1,\dots,f_m) \in [\nabla_{(I,J)} \circ \iota_1 \circ f_1,\dots, \nabla_{(I,J)} \circ \iota_m \circ f_m] \subseteq \left[X_1 \ast \cdots \ast X_m, (\underline{Y},\underline{*})^{\cK_{\nabla(I,J)}} \right]\]
where $\iota_j$ denotes the inclusion of $Y_j$ into the $j$th coordinate.
\end{proposition}

\begin{proof}
    This follows from Proposition~\ref{prop:hwisHWPBasic}, together with the naturality of the higher Whitehead product \cite[Theorem~2.1]{Porter1965HigherOW}.
\end{proof}

Let 
\[
h_w^{\cK}(f_1,\dots,f_m) \colon (\underline{CX},\underline{X})^{\partial \Delta} \longrightarrow (\underline{Y},\underline{\ast})^{\cK}
\]
be the nested higher Whitehead map, as in Definition~\ref{def:nestedhwmap}. We study combinatorial conditions under which folded nested higher Whitehead maps are null-homotopic. We begin with some preparatory results. First, we show that if two simplicial complexes fold into the same folded complex, then the associated folded higher Whitehead maps are equal. 

\begin{lemma}
\label{lem:iffoldsequalthenfoldedhwmapsequal}
Let $\cS_1,\dots,\cS_m$ be simplicial complexes on $[l_1],\dots,[l_m]$, respectively. Let $f_i \colon \Sigma X_i \longrightarrow (\underline{Y},\underline{\ast})^{\cS_i}$ be maps for $i=1,\dots,m$.

Let $\cK, \cK'$ be simplicial complexes on $[l] = [l_1] \sqcup \cdots \sqcup [l_m]$ containing $\partial \Delta \langle \cS_1,\dots,\cS_m \rangle$ and let $\psi \colon I \longrightarrow J$ be a fold of both $\cK$ and $\cK'$ such that $\cK_{\nabla(I,J)} = \cK'_{\nabla(I,J)}$. If $Y_j$ is an associative $H$-space for each $j \in J$, and $Y_i = Y_j$ for each $i \in I_j$, then
\[
\nabla_{(I,J)}  h_w^{\cK}(f_1,\dots,f_m) = \nabla_{(I,J)}  h_w^{\cK'}(f_1,\dots,f_m).
\]
\end{lemma}

\goodbreak

\begin{proof}
Since $\cK_{\nabla(I,J)} = \cK'_{\nabla(I,J)}$, there is a commutative diagram
\begin{equation*}
    \begin{tikzcd}
    \partial \Delta \langle \cS_1,\dots,\cS_m \rangle \ar[r] \ar[d] & \cK \ar[d,"{\overline{\psi}}"] \\
    \cK' \ar[r,"{\overline{\psi}}"] & \cK_{\nabla(I,J)}.
    \end{tikzcd}
\end{equation*}
Therefore, by functoriality of the polyhedral product with respect to simplicial inclusions, the following diagram commutes
\begin{equation*}
    \begin{tikzcd}
    X_1 \ast \cdots \ast X_m \ar[r] &  (\underline{Y},\underline{*})^{\partial \Delta \langle \cS_1,\dots, \cS_m \rangle} \ar[r] \ar[d] & (\underline{Y},\underline{*})^{\cK} \ar[d] \\
    & (\underline{Y},\underline{*})^{\cK'} \ar[r] & (\underline{Y},\underline{*})^{\cK_{\nabla(I,J)}}
    \end{tikzcd}
\end{equation*}
where the composite around the top of the square is $\nabla_{(I,J)}  h_w^{\cK}(f_1,\dots,f_m)$ and the composite around the bottom of the square is $\nabla_{(I,J)}  h_w^{\cK'}(f_1,\dots,f_m)$.
\end{proof}

%\todo[inline,color=green]{The triviality of the folded map does not depend solely on the complex $\cK$, since it can factor through a different complex $\cK'$ which trivialises the map}

The folded higher Whitehead map is null-homotopic in the case that the fold $\nabla \colon \cK \longrightarrow \cK_{\nabla(I,J)}$ factors through an inclusion $\cK \longrightarrow \cK'$ such that $h_w^{\cK'}(f_1,...,f_m)$ is null-homotopic. This observation enables us to obtain a condition for when a folded higher Whitehead map $\nabla_{(I,J)}h_w^{\cK}(f_1,\dots,f_m)$ is null-homotopic. We begin by defining a simplicial complex $\cL_{\psi}$ containing all complexes which fold onto $\cK_\nabla(I,J)$ by $\psi \colon I \longrightarrow J$.

\begin{definition} \label{def:MaxFoldingComplex}
Let $\cK$ be a simplicial complex on $[m]$, and let $\psi \colon I \longrightarrow J$ be a fold of $\cK$ for $I,J \subset [m]$. Denote by $\{k_1,\dots,k_l\}$ the vertex set of the folded complex $\cK_{\nabla(I,J)}$, and let $I_j = \psi^{-1} (j)$ for each $j \in J$. Define the simplicial complex $\cL_{\psi}$ on $[m]$ as the substitution complex, see~\eqref{eq:subsaspolyjoin}, given by
\[
\cL_{\psi} = \cK_{\nabla(I,J)} \langle \Delta[\{k_1\} \sqcup I_{k_1}],\dots,\Delta[\{k_l\} \sqcup I_{k_l}] \rangle
\]
where $I_k = \emptyset$ if $k \notin J$.
\end{definition}

%\todo[color=green, inline]{GS: Reviewer wants a reference back to substitution complexes, todo once we decide if we're putting it as a definition. MS: Done. Referred back to equation.}

\begin{example} \label{ex:FoldingExample3}
    Consider the simplicial complex $\cK$ and the fold map $\psi \colon I \longrightarrow J$ from Example~\ref{ex:FoldingExample1}. The complex $\cK_{\nabla(I,J)}$ is shown on the left of Figure~\ref{fig:FoldingExample3}, while the complex $\cL_{\psi} = \cK_{\nabla(I,J)} \langle \Delta[1,4],2,3 \rangle$ is shown on the right.    
\begin{figure}[h]
    \centering
    \begin{tikzcd}
        \coordinate [label=right:{1}] (1) at (0,0,0);
        \coordinate [label=left:{2}] (2) at (-3,0,1.72);
        \coordinate [label=above left:{3}] (3) at (-3,0,-1.72);
        \draw (1) -- (2) -- (3) -- (1);
        \foreach \point in {1,2,3}
            \fill [black] (\point) circle (1.5 pt);
    \end{tikzcd}
    \hspace{0.1 \textwidth}
    \begin{tikzcd}
        \coordinate [label=right:{1}] (1) at (0,0,0);
        \coordinate [label=left:{2}] (2) at (-3,0,1.72);
        \coordinate [label=above left:{3}] (3) at (-3,0,-1.72);
        \coordinate [label=right:{4}] (4) at (0,1,0);
        \filldraw[fill=black!30] (3) -- (4) -- (1);
        \filldraw[fill=black!30,opacity=0.5] (2) -- (4) -- (1);
        \draw (1) -- (2) -- (3) -- (1);
        \draw (4) -- (2) -- (3) -- (4) -- (1);
        \foreach \point in {1,2,3,4}
            \fill [black] (\point) circle (1.5 pt);
    \end{tikzcd}
    \caption{Construction of the simplicial complex $\cL_{\psi}$.}
    \label{fig:FoldingExample3}
\end{figure}
\end{example}

The complex $\cL_{\psi}$ is the largest simplicial complex which folds to $\cK_{\nabla(I,J)}$ under $\psi \colon I \longrightarrow J$, and therefore in particular contains $\cK$.

\begin{lemma}
\label{lem:iffoldL'containedinfoldK'thenL'containedinL}
Suppose that $\cL'$ is such that $\cL'_{\nabla(I,J)} \subseteq \cK_{\nabla(I,J)}$. Then $\cL' \subseteq \cL_{\psi}$.
\end{lemma}

\begin{proof}
Let $J = \{j_1,\dots,j_r\}$. For any $\sigma \in \cL'$, write $\sigma = \sigma_{j_1} \sqcup \cdots \sqcup \sigma_{j_r} \sqcup \sigma'$,  where $\sigma_{j_k} \subseteq I_{j_k}$ and $\sigma' \in [m] \setminus I$. Since $\cL'_{\nabla(I,J)} \subseteq \cK_{\nabla(I,J)}$, it follows that $\overline{\psi} (\sigma) = (j_1 \sqcup \cdots \sqcup j_r) \cup \sigma' \in \cK_{\nabla(I,J)}$. By construction of $L_{\psi}$, $(I_{j_1} \sqcup \cdots \sqcup I_{j_r}) \cup (j_1 \sqcup \cdots \sqcup j_r) \cup \sigma' \in \cL_{\psi}$. As $\sigma_{j_k} \subseteq I_{j_k}$, it follows that $\sigma \in \cL_{\psi}$. 
\end{proof}

\goodbreak

Finally, we state combinatorial conditions under which the folded higher Whitehead map is null-homotopic.

\begin{proposition} \label{prop:hwFoldedTrivialityConditions}
Let $\cS_1,\dots,\cS_m$ be simplicial complexes on $[l_1],\dots,[l_m]$, respectively. Let $f_i \colon \Sigma X_i \longrightarrow (\underline{Y},\underline{\ast})^{\cS_i}$ be maps for $i=1,\dots,m$. Let $\cK$ be a simplicial complex on $[l_1] \sqcup \cdots \sqcup [l_m]$ containing $\partial \Delta \langle \cS_1,\dots,\cS_m \rangle$ and let $\psi \colon I \longrightarrow J$ be a fold of $\cK$. 

The folded higher Whitehead map
\[\nabla_{(I,J)}  h_w^{\cK}(f_1,\dots,f_m) \colon (\underline{CX},\underline{X})^{\partial \Delta} \longrightarrow  (\underline{Y},\underline{\ast})^{\cK_{\nabla(I,J)}}\]
is null-homotopic if one of the following holds:
\begin{enumerate}[(i)]
    \item $\cK_{\nabla(I,J)}$ contains the folded complex $\partial \Delta \langle \cS_1,\dots,\cS_i',\dots,\cS_m \rangle_{\nabla(I,J)}$, for $\cS_i'$ such that $f_i \colon \Sigma X_i \longrightarrow (\underline{Y},\underline{\ast})^{\cS_i}$ is trivial;
    \item $\cK_{\nabla(I,J)}$ contains the folded complex $\Delta \langle \cS_1,\dots,\cS_m \rangle_{\nabla(I,J)}$;
    \item the map
    \[
    h_w^{\cL_{\psi}} (f_1,\dots,f_m) \colon (\underline{CX},\underline{X})^{\partial \Delta} \longrightarrow (\underline{Y},\underline{\ast})^{\cL_{\psi}}
    \]
    is null-homotopic.
\end{enumerate}
\end{proposition}

\begin{proof}
Suppose that (i) holds. Then $h_w^{\cK'} (f_1,\dots,f_m)$ is null-homotopic by Proposition~\ref{prop:basicnestedhwtriviality}. By  Lemma~\ref{lem:iffoldsequalthenfoldedhwmapsequal}, we have $\nabla_{(I,J)} h_w^{\cK}(f_1,\dots,f_m) = \nabla_{(I,J)} h_w^{\cK'} (f_1,\dots,f_m)$, and so the map $\nabla_{(I,J)} h_w^{\cK}(f_1,\dots,f_m)$ is null-homotopic. The same argument holds if instead (ii) holds, since $\cL_{\psi}$ folds to $\cK_{\nabla(I,J)}$ by Lemma~\ref{lem:iffoldL'containedinfoldK'thenL'containedinL}.
\end{proof}

\begin{example}
Consider the simplicial complex $\cK$ and the fold map $\psi \colon I \longrightarrow J$ from Example~\ref{ex:FoldingExample1}. Let $f_i \colon \Sigma X_i \longrightarrow Y_i$ be maps for $i=1,\dots,4$. The nested higher Whitehead maps $h_w^{\cK}(h_w(f_1,f_2,f_3),f_4)$ and $h_w^{\cK}(h_w(f_1,f_4),f_2,f_3)$ are defined in $(\underline{Y},\underline{\ast})^{\cK}$. Therefore the folded higher Whitehead maps $\nabla_{(I,J)}h_w^{\cK}(h_w(f_1,f_2,f_3),f_4)$ and $\nabla_{(I,J)}h_w^{\cK}(h_w(f_1,f_4),f_2,f_3)$ are defined in $(\underline{Y},\underline{\ast})^{\cK_{\nabla(I,J)}}$. By Proposition~\ref{prop:hwFoldedTrivialityConditions}, the folded higher Whitehead map $\nabla_{(I,J)}h_w^{\cK}(h_w(f_1,f_4),f_2,f_3)$ is null-homotopic since $h_w^{\cL_{\psi}}(h_w(h_w(f_1,f_4),f_2,f_3))$ is null-homotopic by Proposition~\ref{prop:basicnestedhwtriviality}. 

By Proposition~\ref{prop:foldedhwmapishwprod}, for $\cL = \cK_{\nabla(I,J)}$, these folded higher Whitehead maps are elements of the higher Whitehead products $[h_w^{\cL}(f_1,f_2,f_3),f_1^{\cL}]$ and $[\nabla_{(I,J)} h_w^{\cK}(f_1,f_4),f_2^{\cL},f_3^{\cL}]$, respectively. Furthermore, by Proposition~\ref{prop:foldedhwmapishwprod}, the map $\nabla_{(I,J)} h_w^{\cK}(f_1,f_4)$ is the composite of $[f_1,f_1]$ with the inclusion $Y_1 \longrightarrow (\underline{Y},\underline{\ast})^{\cL}$.

    %On the other hand, the map $\nabla_{(I,J)} h_w^{\cK}(h_w(f_1,f_2,f_3),f_4)$ is not null-homotopic (in the case of $\DJK$, by Zhuravleva,..?). In Example~\ref{ex:TorsionExample}, we will show that this map is $2$-torsion (again in $\DJK$).

    %\todo[inline]{Finish example, ie. when and how can we say it's not trivial - we currently don't have any way to do this, even for $\DJK$.}
\end{example}

\section{Relations}

%\todo[inline,color=green]{General introduction: the set $[\Sigma X, Y]$ is a group, and we analyse relations among higher Whitehead maps. Motivated by Jacobi and generalisations to Hardie.}

%The set of homotopy classes of maps $[\SX,Y]$ inherits a group structure from the comultiplication in $\SX$. The Jacobi identity, and Hardie's generalisation \eqref{eq:HardieIdentity2}, are relations among Whitehead products and exterior Whitehead products, respectively. These examples motivate us to analyse relations among higher Whitehead maps more generally.
%\todo[inline]{An introductory paragraph?}
We recall the relations among exterior Whitehead products (see~\eqref{eq:Hardiehwmap}) due to Hardie~\cite{Hardie61}. Let $f_i \colon S^{q_i} \longrightarrow Y_i$ be maps for $i = 1,\dots,m$, where $q_i \geqslant 2$, and let
\begin{align*}
Z &= \bigcup_{i=1}^m Y_i \vee FW(Y_1,\dots,Y_{i-1},Y_{i+1},\dots,Y_m) 
\end{align*}
denote the union of subspaces in $Y_1 \times \cdots \times Y_m$. Then
\begin{align} \label{eq:HardieIdentity2}
& \sum_{i=1}^m (-1)^{\eta(i)} [\theta_i f_i,\psi_i \{f_1,\dots,f_{i-1},f_{i+1},\dots,f_m\}] = 0 \in \pi_{q-2}(Z)
\end{align}
where $q= q_1 + \cdots + q_m -1$, the maps $\theta_i \colon Y_i \longrightarrow Z$ and $\psi_i \colon FW(Y_1,\dots,Y_{i-1},Y_{i+1},\dots,Y_m) \longrightarrow Z$ are inclusions, and $\eta(i) = q_i (q_1 + \cdots + q_i) + 1$.

We rewrite relation~\eqref{eq:HardieIdentity2} using polyhedral products, enabling detection of the combinatorial structure which captures it. Each summand $\{\theta_i f_i,\psi_i \{f_1,\dots,f_{i-1},f_{i+1},\dots,f_m\}\}$ is a map 
\begin{comment}
    of polyhedral products with domain $(\underline{D},\underline{S})^{\partial \Delta^1} = \Sigma^{m-2} \bigwedge_{i=1}^m S^{q_i-1}$ and
\end{comment} 
with codomain $Z$ expressed using polyhedral products as
%\todo[inline, color=green]{GS: Before we tried to express the domain of each Hardie map as a polyhedral product, and therefore each map as a ``map of polyhedral products''. I don't believe this added anything, and it caused confusion for the reviewer.}
\begin{align*}
Z &= \bigcup_{i=1}^m Y_i \vee FW(Y_1,\dots,Y_{i-1},Y_{i+1},\dots,Y_m)  \\
& = \left\{ (y_1,\dots,y_m) \in Y_1 \times \cdots \times Y_m \mid y_i = y_j = \ast \text{ for } 1 \leqslant i < j \leqslant m \right\} \\
& = (\underline{Y},\underline{\ast})^{\sk^{m-3} \Delta^{m-1}}.
\end{align*}

Each minimal missing face $\{[m] \setminus \{i\} \mid i = 1,\dots,m\}$ of $\sk^{m-3} \Delta^{m-1}$ gives rise to the higher Whitehead map $h_w(f_1,\dots,f_{i-1},f_{i+1},\dots,f_m)$ in $(\underline{Y},\underline{\ast})^{\sk^{m-3} \Delta^{m-1}}$. Summing, we rewrite relation~\eqref{eq:HardieIdentity2} as a sum of nested higher Whitehead maps
\begin{equation} \label{eq:HardieIdentityAsHw}
\sum_{i=1}^m (-1)^{\eta(i)} h_w^{\sk^{m-3} \Delta^{m-1}} (h_w(f_1,\dots,f_{i-1},f_{i+1},\dots,f_m),f_i) = 0.
\end{equation}

%\todo[inline]{Emphasise that composing into $\sk$ doesn't change the higher Whitehead map -- what does this mean?}

We generalise this relation to a family of relations among terms of differing forms, replacing the spherical maps $f_i$ with maps $\Sigma X_i \longrightarrow Y_i$, where $X_i$ is itself a suspension.

To enable the study of relations among terms of different forms, we propagate the missing face structure of $\sk^{m-3} \Delta^{m-1} = \sk^{k-3} \Delta^{k-1} (\emptyset,\dots,\emptyset)$ (see~\eqref{eq:compaspolyjoin}) by replacing $\emptyset = \partial \Delta^0$ with $\partial \Delta^n$ for various $n > 0$.

%\todo[color=green,inline]{GS:Again, we're wanted to reference the definition of composition here MS: Done}

%\todo[inline,color=green]{We got to here 24/01/23}

Given a vertex set $[m]$, a $k$-partition $\Pi$ is a collection of pairwise disjoint subsets $\{P_1,\dots,P_k\}$ of $[m]$ such that $\bigcup_{i=1}^k P_i = [m]$.

\begin{definition} \label{def:identitycomplex}
Let $\Pi = \{P_1,\dots,P_k\}$ be a $k$-partition of $[m]$, where $P_j = \{i_1,\dots,i_{n_j} \}$  for $j = 1,\dots,k$. The \textit{identity complex} $\cK_{\Pi}$ associated to $\Pi$ is the composition complex (see~\eqref{eq:compaspolyjoin}) defined by
\[ 
\cK_{\Pi} = \sk^{k-3} \Delta^{k-1} (\partial \Delta[P_1],\dots,\partial \Delta[P_k]).
\]
\end{definition}

%\todo[inline,color=green]{Short: start with dual to Hardie being disjoint points, then show that alexander dual of $\cK_{\Pi}$ is disjoint simplices. Reference this in Prop 5.4 proof.}

The identity complex $\sk^{m-3} \Delta^{m-1}$ associated to Hardie's identity corresponds to the partition $\Pi = \{\{1\},\dots,\{m\}\}$ of $[m]$.

Identity complexes admit a similar description in terms of their Alexander duals, which we describe next. The Alexander dual $\hat{\cK}$ of a simplicial complex $\cK$ on $[m]$ is the simplicial complex with simplices $\{[m] \setminus \sigma \mid \sigma \notin \cK\}$. In particular, the maximal faces of $\hat{\cK}$ are the complements in $[m]$ of the minimal missing faces of $\cK$. Since the minimal missing faces of $\sk^{m-3} \Delta^{m-1}$ are $\{[m] \setminus \{i\} \mid i=1,\dots,m\}$, it is Alexander dual to the simplicial complex consisting of $m$ disjoint vertices,
\[\sk^{m-3} \Delta^{m-1} = \widehat{\bullet_{[m]}}.\]

 For a simplicial complex $\cK$ let $MF(\cK)$ denote the set of its minimal missing faces. %Let $Q_i = \{j_1,\dots,j_{r_i}\} = [m] \setminus P_i$ for each $i=1,\dots,k$. 

\begin{proposition} \label{prop:MFsofKPi}
Let $\Pi = \{P_1,\dots,P_k\}$ be a $k$-partition of $[m]$. Then
\begin{equation} \label{eq:MFSofKPi}
MF(\cK_{\Pi}) = \{[m] \setminus P_i \mid i = 1,\dots,k\}.
\end{equation}
\end{proposition}

To prove Proposition~\ref{prop:MFsofKPi}, we first compute the minimal missing faces of the polyhedral join product. 

%\todo[inline]{Cite this from Duality paper, and put directly in proof for 5.2. Proof for 5.2 follows the statement of 5.2 immediately -- it is not proved in the duality paper}

\begin{proposition} \label{prop:missingfacespolyhedraljoin}
Let $\cK$ be a simplicial complex on $[m]$, and let $(\cS_1, \cT_1),...,(\cS_m, \cT_m)$ be simplicial pairs on vertex sets $[l_1],...,[l_m]$, respectively. Then
\[MF\left( (\underline{\cS},\underline{\cT})^{* \cK} \right) = \left\{\tau \in MF(\cS_i) \; | \; i \in \cK \right\} \sqcup \left\{\bigsqcup_{i \in \kappa} \tau_i \; | \; \kappa \in MF(\cK), \tau_i \in MF(\cT_i), \tau_i \in \cS_i \right\}.\]
\end{proposition}

\begin{proof}
We first show that 

\[
\begin{aligned} \{J \in MF(\cS_i) \; | \; i \in \cK \} &\sqcup \left\{ \bigsqcup_{i \in L} J_i \; | \; L \in MF(\cK), J_i \in MF(\cT_i), J_i \in \cS_i \right\} \\ & \quad \subseteq MF\left( (\underline{\cS},\underline{\cT})^{* \cK} \right).
\end{aligned}\]

For any $\{i\} \in \cK$, $J \in MF(\cS_i)$ implies that $J \in MF\left( (\underline{\cS},\underline{\cT})^{* \cK} \right)$. Now consider $\bigsqcup_{i \in L} J_i$, where $L \in MF(\cK)$, $J_i \in MF(\cT_i)$ and $J_i \in \cS_i$ for all $i \in L$. This is a missing face of $(\underline{\cS},\underline{\cT})^{* \cK}$ by definition of the polyhedral join. Moreover, it is minimal since for any $i \in L$ and $s \in J_i$
\[\bigsqcup_{i \neq k \in L} J_k \sqcup (J_i - \{s\}) = \bigsqcup_{k \in \tau} J_k \sqcup \sigma_i \in (\underline{\cS},\underline{\cT})^{* \tau} \]
where $\tau \in \cK$ since $L$ is a minimal missing face, and $\sigma_i \in \cT_i$ since $J_i$ is a minimal missing face.

Now we show that 
\[ 
\begin{aligned}
MF\left( (\underline{\cS},\underline{\cT})^{* \cK} \right) &\subseteq \{J \in MF(\cS_i) \; | \; i \in \cK \}\\ & \quad \sqcup \left\{ \bigsqcup_{i \in L} J_i \; | \; L \in MF(\cK), J_i \in MF(\cT_i), J_i \in \cS_i \right\}.
\end{aligned}
\]

Let $F \in MF((\underline{\cS},\underline{\cT})^{*\cK})$. We show that either $F \in \{J \in MF(\cS_i) \; | \; i \in \cK \}$ or $F \in \{\bigsqcup_{i \in L} F_i \; | \; L \in MF(\cK), F_i \in MF(\cT_i), F_i \in \cS_i \}$. 

Denote by $F_i = F|_{[l_i]}$ for each $i=1,\dots,m$ . If $F_i \in MF(\cS_i)$, then $F = F_i \in MF(\cK)$. Otherwise $F$ would not be minimal, since $F_i \subseteq F$ is a missing face. 

On the other hand, suppose that $F_i \in \cS_i$ for all $i$. Denote by $\sigma = \{i \in [m] \; | \; F|_i \neq \emptyset \}$. Firstly, $\sigma \notin \cK$, as otherwise $F \in \cK$ since $F = \bigsqcup_{i \in \sigma} F_i$, where $F_i \in \cS_i$ for all $i$. For all $i \in \sigma$, $F_i$ is a non-face of $\cT_i$. Otherwise, by removing vertices from $F_i$, we obtain a smaller non-face of $(\underline{\cS},\underline{\cT})^{*\cK}$. For such $i$, $F_i \in MF(\cT_i)$. It follows that $\sigma \in MF(\cK)$, as otherwise we restrict to a minimal missing face $\tau \in MF(\cK)$ with $\tau \subseteq \sigma$, and $\bigsqcup_{i \in \tau} F_i$ is a missing face of $(\underline{\cS},\underline{\cT})^{*\cK}$. Finally, for all $i$, $F_i$ is a minimal missing face of $\cT_i$. This follows by observing that otherwise for $\widehat{F_i} \subsetneq F_i$ with $\widehat{F_i} \in MF(\cS_i)$, $(F - F_i) \sqcup \widehat{F_i} \in (\underline{\cS},\underline{\cT})^{*\cK}$, which is a contradiction.
\end{proof}

\begin{proof}[Proof of Proposition~\ref{prop:MFsofKPi}]
By Proposition~\ref{prop:missingfacespolyhedraljoin}, the minimal missing faces of the composition complex $\cK(\cT_1,\dots,\cT_m)$ are
\[
MF(\cK(\cT_1,\dots,\cT_m)) = \left\{ \bigsqcup_{i \in L} J_i \mid L \in MF(\cK), J_i \in MF(\cT_i) \right\}
\]
so that
\begin{equation*} \pushQED{\qed}
    MF(\sk^{k-3} \Delta^{k-1} (\partial \Delta[I_1],\dots,\partial \Delta[I_k]))=\left\{ \{[m] \setminus \{i_1,\dots,i_{n_i}\}\} \mid i=1,\dots,k \right\}. \qedhere
    \popQED
\end{equation*}
\phantom \qedhere
\end{proof}

%\todo[inline,color=green]{Got to here 26/01/23}

Since a simplicial complex is determined by the set of minimal missing faces, $\cK_{\Pi}$ is the simplicial complex on $[m]$ given by
\[
MF(\cK_{\Pi}) = \{[m] \setminus P_i \mid i = 1,\dots,k\}.
\]
It follows that the Alexander dual of $\cK_{\Pi}$ is given by
\begin{equation} \label{eq:AlexDualofIdentityComplex}
\widehat{\cK_{\Pi}} = \bigsqcup_{i=1}^k \Delta[P_i].
\end{equation}

%\todo[inline]{Change $-$ for $\setminus$ throughout}

%\todo[inline,color=green]{Change `definintion' to `properties'}

We use this description of the Alexander dual to decompose $\cK_{\Pi}$ as the union of subcomplexes.

\begin{proposition} \label{prop:identitycomplexisacomposition}
Let $\Pi = \{P_1,\dots,P_k\}$ be a $k$-partition of $[m]$ and denote $P_i = \{i_1,\dots,i_{p_i}\}$ and $[m] \setminus P_i = \{j_1,\dots,j_{q_i}\}$. Then
\begin{equation} \label{eq:altdescofKpi}
    \cK_{\Pi} = \bigcup_{i=1}^k \partial \Delta \langle \partial \Delta[j_1,\dots,j_{q_i}],i_1,\dots,i_{p_i} \rangle.
\end{equation}
\end{proposition}

\begin{proof}
Two finite simplicial complexes are isomorphic if and only if their Alexander duals are isomorphic. By~\eqref{eq:AlexDualofIdentityComplex}, the Alexander dual of $\cK_{\Pi}$ has maximal faces $\bigsqcup_{i=1}^k \Delta[i_1,\dots,i_{p_i}]$. On the other hand, if $\cK^i = \partial \Delta \langle \partial \Delta[j_1,\dots,j_{q_i}],i_1,\dots,i_{p_i} \rangle$, then by Proposition~\ref{prop:missingfacespolyhedraljoin}, the maximal faces $\widehat{\cK^i}_{max}$ of $\widehat{\cK^i}$ are given by
\[
\{[m] \setminus \{j_1,\dots,j_{q_i}\} \sqcup \{[m] \setminus \{j,i_1,\dots,i_{p_i} \} \; | \; j \in Q_i \} = \{i_1,\dots,i_{p_i}\} \sqcup \{Q_i \setminus j \; | \; j \in Q_i \}
\]
where $Q_i = [m] \setminus P_i = \{j_1,\dots,j_{q_i}\}$. Then
\[
\widehat{\bigcup_{i=1}^k \cK^i} = \bigcap_{i=1}^k \widehat{\cK^i} = \bigcap_{i=1}^k \left(\widehat{\cK^i}_{max} \right) 
= \bigsqcup_{i=1}^k \Delta[i_1,\dots,i_{p_i}]
\]
and the result follows.
\end{proof}

%\todo[inline,color=green]{Do two triangles back to back to illustrate the 3 decompositions of $\cK_{\Pi}$}

We illustrate the three different decompositions of the complex $\cK_{\Pi}$.

\begin{example} \label{ex:FirstIdentityComplex}
    Let $\Pi = \{\{1\},\{2,3\},\{4\}\}$ be a $3$-partition of $[4]$. The complex $\cK_{\Pi} = \sk^0 \Delta^2(\circ[1],\partial \Delta [2,3],\circ[4])$ is given by
    \[
    \cK_{\Pi} = \{1\} \ast \partial \Delta[2,3] \cup \partial \Delta[2,3] \ast \{4\} \cup \Delta[2,3].
    \]
    This decomposition is shown in the following figure.
\begin{figure}[h]
\centering
\begin{tikzpicture}[scale=1.5]
    \coordinate [label=left:{$1$}] (1) at (0,0.5);
    \coordinate [label=below:{$3$}] (2) at (0.7,0);
    \coordinate [label=above:{$2$}] (3) at (0.7,1);
    \coordinate [label=right:{$4$}] (4) at (1.4,0.5);
    \coordinate [label={$=$}] (5) at (2.1,0.35);
    \draw (3) -- (1) -- (2) -- (3) -- (4) -- (2);
    \foreach \point in {1,2,3,4}
        \fill [black] (\point) circle (1 pt);
\end{tikzpicture} \hspace{0.025 \textwidth}
\begin{tikzpicture}[scale=1.5]
    \coordinate [label=left:{$1$}] (1) at (0,0.5);
    \coordinate [label=below:{$3$}] (2) at (0.7,0);
    \coordinate [label=above:{$2$}] (3) at (0.7,1);
    %\coordinate [label={$\cup$}] (4) at (1.4,0.35);
    \draw (2) -- (1) -- (3);
    \foreach \point in {1,2,3}
        \fill [black] (\point) circle (1 pt);
\end{tikzpicture} \hspace{0.05 \textwidth}
\begin{tikzpicture}[scale=1.5]
    \coordinate [label=below:{$3$}] (2) at (0.7,0);
    \coordinate [label=above:{$2$}] (3) at (0.7,1);
    \draw (2) -- (3);
    \foreach \point in {2,3}
        \fill [black] (\point) circle (1 pt);
\end{tikzpicture} \hspace{0.05 \textwidth}
\begin{tikzpicture}[scale=1.5]
    %\coordinate [label={$\cup$}] (1) at (0,0.35);
    \coordinate [label=below:{$3$}] (2) at (0.7,0);
    \coordinate [label=above:{$2$}] (3) at (0.7,1);
    \coordinate [label=right:{$4$}] (4) at (1.4,0.5);
    \draw (2) -- (4) -- (3);
    \foreach \point in {2,3,4}
        \fill [black] (\point) circle (1 pt);
\end{tikzpicture}
\caption{Decomposition of the identity complex $\cK_{\{\{1\},\{2,3\},\{4\}\}}$.}
\end{figure}
%\todo[inline,color = green]{Caption}
    The minimal missing faces of $\cK_{\Pi}$ are $Q_1 = \{2,3,4\}$, $Q_2 = \{1,4\}$, and $Q_3 = \{1,2,3\}$. Finally, since the complex $\cK_{\Pi}$ is also given by $\partial \Delta \langle \partial \Delta[1,4],2,3 \rangle$, we obtain the decomposition
    \[
    \cK_{\Pi} = \partial \Delta \langle \partial \Delta[2,3,4] , 1 \rangle \cup \partial \Delta \langle \partial \Delta[1,4],2,3 \rangle \cup \partial \Delta \langle \partial \Delta[1,2,3],4 \rangle.
    \]
\end{example}

\subsection{Relations among higher Whitehead maps}

Let $\Pi = \{P_1,\dots,P_k\}$ be a $k$-partition of $[m]$ and let $f_i \colon \Sigma X_i \longrightarrow Y_i$ be maps for $i=1,\dots,m$. By Proposition~\ref{prop:identitycomplexisacomposition}, $\cK_{\Pi} = \bigcup_{i=1}^k \partial \Delta \langle \partial \Delta[j_1,\dots,j_{q_i}],i_1,\dots,i_{p_i} \rangle$. Therefore, the space $(\underline{Y},\underline{\ast})^{\cK_{\Pi}}$ contains the codomains of the higher Whitehead maps
\begin{align} \label{eq:hwMapsInRelation}
& h_w \left( h_w \left( f_{j_1},\dots,f_{j_{q_i}} \right) ,f_{i_1},\dots,f_{i_{p_i}} \right) \colon\\
& \qquad \qquad \qquad \Sigma^{p_i} \left( \Sigma^{q_i-2} X_{j_1} \wedge \cdots \wedge X_{j_{q_i}} \right) \wedge X_{i_1} \wedge \cdots \wedge X_{i_{p_i}} \longrightarrow (\underline{Y},\underline{\ast})^{\cK^i} \nonumber
\end{align}
for $i = 1,\dots,k$. 

We generalise the relation of Hardie~\eqref{eq:HardieIdentity2} to relations among the higher Whitehead maps~\eqref{eq:hwMapsInRelation} in polyhedral products over the identity complex $\cK_{\Pi}$.

%\todo[inline]{Change $i_{n_i}$ to $i_{p_i}$ and $j_{r_i}$ to $j_{q_i}$}

\begin{theorem} \label{thm:maintheorem}
For $i=1,\dots,m$, let $f_i \colon \Sigma X_i \longrightarrow Y_i$ be maps. Let $\Pi = \{P_1,\dots,P_k\}$ be a $k$-partition of $[m]$ for $k \geqslant 3$ and denote by  $P_i = \{i_1,\dots,i_{p_i}\}$ and $Q_i = [m] \setminus P_i = \{j_1,\dots,j_{q_i}\}$. If $X_i$ is a suspension for each $i=1,\dots,m$, then
\begin{equation} \label{eq:maintheorem}
\sum_{i=1}^k h_w^{\cK_{\Pi}} \left( h_w \left( f_{j_1},\dots,f_{j_{q_i}} \right),f_{i_1},\dots,f_{i_{p_i}} \right) \circ \sigma_i = 0
\end{equation}
in $\left[ \Sigma^{m-2} X_1 \wedge \cdots \wedge X_m, (\underline{Y},\underline{\ast})^{\cK_{\Pi}}\right]$, where
\[
\sigma_i \colon \Sigma^{m-2} X_1 \wedge \cdots \wedge X_m \longrightarrow \Sigma^{p_i} (\Sigma^{q_i-2} X_{j_1} \wedge \cdots \wedge X_{j_{q_i}}) \wedge X_{i_1} \wedge \cdots \wedge X_{i_{p_i}}
\]
is induced by the coordinate permutation
\[
X_1 \times \cdots \times X_m \longrightarrow X_{j_1} \times \cdots \times X_{j_{q_i}} \times X_{i_1} \times \cdots \times X_{i_{p_i}}.
\]
\end{theorem}

%\todo[color=green, inline]{GS: We're asked to clarify how $\Sigma^{m-2} X_1 \wedge \cdots \wedge X_m$ sits inside $CX_1 \times \cdots \times CX_m$. This is also a similar problem to clarifying how \eqref{eq:hwcofibdiag} follows from Porter. George to figure out how $(\underline{D},\underline{S})^{\partial \Delta^1}$. $(t_1,x_1),...,(t_m, x_m)$. Distribute suspension across first number of smashes. See if Porter/Selick/Whitehead's yellow book covers this? MS: I've changed it so that it just says that it's induced by a coordinate permutation (which it is).}

We delay the proof of Theorem~\ref{thm:maintheorem} to Section~\ref{sec:newmainproof}.
Now, we illustrate how Theorem~\ref{thm:maintheorem} generalises Hardie's identity to higher Whitehead maps and obtain examples of new relations among higher Whitehead maps. 

\begin{example} \label{ex:MainIdentityExamples}
For $i=1,\dots,m$, let $f_i \colon \Sigma X_i \longrightarrow Y_i$ be maps such that $X_i$ is a suspension. 
\begin{enumerate}[(i)]
\item %1
Let $\Pi = \{\{1\}, \{2\}, \{3\}\}$. Then $Q_1 = \{2,3\}$, $Q_2 = \{1,3\}$ and $Q_3 = \{1,2\}$, and $\cK_{\Pi}= \bullet_{[3]}$. We obtain the relation
\[
h_w^{\cK_{\Pi}} (h_w(f_2,f_3),f_1) \circ \sigma_1 + h_w^{\cK_{\Pi}} (h_w(f_1,f_3),f_2) \circ \sigma_2 + h_w^{\cK_{\Pi}} (h_w(f_1,f_2),f_3) \circ \sigma_3 
= 0
\]
in $\left[\Sigma X_1 \wedge X_2 \wedge X_3, Y_1 \vee Y_2 \vee Y_3\right]$. This recovers the generalised Jacobi identity, see~\cite{Arkowitz1962paper}.
\item %2
Generalising the previous example, let $\Pi = \{\{1\} , \dots ,\{m\}\}$. This implies that $Q_i = \{1,\dots,\hat{i},\dots,m\}$ for $i = 1,\dots,m$ and $\cK_{\Pi} = \sk^{m-3} \Delta^{m-1}$, giving rise to the relation
\[
h_w^{\cK_{\Pi}} \left(h_w \left(f_2,\dots,f_m \right), f_1 \right) \circ \sigma_1 + \cdots + h_w^{\cK_{\Pi}} \left(h_w \left(f_1,\dots,f_{m-1} \right), f_m \right) \circ \sigma_m = 0
\]
in $\left[\Sigma^{m-2} X_1 \wedge \cdots \wedge X_m, (\underline{Y},\underline{\ast})^{\sk^{m-3} \Delta^{m-1}} \right]$, generalising Hardie's identity~\eqref{eq:HardieIdentity2}.
\item %3 
Define a $3$-partition $\Pi = \{\{1\},\{2,3\},\{4\}\}$ of the vertex set $[4]$. In this case $Q_1 = \{2,3,4\}$, $Q_2 = \{1,4\}$ and $Q_3 = \{1,2,3\}$ and $\cK_{\Pi}$ is the simplicial complex shown in Example~\ref{ex:FirstIdentityComplex}. Then
\begin{align*}
h_w^{\cK_{\Pi}} (h_w(f_2,f_3,f_4),f_1)  \circ \sigma_1 + h_w^{\cK_{\Pi}} & (h_w(f_1,f_4),f_2,f_3) \circ \sigma_2 \\ &+ h_w^{\cK_{\Pi}} (h_w(f_1,f_2,f_3),f_4) \circ \sigma_3
=0
\end{align*}
in $\left[\Sigma^2 X_1 \wedge X_2 \wedge X_3 \wedge X_4, (\underline{Y},\underline{\ast})^{\cK_{\Pi}} \right]$. In contrast to Hardie's identity~\eqref{eq:HardieIdentity2}, the summands of this relation include both $2$-ary and $3$-ary higher Whitehead maps.
\end{enumerate}
\end{example}

When the maps $f_i$ are spherical, the permutation maps $\sigma_i$ are degree maps. Recall the definition of the Koszul sign $\epsilon(\sigma)$ of a permutation $\sigma$ as given in Proposition~\ref{prop:hwSymmetry}.

\begin{corollary} \label{cor:MainTheoremSpherical}
Let $f_i \in \pi_{r_i}(Y_i)$ for $i=1,\dots,m$. Let $\Pi = \{P_1,\dots,P_k\}$ be a $k$-partition of $[m]$ for $k \geqslant 3$ and denote $P_i = \{i_1,\dots,i_{p_i}\}$ and $J_i = [m] \setminus P_i = \{j_1,\dots,j_{q_i}\}$. If $r_i \geqslant 2$ for each $i=1,\dots,m$, then there is a relation
\[
\sum_{i=1}^k \epsilon(\sigma_i) h_w^{\cK_\Pi} \left( h_w \left(f_{j_1},\dots,f_{j_{q_i}} \right), f_{i_1},\dots,f_{i_{p_i}} \right) = 0
\]
in $\pi_{r_1+\cdots+r_m-2} \left( (\underline{Y},\underline{\ast})^{\cK_{\Pi}} \right)$, where $\epsilon(\sigma_i)$ is the Koszul sign of the permutation\\ $\sigma_i \colon (1,\dots,m) \mapsto (j_1,\dots,j_{q_i},i_1,\dots,i_{p_i})$.
\end{corollary}

\begin{proof}
The map $\sigma_i \colon S^{r_1+\cdots+r_m-2} \longrightarrow S^{r_1+\cdots+r_m-2}$ from the statement of Theorem~\ref{thm:maintheorem} has degree $\epsilon(\sigma_i)$ and the result follows.
\end{proof}

\goodbreak

\begin{example}
%\todo[inline,color=green]{Just do Hardie, and then say it reduces to Jacobi}
In some cases, $\epsilon(\rho)$ can be calculated explicitly. We revisit the relations from Example~\ref{ex:MainIdentityExamples} when $f_i \colon S^{r_i} \longrightarrow Y_i$, with $r_i \geqslant 2$ for $i=1,\dots,m$.
\begin{enumerate}[(i)]
    \item Let $\Pi = \{\{1\},\dots,\{m\}\}$. Denote by $\iota_i$ the inclusion $Y_i \longrightarrow (\underline{Y},\underline{\ast})^{\cK_{\Pi}}$ and $\kappa_i$ the inclusion $FW(Y_1,\dots,Y_{i-1},Y_{i+1},\dots,Y_m) \longrightarrow (\underline{Y},\underline{\ast})^{\cK_{\Pi}}$. We recover relation~\eqref{eq:HardieIdentity2} of Hardie,
    \[
    \sum_{i=1}^m (-1)^{p_i(p_{i+1}+\cdots+p_m)} [\theta_i h_w(f_1,\dots,f_{i-1},f_{i+1},\dots,f_m),\psi_i f_i] = 0
    \]
    in $\pi_{p_1 + \cdots + p_m - 2} \left( (\underline{Y},\underline{\ast})^{\cK_{\Pi}} \right)$. Moreover, if $k=3$, let $g_j = \iota_j \circ f_j$ for $j=1,2,3$. Multiplying by $(-1)^{p_1 p_3}$ and applying Proposition~\ref{prop:hwisHWPBasic}, we recover the graded Jacobi identity for Whitehead products
    \[
    (-1)^{p_1 p_2} [[g_2,g_3],g_1] + (-1)^{p_2 p_3} [[g_1,g_3],g_2] + (-1)^{p_1 p_3} [[g_1,g_2],g_3] = 0.
    \]
    \item for $\Pi = \{\{1\},\{2,3\},\{4\}\}$, there is a relation
    \begin{align*}
    (-1)^{p_1(p_2+p_3)} h_w^{\cK_{\Pi}} (h_w(f_2,f_3,f_4),f_1)  + & (-1)^{(p_2+p_3)p_4} h_w^{\cK_{\Pi}}  (h_w(f_4,f_1),f_2,f_3) \\ &+ (-1)^{p_1 p_4} h_w^{\cK_{\Pi}} (h_w(f_1,f_2,f_3),f_4)
    =0
\end{align*}
in $\pi_{p_1+\cdots+p_4-2}\left( (\underline{Y},\underline{\ast})^{\cK_{\Pi}} \right)$.
\end{enumerate}
\end{example}

As a further consequence of Theorem~\ref{thm:maintheorem}, we expand the family of simplicial complexes in which we can identify relations among higher Whitehead maps.

\begin{corollary} \label{cor:MainThmSubstituted}
    For $i=1,\dots,m$, let $\cS_i$ be a simplicial complex on $[l_i]$ and let $f_i \colon \Sigma X_i \longrightarrow (\underline{Y},\underline{\ast})^{\cS_i}$ be maps. Let $\Pi = \{P_1,\dots,P_k\}$ be a $k$-partition of $[m]$ for $k \geqslant 3$ and denote by $P_i = \{i_1,\dots,i_{p_i}\}$ and $Q_i = [m] \setminus P_i = \{j_1,\dots,j_{q_i}\}$. Suppose that $\cL$ is a simplicial complex containing $\cK_{\Pi} \langle \cS_1,\dots,\cS_m \rangle$. Then if $X_i$ is a suspension for each $i=1,\dots,m$,
\begin{equation*}
\sum_{i=1}^k h_w^{\cL} \left( h_w \left( f_{j_1},\dots,f_{j_{q_i}} \right),f_{i_1},\dots,f_{i_{p_i}} \right) \circ \sigma_i = 0
\end{equation*}
in $\left[ \Sigma^{m-2} X_1 \wedge \cdots \wedge X_m, (\underline{Y},\underline{\ast})^{\cL} \right]$, where
\[
\sigma_i \colon \Sigma^{m-2} X_1 \wedge \cdots \wedge X_m \longrightarrow \Sigma^{p_i} (\Sigma^{q_i-2} X_{j_1} \wedge \cdots \wedge X_{j_{q_i}}) \wedge X_{i_1} \wedge \cdots \wedge X_{i_{p_i}}
\]
is induced by the coordinate permutation
\[
X_1 \times \cdots \times X_m \longrightarrow X_{j_1} \times \cdots \times X_{j_{q_i}} \times X_{i_1} \times \cdots \times X_{i_{p_i}}.
\]
\end{corollary}

%\todo[inline,color=green]{Corollary of main theorem describing propagation. Increase the family of simplicial complexes in which we can identity relations among higher Whitehead maps}

\begin{example}
%\todo[inline,color=green]{We answer whether the relations are non-trivial/have non-trivial elements}
In some cases, the summands in relation~\eqref{eq:maintheorem} can be shown to be non-trivial. 
Suppose that for $i=1,\dots,m$, $f_i \colon S^2 \longrightarrow \CP$ is the inclusion of the bottom cell. Let $\Pi$ be a $k$-partition of $[m]$. By Corollary~\ref{cor:MainTheoremSpherical}
\begin{equation} \label{eq:RelationForPanov}
\sum_{i=1}^k \epsilon(\rho) h_w^{\cK_\Pi} \left( h_w \left(f_{j_1},\dots,f_{j_{q_i}} \right), f_{i_1},\dots,f_{i_{p_i}} \right) = 0.
\end{equation}

We show that none of the maps in~\eqref{eq:RelationForPanov} are null-homotopic. Abramyan and Panov showed in \cite[Theorem~5.2]{AbramyanPanov} that the map $h_w^{\cK_{\Pi}} \left( h_w \left( f_{j_1},\dots,f_{j_{q_i}} \right),f_{i_1},\dots,f_{i_{p_i}} \right)$ is null-homotopic if and only if $\Delta \langle \partial \Delta[j_1,\dots,j_{q_i}],i_1,\dots,i_{p_i} \rangle \subseteq \cK_{\Pi}$. The complex $\Delta \langle \partial \Delta[j_1,\dots,j_{q_i}],i_1,\dots,i_{p_i} \rangle \subseteq \cK_{\Pi}$ contains the simplex $Q_i = [m] \setminus P_i$. On the other hand, by the definition of the complex $\cK_{\Pi}$ associated to the partition $\Pi$, $Q_i$ is a missing face of $\cK_{\Pi}$. Therefore, $\Delta \langle \partial \Delta[j_1,\dots,j_{q_i}],i_1,\dots,i_{p_i} \rangle$ is not a subcomplex of $\cK_{\Pi}$. Hence each term in the relation
\[
\sum_{i=1}^k \epsilon(\rho) h_w^{\cK_\Pi} \left( h_w \left(f_{j_1},\dots,f_{j_{q_i}} \right), f_{i_1},\dots,f_{i_{p_i}} \right) = 0
\]
is non-trivial.
\end{example}

\section{Relations among folded higher Whitehead maps}

Folded higher Whitehead maps provide a combinatorial framework to study elements of higher Whitehead products with repeated factors. By propagating the relations among higher Whitehead maps given in Theorem~\ref{thm:maintheorem}, we study relations among folded higher Whitehead maps.

Relations among folded higher Whitehead maps have been studied rationally by Zhuravleva~\cite{Zhur21}. Let $\cK = \partial \Delta \langle \partial \Delta \langle 1,2,3 \rangle, 4,5 \rangle$ and for $i=1,\dots,m$, let $u_i \colon S^1 \longrightarrow \Omega \CP$ be the adjoint of the inclusion $\mu_i \colon S^2 \longrightarrow \CP$ of the bottom cell. Zhuravleva showed that the relation
\begin{equation} \label{eq:ElizavetaRelation}
[[u_1,u_2,u_3],[u_1,u_4,u_5]] + [[[u_1,u_2,u_3],u_4,u_5],u_1] = 0
\end{equation}
holds in $\pi_8 (\Omega \DJK) \otimes \mathbb{Q}$. We improve this result by showing that relation~\eqref{eq:ElizavetaRelation} holds integrally in $\pi_*(\Omega \DJK)$.

\subsection{Folds of identity complexes}

For $i=1,\dots,m$, let $f_i \colon \Sigma X_i \longrightarrow Y_i$ be maps. Let $\Pi = \{P_1,\dots,P_k\}$ be a $k$-partition of $[m]$, and denote by $P_i = \{i_1,\dots,i_{p_i}\}$, and $[m] \setminus P_i = \{j_1,\dots,j_{q_i}\}$. By Theorem~\ref{thm:maintheorem},
\begin{equation} \label{eq:MainThmForFold}
\sum_{i=1}^k h_w^{\cK_{\Pi}} \left( h_w \left( f_{j_1},\dots,f_{j_{q_i}} \right), f_{i_1},\dots,f_{i_{p_i}} \right) \circ \sigma_i = 0 
\end{equation}
in $\left[ \Sigma^{m-2} X_1 \wedge \cdots \wedge X_m, (\underline{Y},\underline{\ast})^{\cK_{\Pi}}\right]$. Let $\psi \colon I \longrightarrow J$ be a fold of $\cK_{\Pi}$ and suppose that $Y_j$ is an associative $H$-space for each $j \in J$ , and that $Y_i = Y_j$ for every $i \in I_j = \psi^{-1} (j)$. By composing each summand in~\eqref{eq:MainThmForFold} with the fold $
\nabla_{(I,J)} \colon (\underline{Y},\underline{\ast})^{\cK_{\Pi}} \longrightarrow (\underline{Y},\underline{\ast})^{(\cK_{\Pi})_{\nabla(I,J)}}$,
we obtain the relation 
\begin{equation} \label{eq:MainThmFolded}
    \sum_{i=1}^k \nabla_{(I,J)}  h_w^{\cK_{\Pi}} \left( h_w \left( f_{j_1},\dots,f_{j_{q_i}} \right), f_{i_1},\dots,f_{i_{p_i}} \right) \circ \sigma_i = 0
\end{equation}
among folded higher Whitehead maps in $\left[ \Sigma^{m-2} X_1 \wedge \cdots \wedge X_m, (\underline{Y},\underline{\ast})^{(\cK_{\Pi})_{\nabla(I,J)}} \right]$. We analyse the combinatorial properties of the folded complex $(\cK_{\Pi})_{\nabla(I,J)}$ for various partitions $\Pi$ and folds $\psi \colon I \longrightarrow J$, which govern the form of relation~\eqref{eq:MainThmFolded}.

%\todo[inline,color=green]{Linking sentence: We analyse the different type of relations obtained from different partitions $\Pi$ and folds $\psi \colon I \longrightarrow J$.}

%We start by considering the case that $\cK_i = \bullet$ for each $i=1,\dots,m$, and give a full characterisation of fold maps $\psi \colon I \longrightarrow J$ which produce relations among folded higher Whitehead maps. These relations detect interesting structures in homotopy groups, for example the existence of $2$-torsion.

We determine, in terms of the partition $\Pi$ and the fold $\psi \colon I \longrightarrow J$, which terms of relation~\eqref{eq:MainThmFolded} are null-homotopic after folding. By Proposition~\ref{prop:hwFoldedTrivialityConditions}, the map $\nabla_{(I,J)}  h_w^{\cK_{\Pi}} \left( h_w \left(f_{j_1},\dots,f_{j_{q_i}} \right),f_{i_1},\dots,f_{i_{p_i}} \right)$
is null-homotopic if\\ $h_w^{\cL} \left( h_w \left(f_{j_1},\dots,f_{j_{q_i}} \right),f_{i_1},\dots,f_{i_{p_i}} \right)$ is null-homotopic, where
\begin{equation} \label{eq:DefLRecall}
\cL = (\cK_{\Pi})_{\nabla(I,J)} \langle \Delta[\{k_1\} \sqcup I_{k_1}],\dots,\Delta[\{k_l\} \sqcup I_{k_l}] \rangle
\end{equation}
where $I_k = \emptyset$ if $k \notin J$, and $\{k_1,\dots,k_l\}$ is the vertex set of $(\cK_{\Pi})_{\nabla(I,J)}$.

By first describing the folded identity complex $(\cK_{\Pi})_{\nabla(I,J)}$ in terms of the partition $\Pi$ and the fold $\psi \colon I \longrightarrow J$, we are able to express the complex $\cL$ associated to $\cK_{\Pi}$ in terms of $\Pi$ and $\psi$. Applying Proposition~\ref{prop:hwFoldedTrivialityConditions}, we identify the null-homotopic terms in relation~\eqref{eq:MainThmFolded}.

\begin{lemma} \label{lm:IdentityComplexFold}
Let $\Pi = \{P_1,\dots,P_k\}$ be a $k$-partition of $[m]$. Then the following hold:
\begin{enumerate}[(i)]
    \item if $I \cup J \subseteq P_l$ for some $l \in \{1,\dots,k \}$, then
    \[
    (\cK_{\Pi})_{\nabla(I,J)} = \partial \Delta[[m] \setminus P_l] \ast \Delta[P_l \setminus I];
    \]
    \item if $|I| = |J| = 1$ and $I \subseteq P_i$ and $J \subseteq P_j$ for $i \ne j$, then
    \[
    (\cK_{\Pi})_{\nabla(I,J)} = \partial \Delta[[m] \setminus I];
    \]
    \item otherwise, $(\cK_{\Pi})_{\nabla(I,J)} = \Delta[[m] \setminus I]$.
\end{enumerate}
\end{lemma}

\begin{proof}
We first prove (ii). Let $(\underline{\Delta},\underline{\partial \Delta})$ denote the tuple of pairs\\ $((\Delta[P_1],\partial \Delta[P_1]),\dots,(\Delta[P_k],\partial \Delta[P_k]))$. By definition, $\cK_{\Pi} = (\underline{\Delta},\underline{\partial \Delta})^{*sk^{k-3} \Delta[1,\dots,k]}$. For $i \ne j$, this polyhedral join decomposes into the union
\begin{align*}
\cK_{\Pi} &= \Delta[P_i] \ast \Delta[P_j] \ast \left( (\underline{\Delta},\underline{\partial \Delta})^{*sk^{k-5} \Delta[[k] \setminus \{i,j\}]}  \right) \\
& \qquad \cup (\Delta[P_i] \ast \partial \Delta[P_j] \cup \partial \Delta[P_i] \ast \Delta[P_j]) \ast \left( (\underline{\Delta},\underline{\partial \Delta})^{* \partial \Delta[[k] \setminus \{i,j\}]}  \right) \\
& \qquad \cup \partial \Delta[P_i] \ast \partial \Delta[P_j] \ast \left( (\underline{\Delta},\underline{\partial \Delta})^{*\Delta[[k] \setminus \{i,j\}]}  \right).
\end{align*}
Suppose that $I = \{u\}$ and $J = \{v\}$ such that $u \in P_i$ and $v \in P_j$. Then
\begin{align*}
    &\nabla_{(u,v)}(\Delta[P_i] \ast \Delta[P_j]) = \Delta[P_i \sqcup P_j \setminus \{u\}] \\
    &\nabla_{(u,v)}(\Delta[P_i] \ast \partial \Delta[P_j])  = \Delta[P_i \sqcup P_j \setminus \{u\}] \\
    &\nabla_{(u,v)}(\partial \Delta[P_i] \ast \Delta[P_j])  = \Delta[P_i \sqcup P_j \setminus \{u\}] \\
    &\nabla_{(u,v)}(\partial \Delta[P_i] \ast \partial \Delta[P_j])  = \partial \Delta[P_i \sqcup P_j \setminus \{u\}].
\end{align*}
Since $sk^{k-5} \Delta[[k] \setminus \{i,j\}]$ is a subcomplex of $\partial \Delta[[k] \setminus \{i,j\}]$, it follows that
\begin{align*}
    (\cK_{\Pi})_{\nabla(I,J)} &= \Delta[P_i \sqcup P_j \setminus \{u\}] \ast \left( (\underline{\Delta},\underline{\partial \Delta})^{* \partial \Delta[[k] \setminus \{i,j\}]}  \right) \\
    & \qquad \cup \partial \Delta[P_i \sqcup P_j \setminus \{u\}] \ast \left( (\underline{\Delta},\underline{\partial \Delta})^{*\Delta[[k] \setminus \{i,j\}]}  \right) \\
    &= ((\Delta[P_1],\partial \Delta[P_1]),\dots,(\Delta[P_i \sqcup P_j \setminus \{u\}], \partial \Delta[P_i \sqcup P_j \setminus \{u\}]),\dots \\
    & \hspace{10em} \dots,(\Delta[P_k],\partial \Delta[P_k]))^{\ast \partial \Delta[[k] \setminus i]} \\
    &= \partial \Delta[[m] \setminus \{u\}]
\end{align*}
proving (ii).

%\todo[inline,color=green]{Fill in $\cdots$, probably mathematically, or in words. Fix conclusion.}

We now establish (i) and (iii). Assume that $u,v \in P_l$. The folding sends $\Delta[P_l] \mapsto \Delta[P_l\setminus\{u\}]$ and $\partial \Delta[P_l] \mapsto \Delta[P_l\setminus\{u\}]$.
%\todo[inline,color=green]{Add maths to justify}
Then since
\[
\cK_{\Pi} = \Delta[P_l] \ast \left( (\underline{\Delta},\underline{\partial \Delta})^{* sk^{k-4} \Delta[[k] \setminus \{l\}]} \right) \cup \partial \Delta[P_l] \ast \left( (\underline{\Delta},\underline{\partial \Delta})^{* \partial \Delta[[k] \setminus \{l\}]} \right)
\]
we obtain 
\[
    (\cK_{\Pi})_{\nabla(I,J)} = \Delta[P_l\setminus\{u\}] \ast \left( (\underline{\Delta},\underline{\partial \Delta})^{* \partial \Delta[[k] \setminus \{l\}]} \right) = \Delta[P_l \setminus \{u\}] \ast \partial \Delta[[m] \setminus P_l].
\]
Next, for any $I$ and $J$ with $u \in I$ and $v \in J$, by Proposition~\ref{prop:propsoffoldofsimpcmplx},
\[
(\cK_{\Pi})_{\nabla(I,J)} = \left( (\cK_{\Pi})_{\nabla(u,v)} \right)_{\nabla(I \setminus \{u\},J \setminus \{v\})} = \left( \Delta[P_l \setminus \{u\}] \ast \partial \Delta[[m] \setminus P_l] \right)_{\nabla(I \setminus \{u\},J \setminus \{v\})} 
\]
Claim (i) then follows since for $I,J \subseteq P_l$, $\left( \Delta[P_l \setminus \{u\}] \ast \partial \Delta[[m] \setminus P_l] \right)_{\nabla(I \setminus \{u\},J \setminus \{v\})} = \Delta[P_l \setminus I] \ast \partial \Delta[[m] \setminus P_l]$.

Alternatively, claim (iii) follows since for any fold not satisfying the hypotheses of (i) or (ii), there is $u' \in I$ with $u' \in P_{l'}$ for $l' \ne l$. Then\\ $\left( \Delta[P_l \setminus \{u\}] \ast \partial \Delta[[m] \setminus P_l] \right)_{\nabla(I \setminus \{u\},J \setminus \{v\})} = \Delta[[m] \setminus I]$ since any fold of $\partial \Delta[[m] \setminus P_l]$ is a simplex.
\end{proof}

\begin{lemma} \label{lm:IdentityComplexUnfolds}
Let $\cL$ denote the complex in~\eqref{eq:DefLRecall}. Then the following hold:
\begin{enumerate}[(i)]
    \item if $I \cup J \subseteq P_l$ for some $l \in \{1,\dots,k\}$, then
    \[
    \cL = \partial \Delta[[m] \setminus P_l] \ast \Delta[P_l];
    \]
    \item if $|I|=|J|=1$ and $I \subseteq P_i$ and $J \subseteq P_j$ for $i \ne j$, then $\cL$ has minimal missing faces $[m] \setminus I$ and $[m] \setminus J$;
    \item otherwise, $\cL = \Delta[[m]]$.
\end{enumerate}
\end{lemma}

\begin{proof}
%\todo[inline,color=red]{Change to eqref from updated Prop -- We now need the proposition and its proof.}
By Proposition~\ref{prop:missingfacespolyhedraljoin},
\[
MF(\cL) = \left\{ \bigsqcup_{i \in \kappa} \tau_i \mid \kappa \in MF \left( (\cK_{\Pi})_{\nabla(I,J)}  \right), \tau_i \in \Delta[\{k_i\} \sqcup I_{k_i}] \right\}.
\]
Statement (i) then follows from Lemma~\ref{lm:IdentityComplexFold} since the missing face of $\partial \Delta[[m] \setminus P_l] \ast \Delta[P_l \setminus I]$ is not supported by any vertex in $J$. Similarly, statements (ii) and (iii) also follow from Lemma~\ref{lm:IdentityComplexFold}.
\end{proof}

%\todo[inline,color=green]{Fix this in English}

Therefore the only folds of $\cK_{\Pi}$ that do not trivialise the map\\ $\nabla_{(I,J)}  h_w^{\cK_{\Pi}} \left( h_w \left(f_{j_1},\dots,f_{j_{q_i}} \right), f_{i_1},\dots,f_{i_{p_i}} \right)$ are those identifying some $i_s$ with some $j_r$, or vice-versa. The relations arising from folds of this form are described in the following theorem.

\begin{theorem} \label{thm:FoldsofIdentityComplexes}
Let $\Pi = \{P_1,\dots,P_k\}$ be a $k$-partition of $[m]$ and let $f_i \colon \Sigma X_i \longrightarrow Y_i$ be maps for $i=1,\dots,m$, where each $X_i$ is a suspension. Fix $u,v \in [k]$ with $u \ne v$ and write $P_u = \{i_1,\dots,i_{p_i}\}$ and $P_v = \{j_1,\dots,j_{p_j}\}$. Let $i \in P_u$ and $j \in P_v$ and let $\psi \colon \{i\} \longrightarrow \{j\}$ be a fold of $\cK_{\Pi}$. If $Y_i = Y_j$ is an $H$-space, then there is a relation
\begin{align*}
& \nabla_{(i,j)}  h_w^{\cK_{\Pi}}(h_w(f_{i'_1},\dots,f_{i'_{q_i}}),f_{i_1},\dots,f_{i_{p_i}}) \circ \sigma_i \\
& \qquad + \nabla_{(i,j)}  h_w^{\cK_{\Pi}}(h_w(f_{j'_1},\dots,f_{j'_{q_j}}),f_{j_1},\dots,f_{j_{p_j}}) \circ \sigma_j = 0
\end{align*}
in $ \left[ \Sigma^{m-2} X_1 \wedge \cdots \wedge X_m, (\underline{Y},\underline{\ast})^{\partial \Delta[1,\dots,i-1,i+1,\dots,m]} \right]$, where $[m]  \setminus  P_u = \{i'_1,\dots,i'_{q_i}\}$ and $[m] \setminus  P_v = \{j'_1,\dots,j'_{q_j}\}$.
\end{theorem}

\begin{proof}
We show that if either $I$ or $J$ contain more than one element, or if $I$ and $J$ are singletons both contained in one of $P_i = \{i_1,\dots,i_{p_i}\}$ or $[m]  \setminus  P_i = \{j_1,\dots,j_{q_i}\}$, then the folded higher Whitehead map
\[
H_i = \nabla_{(I,J)}  h_w^{\cK_{\Pi}} \left( h_w \left(f_{j_1},\dots,f_{j_{q_i}} \right), f_{i_1},\dots,f_{i_{p_i}} \right)
\]
is null-homotopic.

Suppose that at least one of $I$ and $J$ contains more than one element. Firstly, if $I \sqcup J$ is not contained in some $P_r$ for $r \in \{1,\dots,k\}$, then by Lemma~\ref{lm:IdentityComplexUnfolds}(iii), $\cL = \Delta[m]$. By Proposition~\ref{prop:hwFoldedTrivialityConditions}, it follows that the map $H_i$ is null-homotopic for each $i=1,\dots,k$.

Secondly, if $I \sqcup J \subseteq P_r$ for some $r=1,\dots,k$, then by Lemma~\ref{lm:IdentityComplexUnfolds}(i), $\cL = \partial \Delta[[m] \setminus P_r] \ast \Delta[P_r]$. Then $H_r$ is null-homotopic, and furthermore for $i \ne r$, since $\partial \Delta \langle \Delta[j_1,\dots,j_{q_i}],i_1,\dots,i_{p_i} \rangle \subseteq \cL$, the map $H_i$ is also null-homotopic for $i \ne r$.

Now suppose that $|I| = |J| = 1$, and that $I = \{u\} \in P_s$ and $J = \{v\} \in P_t$ for $s \ne t$. By Lemma~\ref{lm:IdentityComplexUnfolds}(ii), $\cL$ has minimal missing faces $[m] \setminus \{u\}$ and $[m] \setminus \{v\}$. Therefore if $u = j_l$ and $v = j_k$, then $\partial \Delta \langle \Delta[j_1,\dots,j_{q_i}],i_1,\dots,i_{p_i} \rangle \subseteq \cL$, and hence $H_i$ is null-homotopic. Observe that since $s \ne t$, it is not possible for $u = i_l$ and $v = i_k$ in this case. 

The only remaining possibility is that $|I| = |J| = 1$ and that $I \subseteq P_u$ and $J \subseteq P_v$ for $u \ne v$, and so the theorem follows.
\end{proof}

%\todo[inline,color=green]{State theorem, and use the proof of 6.3. Remove 6.3.}

%It follows from Proposition~\ref{prop:FoldedIdentityMaps} that combinatorial conditions imply that every summand of relation~\eqref{eq:MainThmFolded} is null-homotopic unless $|I|=|J|=1$ with $I \subseteq P_i$ and $J \subseteq P_j$ for $i \neq j$, in which case we obtain a new relation among folded higher Whitehead maps.

A special case of Theorem~\ref{thm:FoldsofIdentityComplexes} gives an example of when a folded $n$-ary higher Whitehead map is identified with a folded $2$-ary higher Whitehead map.

\begin{corollary} \label{cor:HigherToTwo}
Let $f_i \colon \Sigma X_i \longrightarrow Y_i$ be maps for $i=1,\dots,m$, where each $X_i$ is a suspension. Fix $i,j \in [m]$ with $i \ne j$ and let $\{i_1,\dots,i_p\} \subseteq [m]  \setminus  \{i,j\}$ be any non-empty proper subset. If $Y_i=Y_j$ is an $H$-space, then there is a relation
\begin{align*}
& \nabla_{(i,j)}  h_w^{\cK_{\Pi}}(h_w(f_{j_1},\dots,f_{j_q},f_i),f_{i_1},\dots,f_{i_p},f_j) \circ \sigma \\
& \qquad + \nabla_{(i,j)}  h_w^{\cK_{\Pi}}(h_w(f_1,\dots,f_{i-1},f_{i+1},\dots,f_j,\dots,f_m),f_i) \circ \tau = 0
\end{align*}
in $\left[ \Sigma^{m-2} X_1 \wedge \cdots \wedge X_m, (\underline{Y},\underline{\ast})^{\partial \Delta[1,\dots,i-1,i+1,\dots,m]}  \right]$, where $[m] \setminus \{i,j,i_1,\dots,i_p\} = \{j_1,\dots,j_q\}$ and $\sigma$ and $\tau$ are appropriate permutation maps.
\end{corollary} 

\begin{proof}
The result follows by applying Theorem~\ref{thm:FoldsofIdentityComplexes} to the $3$-partition
\[
\Pi = \{\{i\},\{j,i_1,\dots,i_p\},\{j_1,\dots,j_q\}\}
\]
of $[m]$.
\end{proof}

%\begin{example}
%Some folded higher Whitehead maps coincide with $2$-ary Whitehead products. For example let $m=5$. Then combining Theorem~\ref{thm:FoldsofIdentityComplexes} with Proposition~\ref{prop:foldedhwmapishwprod}, the folded higher Whitehead maps $\nabla_{(5,1)}  h_w(h_w(f_1,f_2,f_3),f_4,f_5)$ and $\nabla_{(5,1)}  h_w(h_w(f_1,f_2),f_3,f_4,f_5)$ are both identified, up to an appropriate permutation map, with the Whitehead product
%\[
%[h_w(f_1,f_2,f_3,f_4),\iota \circ f_1]
%\]
%for the inclusion map $\iota \colon Y_1 \longrightarrow FW(Y_1,Y_2,Y_3,Y_4)$. 
%\end{example}

We demonstrate Corollary~\ref{cor:HigherToTwo} by identifying new relations in homotopy groups.

\begin{example} \label{ex:TorsionExample}
Let $\Pi = \{\{1\},\{2,\dots,m-1\},\{m\}\}$ and let $f_i \colon S^{p_i} \longrightarrow Y_i$, $p_i \geqslant 2$, be maps for $i=1,\dots,m$ with $f_1 = f_m$ and $p_1 = p_m$ even. Consider the maps $h_w^{\cK_{\Pi}}(h_w(f_1,\dots,f_{m-1}),f_m)$ and $h_w^{\cK_{\Pi}}(h_w(f_2,\dots,f_m),f_1)$. Applying the fold $\psi \colon \{m\} \longrightarrow \{1\}$ of $\cK_{\Pi}$ and using Propositions~\ref{prop:hwSymmetry} and~\ref{prop:foldedhwmapishwprod},
\[
\nabla_{(m,1)} h_w^{\cK_{\Pi}}(h_w(f_1,\dots,f_{m-1}),f_m) = \nabla_{(m,1)} h_w^{\cK_{\Pi}}(h_w(f_2,\dots,f_m),f_1).
\]

On the other hand, by Corollary~\ref{cor:HigherToTwo}, we obtain
\begin{align*}
2 \nabla_{(m,1)} h_w^{\cK_{\Pi}} & (h_w(f_1, \dots,f_{m-1}),f_m) \\
&= \nabla_{(m,1)} h_w^{\cK_{\Pi}} (h_w(f_1,\dots,f_{m-1}),f_m) + \nabla_{(m,1)} h_w^{\cK_{\Pi}} (h_w(f_2,\dots,f_{m}),f_1) = 0
\end{align*}
in $\pi_{p_1 + \cdots + p_m - 2}((\underline{Y},\underline{\ast})^{\partial \Delta})$. Therefore our methods allow us to study, via higher Whitehead maps, elements of homotopy groups which cannot be seen with rational methods. We study when the map $\nabla_{(m,1)} h_w^{\cK_{\Pi}} (h_w(f_1, \dots,f_{m-1}),f_m)$ is $2$-torsion, that is, when it is not null-homotopic. In general, this depends on the properties of the space $Y_i$. To demonstrate this, let $p_i = 2$ for $i=1,\dots,m$, so that $f_i \colon S^2 \longrightarrow Y_i$.

We first suppose that $Y_i = S^2$ and $f_i \colon S^2 \longrightarrow S^2$ are identity maps for $i=1,\dots,m$. While $S^2$ is not an $H$-space, we can still apply the fold $\psi$ to the relation
\[
h_w^{\cK_{\Pi}}(h_w(f_1,\dots,f_{m-1}),f_m) + h_w^{\cK_{\Pi}}(h_w(f_2,\dots,f_m),f_1) \circ + h_w^{\cK_{\Pi}}(h_w(f_1,f_m),f_2,\dots,f_{m-1}) = 0
\]
to obtain that
\[
2 \nabla_{(m,1)} h_w^{\cK_{\Pi}}(h_w(f_1,\dots,f_{m-1}),f_m) + \nabla_{(m,1)} h_w^{\cK_{\Pi}}(h_w(f_1,f_m),f_2,\dots,f_{m-1}) = 0.
\]
By naturality of the higher Whitehead map and Proposition~\ref{prop:foldedhwmapishwprod}, we have
\[
\nabla_{(m,1)} h_w(h_w(f_1,f_m),f_2,\dots,f_{m-1}) = h_w([f_1,f_1],f_2,\dots,f_{m-1}).
\]
Since $[f_1,f_1] = 2 \eta$, where $\eta \in \pi_3(S^2)$ is the Hopf map, we obtain from linearity of the higher Whitehead map that
\begin{equation} \label{eq:FoldedHwIsHopf}
\nabla_{(m,1)} h_w(h_w(f_1,\dots,f_{m-1}),f_m) = - h_w(\eta,f_2,\dots,f_{m-1}).
\end{equation}
Therefore the folded higher Whitehead map $\nabla_{(m,1)} h_w(h_w(f_1,\dots,f_{m-1}),f_m)$ is not null-homotopic since the higher Whitehead map $h_w(\eta,f_2,\dots,f_{m-1})$ is not null-homotopic, see Example~\ref{ex:HopfHigherWhitehead}.

Now suppose that $Y_i = \CP$ and $f_i \colon S^2 \longrightarrow \CP$ is the inclusion of the bottom cell for $i=1,\dots,m$. Applying the map $(S^2,\ast)^{\partial \Delta} \longrightarrow (\CP,\ast)^{\partial \Delta}$ to each side of~\eqref{eq:FoldedHwIsHopf} trivialises the map $h_w(\eta,f_2,\dots,f_{m-1})$, since $\eta = 0$ in $\pi_3(\CP)$. Therefore by naturality of the higher Whitehead map, the folded map $\nabla_{(m,1)} h_w(h_w(f_1,\dots,f_{m-1}),f_m)$ is trivial in $\pi_{2m-2}((\CP,\ast)^{\partial \Delta})$.

Alternatively, let $Y_i = \Omega S^3$ and $f_i \colon S^2 \longrightarrow \Omega S^3$ be the suspension for $i=1,\dots,m$. Then since the composite $S^3 \xlongrightarrow{\eta} S^2 \longrightarrow \Omega \Sigma S^2$ is adjoint to $\Sigma \eta$, which generates $\pi_4(S^3)$, then a similar argument to Example~\ref{ex:HopfHigherWhitehead} shows that applying the map $(S^2,\ast)^{\partial \Delta} \longrightarrow (\Omega S^3,\ast)^{\partial \Delta}$ to~\eqref{eq:FoldedHwIsHopf} does not trivialise the map $h_w(\eta,f_2,\dots,f_{m-1})$. Therefore, the folded map $\nabla_{(m,1)} h_w(h_w(f_1,\dots,f_{m-1}),f_m)$ is a $2$-torsion element in $\pi_{2m-2}((\Omega S^3,\ast)^{\partial \Delta})$.
\end{example}

\subsection{Folding and substitution}

By introducing richer combinatorial structure into the identity complex $\cK_{\Pi}$ before folding,  we obtain a broader family of relations. Suppose that $\cK_1,\dots,\cK_m$ are simplicial complexes on $[l_1],\dots,[l_m]$, respectively. Denote by $\cS_{\Pi} = \cK_{\Pi} \langle \cK_1,\dots,\cK_m \rangle$ the substitution complex on the vertex set $[l] = [l_1] \sqcup \cdots \sqcup [l_m]$ of $\cK_1,\dots,\cK_m$ into the identity complex $\cK_{\Pi}$. By analysing different partitions $\Pi$ of $[l]$ and folds $\psi \colon I \longrightarrow J$ of $\cS_{\Pi}$, we obtain two further families of relations among folded higher Whitehead maps.

\begin{comment}
Unlike the case with $\cK_i = \bullet$ for each $i$, for an arbitrary fold $\psi \colon I \longrightarrow J$ of $\cK$ the folded complex $\cK_{\nabla(I,J)}$ cannot be determined in general. The main difficulty is that the simplicial folding operation does not commute with substitution.
\end{comment}

We first consider folds of $\cS_{\Pi}$ which are induced by folds within the $\cK_i$ for $i = 1,\dots,m$.

\begin{lemma} \label{lm:FoldingCommuteLemma}
Let $\cL$ be a simplicial complex on $[m]$ and let $\cK_1,\dots,\cK_m$ be simplicial complexes on $[l_1],\dots,[l_m]$, respectively. Let $[l] = [l_1] \sqcup \cdots \sqcup [l_m]$ and suppose that $\psi \colon I \longrightarrow J$ is a fold of $\cL \langle \cK_1,\dots,\cK_m \rangle$ such that $\psi([l_i] \cap I) \subseteq [l_i]$ for each $i=1,\dots,m$. Then
\[
\cL \langle \cK_1,\dots,\cK_m \rangle_{\nabla(I,J)} = \cL \langle \left( \cK_1 \right)_{\nabla(I,J)},\dots,\left( \cK_m \right)_{\nabla(I,J)} \rangle
\]
where $\left( \cK_i \right)_{\nabla(I,J)} = \cK_i$ for $[l_i] \cap I = \emptyset$.
\end{lemma}

\begin{proof}
We extend the fold $\psi \colon I \longrightarrow J$ to a map $\overline{\psi} \colon [l] \longrightarrow [l]$ by setting $\overline{\psi}(i) = \psi(i)$ if $i \in I$ and $\overline{\psi}(i) = i$ otherwise. Then $\psi$ and $\overline{\psi}$ induce the same fold of simplicial complexes. By definition,
\[
\cL \langle \cK_1,\dots,\cK_m \rangle = \left\{ \bigsqcup_{j \in \tau} \sigma_j \mid \sigma_j \in \cK_j, \tau \in \cL \right\}
\]
and therefore
\begin{equation} \label{eq:RHSAbove}
\cL \langle \cK_1,\dots,\cK_m \rangle_{\nabla(I,J)} = \left\{ \bigsqcup_{j \in \tau} \overline{\psi}(\sigma_j) \mid \sigma_j \in \cK_j, \tau \in \cL \right\}.
\end{equation}
Any simplex of $\left( \cK_i \right)_{\nabla(I,J)}$ is either a simplex of $\cK_i$, or $\psi(\sigma_i)$ for some $\sigma_i \in \cK_i$. In either case, every simplex $\left( \cK_i \right)_{\nabla(I,J)}$ can be written as $\overline{\psi}(\sigma_i)$ for some $\sigma_i \in \cK_i$, and so the right-hand side of~\eqref{eq:RHSAbove} is $\cL \langle \left( \cK_1 \right)_{\nabla(I,J)},\dots,\left( \cK_m \right)_{\nabla(I,J)} \rangle$.
\end{proof}

\begin{theorem} \label{thm:FoldsWitinKi}
Let $\cK_1,\dots,\cK_m$ be simplicial complexes on $[l_1],\dots,[l_m]$, respectively, and let $[l] = [l_1] \sqcup \cdots \sqcup [l_m]$. Let $\Pi = \{P_1,\dots,P_k\}$ be a $k$-partition of $[m]$ and denote by $P_i = \{i_1,\dots,i_{p_i}\}$ and $[m] \setminus P_i = \{j_1,\dots,j_{q_i}\}$. Let $\cS_{\Pi} = \cK_{\Pi} \langle \cK_1,\dots,\cK_m \rangle$ and let $\psi \colon I \longrightarrow J$ be a fold of $\cS_{\Pi}$ such that $\psi([l_i] \cap I) \subseteq [l_i]$ for each $i=1,\dots,m$. 

Let $f_i \colon \Sigma X_i \longrightarrow (\underline{Y},\underline{\ast})^{\cK_i}$ be nested higher Whitehead maps for $i=1,\dots,m$. If each $X_i$ is a suspension, then
\begin{equation} \label{eq:FoldInKiRelation}
\sum_{i=1}^k h_w^{\cS_{\Pi}} \left( h_w \left( \nabla_{(I,J)} f_{j_1},\dots \nabla_{(I,J)} f_{j_{q_i}} \right), \nabla_{(I,J)} f_{i_1},\dots,\nabla_{(I,J)} f_{i_{p_i}} \right) \circ \sigma_i = 0
\end{equation}
in $\left[ \Sigma^{m-2} X_1 \wedge \cdots \wedge X_m,(\underline{Y},\underline{\ast})^{(\cS_{\Pi})_{\nabla(I,J)}} \right]$. Furthermore, all summands are null-homotopic if
\[
\nabla_{(I,J)} f_j \colon \, \Sigma X_j \longrightarrow (\underline{Y},\underline{\ast})^{\left( \cK_j \right)_{\nabla(I,J)}}
\]
is null-homotopic for some $j \in \{1,\dots,m\}$.
\end{theorem}

\begin{proof}
By Lemma~\ref{lm:FoldingCommuteLemma}, $(\cS_{\Pi})_{\nabla(I,J)} = \cK_{\Pi} \langle \left(\cK_1 \right)_{\nabla(I,J)},\dots,\left(\cK_m \right)_{\nabla(I,J)} \rangle$. Therefore by Proposition~\ref{prop:naturality},
\begin{align*}
& \nabla_{(I,J)}  h_w^{\cS_{\Pi}} \left( h_w \left( f_{j_1},\dots,f_{j_{q_i}} \right), f_{i_1},\dots,f_{i_{p_i}} \right) \\
& \qquad \qquad = h_w^{\cS_{\Pi}} \left( h_w \left( \nabla_{(I,J)} f_{j_1},\dots \nabla_{(I,J)} f_{j_{q_i}} \right), \nabla_{(I,J)} f_{i_1},\dots,\nabla_{(I,J)} f_{i_{p_i}} \right)
\end{align*}
establishing the claimed relation. Moreover, if $\nabla_{(I,J)} f_j = 0$ for some $j \in \{1,\dots,m\}$, every term in relation~\eqref{eq:FoldInKiRelation} is null-homotopic by Proposition~\ref{prop:basicnestedhwtriviality}.
%Conversely if $\nabla_{(I,J)} \circ f_j$ is non-trivial for every $j = 1,\dots,m$, the map
%\[h_w^{\cK} \left( h_w \left( \nabla_{(I,J)} \circ f_{j_1},\dots \nabla_{(I,J)} \circ f_{j_{r_j}} \right), \nabla_{(I,J)} \circ f_{i_1},\dots,\nabla_{(I,J)} \circ f_{i_{n_i}} \right)\] is non-trivial for each $i=1,\dots,k$ by Proposition~\ref{prop:MainIdentityIsNonZero}.
%\todo[inline,color=green]{Update to topological triviality reference}
\end{proof}

The second type of folds of $\cS_{\Pi} = \cK_{\Pi} \langle \cK_1,\dots,\cK_m \rangle$ we consider are those which fold $\cK_i$ onto $\cK_j$ for $i \neq j$.

\begin{lemma}
Let $\Pi = \{\{1\},\{2,\dots,m-1\},\{m\}\}$ be a $3$-partition of $[m]$, and define $\cS_{\Pi} = \cK_{\Pi} \langle \cK_1,\dots,\cK_m \rangle$. Assume that $\cK_m$ is isomorphic to a full subcomplex of $\cK_1$, with the isomorphism given by a simplicial map $\psi \colon \cK_m \longrightarrow \cK_1$. Let $I = [l_m]$ and $J = \psi([l_m])$. Then
\[
(\cS_{\Pi})_{\nabla(I,J)} = \partial \Delta \langle \cK_1,\dots,\cK_{m-1} \rangle.
\]
\end{lemma}

\begin{proof}
The subcomplexes $\cK_1$ and $\cK_m$ of $\cS_{\Pi}$ are disjoint since there is no edge between the vertices $\{1\}$ and $\{m\}$ in $\cK_{\Pi}$. Therefore, we rewrite $\cS_{\Pi}$ as
\[
\cS_{\Pi} = \cK_{\Pi} \langle \cK_1,\dots,\cK_m \rangle = \partial \Delta \langle \cK_1 \sqcup \cK_m,\cK_2,\dots,\cK_{m-1} \rangle.
\]
The result follows by Lemma~\ref{lm:FoldingCommuteLemma} since $(\cK_1 \sqcup \cK_m)_{\nabla(I,J)} = \cK_1$.
\end{proof}

%\todo[inline,color=magenta]{Remark why the folds don't need $H$ assumption - done in section 4, how do we reference it here?}

For $i=1,\dots,m$, let $f_i \colon \Sigma X_i \longrightarrow (\underline{Y},\underline{\ast})^{\cK_i}$ be nested higher Whitehead maps. If each $X_i$ is a suspension, then there is a relation
\begin{align} \label{eq:ThirdFoldedIdentity}
    & \nabla_{(I,J)}  h_w^{\cS_{\Pi}} \left( h_w(f_2,\dots,f_{m}),f_1 \right) \circ \sigma_1 + \nabla_{(I,J)}  h_w^{\cS_{\Pi}} \left( h_w(f_1,f_m),f_2,\dots,f_{m-1} \right) \circ \sigma_2 \\
    & \qquad + \nabla_{(I,J)}  h_w^{\cS_{\Pi}} \left( h_w(f_1,\dots,f_{m-1}),f_m \right) = 0 \nonumber
\end{align}
in $\left[ \Sigma^{m-2} X_1 \wedge \cdots \wedge X_m, (\underline{Y},\underline{\ast})^{(\cS_{\Pi})_{\nabla(I,J)}} \right]$. By Proposition~\ref{prop:hwFoldedTrivialityConditions}(i), the second summand is null-homotopic if 
\[
(\partial \Delta \langle \Delta \langle \cK_1,\cK_m \rangle,\cK_2,\dots,\cK_{m-1} \rangle)_{\nabla(I,J)} = \partial \Delta \langle (\cK_1 \ast \cK_m)_{\nabla(I,J)},\cK_2,\dots,\cK_{m-1} \rangle
\]
is a subcomplex of $(\cS_{\Pi})_{\nabla(I,J)} = \partial \Delta \langle \cK_1,\dots,\cK_{m-1} \rangle$. Since $\cK_1 \subseteq (\cK_1 \ast \cK_m)_{\nabla(I,J)}$, this condition is equivalent to $\cK_1 = (\cK_1 \ast \cK_m)_{\nabla(I,J)}$. We therefore obtain the following.

\begin{theorem} \label{thm:GeneralUnconnectedFold}
Let $\Pi = \{\{1\},\{2,\dots,m-1\},\{m\}\}$ be a $3$-partition of $[m]$, and define $\cS_{\Pi} = \cK_{\Pi} \langle \cK_1,\dots,\cK_m \rangle$. Assume that $\cK_m$ is isomorphic to a full subcomplex of $\cK_1$, with the isomorphism given by a simplicial map $\psi \colon \cK_m \longrightarrow \cK_1$. Let $I = [l_m]$ and $J = \psi([l_m])$. Let $f_i \colon \Sigma X_i \longrightarrow (\underline{Y},\underline{\ast})^{\cK_i}$ be nested higher Whitehead maps for $i=1,\dots,m$, with each $X_i$ a suspension. Then
\begin{align*}
    & \nabla_{(I,J)}  h_w^{\cS_{\Pi}} \left( h_w(f_2,\dots,f_{m}),f_1 \right) \circ \sigma_1 + \nabla_{(I,J)}  h_w^{\cS_{\Pi}} \left( h_w(f_1,f_m),f_2,\dots,f_{m-1} \right) \circ \sigma_2 \\
    & \qquad + \nabla_{(I,J)}  h_w^{\cS_{\Pi}} \left( h_w(f_1,\dots,f_{m-1}),f_m \right) = 0
\end{align*}
in $\left[ \Sigma^{m-2} X_1 \wedge \cdots \wedge X_m, (\underline{Y},\underline{\ast})^{(\cS_{\Pi})_{\nabla(I,J)}} \right]$, where the second summand is trivial if $\cK_1 = (\cK_1 \ast \cK_m)_{\nabla(I,J)}$. \qed 
\end{theorem} 

\begin{comment}

The following gives a criterion for identifying when $\cK_1 = (\cK_1 \ast \cK_m)_{\nabla(I,J)}$ in terms of the minimal missing faces of $\cK_1$. 

\begin{lemma} \label{lm:FoldingJoinLemma}
The containment $\cK_1 \subseteq (\cK_1 \ast \cK_m)_{\nabla(I,J)}$ is strict if and only if there is $j \in J$ and $L \in MF(\cK_1)$ such that $j \in L$.
\end{lemma}

\begin{proof}
Suppose there exists $j \in J$ and $L \in MF(\cK_1)$ with $j \in L$. Write $L = \{l_1,\dots,l_k,j\}$ and let $i \in I$ be such that $\psi(i) = j$. Then since $\Delta[l_1,\dots,l_k] \in \cK_1$ and $i \in \cK_m$, the simplex $\Delta[l_1,\dots,l_k,i]$ is in $\cK_1 \ast \cK_m$. Therefore $\Delta[l_1,\dots,l_k,j]$ is a simplex of $(\cK_1 \ast \cK_m)_{\nabla(I,J)}$ but not of $\cK_1$.

Conversely suppose that the claimed containment is strict. Then there is $\sigma \in (\cK_1 \ast \cK_m)_{\nabla(I,J)}$ such that $\sigma \notin \cK_1$. Furthermore, $\sigma$ can always be chosen to be a minimal missing face of $\cK_1$. Since $\sigma \notin \cK_1$, then $\sigma = \psi(\tau)$ for some $\tau \in \cK_1 \ast \cK_m$. In particular, $\sigma \cap J$ is non-empty. Therefore any $j \in \sigma \cap J$ gives a $j \in J$ such that $j \in \sigma$.
\end{proof}

\end{comment}

%\todo[inline,color=green]{As an application of Theorem~\ref{thm:GeneralUnconnectedFold} we recover the motivating example by providing the integral relation in homotopy groups of $\DJK$. Reference for example}

As an application of Theorem~\ref{thm:GeneralUnconnectedFold}, we return to relation~\eqref{eq:ElizavetaRelation}, considered at the start of this section. We show that this relation holds integrally.

\begin{example}
Consider the $3$-partition $\Pi = \{\{1\},\{2,3\},\{4\}\}$. Assume that $\cK_1 = \partial \Delta[1_1,1_2,1_3]$ and $\cK_i = \bullet$ for $i=2,3,4$. Let
\[
f_1 = h_w(f_{1_1},f_{1_2},f_{1_3}) \colon \Sigma^2 X_{1_1} \wedge X_{1_2} \wedge X_{1_3} \longrightarrow (\underline{Y},\underline{\ast})^{\cK_1}
\]
and $f_i \colon \Sigma X_i \longrightarrow Y_i$ for $i=2,3,4$ be maps. Let $I = \{4\}$, $J = \{1_1\}$ and $\psi \colon I \longrightarrow J$ be the fold of $\cS_{\Pi} = \cK_{\Pi} \langle \cK_1,\cK_2,\cK_3,\cK_4 \rangle = \partial \Delta \langle \cK_1 \sqcup \cK_4,\cK_2,\cK_3 \rangle$ on $\{1_1,1_2,1_3,2,3,4\}$ given by $\psi(4) = 1_1$. Suppose that $X_i$ is a suspension for $i=2,3,4$. Then since $(\cK_{1} \ast \cK_4)_{\nabla(I,J)} = \Delta[1_1,1_2,1_3] \ne \cK_1$, by Theorem~\ref{thm:GeneralUnconnectedFold} there is a relation
\begin{align} \label{eq:ElizavetaRelationGeneral}
& \nabla_{(4,1_1)}  h_w^{\cS_{\Pi}}(h_w(f_2,f_3,f_4),h_w(f_{1_1},f_{1_2},f_{1_3}) \circ \sigma_1 + \nabla_{(4,1_1)}  h_w^{\cS_{\Pi}}(h_w(h_w(f_{1_1},f_{1_2},f_{1_3}),f_4),f_2,f_3) \circ \sigma_2 \\
& \qquad + \nabla_{(4,1_1)}  h_w^{\cS_{\Pi}}(h_w(h_w(f_{1_1},f_{1_2},f_{1_3}),f_2,f_3),f_4) = 0 \nonumber 
\end{align}
in $\left[ \Sigma^2 X_1 \wedge X_2 \wedge X_3 \wedge X_4, (\underline{Y},\underline{\ast})^{(\cS_{\Pi})_{\nabla(4,1_1)}} \right]$, where $(\cS_{\Pi})_{\nabla(4,1_1)} = \partial \Delta \langle \partial \Delta[1_1,1_2,1_3],2,3 \rangle$.  

Specialising to the case that $X_i = S^1$ and $Y_i = \CP$ for all $i=1,\dots,m$, we recover integrally a relation detected in the rational homotopy groups of the Davis--Januszkiewicz space $DJ_{\cL}$, where $\cL = (\cS_{\Pi})_{\nabla(4,1_1)}$. In particular, let each $f_j \colon S^2 \longrightarrow \CP$ be the inclusion of the bottom cell and let $u_j \colon S^1 \longrightarrow \Omega \CP$ be the adjoint of $f_j$. Zhuravleva~\cite{Zhur21} showed that
\begin{equation} \label{eq:ElizavetaBrackets}
[[u_2,u_3,u_{1_1}],[u_{1_1},u_{1_2},u_{1_3}]] + [[[u_{1_1},u_{1_2},u_{1_3}],u_2,u_3],u_{1_1}] = 0
\end{equation}
in $\pi_8(\Omega DJ_{\cL}) \otimes \mathbb{Q}$, where $[\cdot,\cdot,\cdot]$ is the triple Lie bracket.

Integrally, by relation~\eqref{eq:ElizavetaRelationGeneral}, the Whitehead products $[h_w^{\cL}(f_2,f_3,f_{1_1}),h_w^{\cL}(f_{1_1},f_{1_2},f_{1_3})]$ and $[h_w^{\cL}(h_w(f_{1_1},f_{1_2},f_{1_3}),f_2,f_3),f^{\cL}_{1_1}]$ differ up to sign by
\begin{equation} \label{eq:Elizaveta2TorsionTerm}
\nabla_{(4,1_1)}  h_w(h_w(h_w(f_{1_1},f_{1_2},f_{1_3}),f_4),f_2,f_3) \in [[h_w^{\cL}(f_{1_1},f_{1_2},f_{1_3}),f^{\cL}_{1_1}],f^{\cL}_2,f^{\cL}_3]
\end{equation}
where $f^{\cL}_i \colon \Sigma S^2 \xlongrightarrow{f_i} \CP_i \longrightarrow DJ_{\cL}$.  
By Example~\ref{ex:TorsionExample}, the folded map $\nabla_{(4,1_1)}  h_w(h_w(f_{1_1},f_{1_2},f_{1_3}),f_4)$ is trivial, and therefore so is map~\eqref{eq:Elizaveta2TorsionTerm}. Hence relation~\eqref{eq:ElizavetaRelationGeneral} recovers relation~\eqref{eq:ElizavetaBrackets} integrally. 

Alternatively, the map~\eqref{eq:Elizaveta2TorsionTerm} may not be null-homotopic for different maps $f_i$. By Example~\ref{ex:TorsionExample}, if each $f_i \colon S^2 \longrightarrow \Omega S^3$ is adjoint to the suspension map on $S^2$, then the folded map $\nabla_{(4,1_1)}  h_w(h_w(f_{1_1},f_{1_2},f_{1_3}),f_4)$ is $2$-torsion. It follows that map~\eqref{eq:Elizaveta2TorsionTerm} is also $2$-torsion. We therefore recover further integral relations which cannot be detected rationally.
\end{example}

\section{Proof of main theorem} \label{sec:newmainproof}

In this section we prove Theorem~\ref{thm:maintheorem}, which we now recall.

\newtheorem*{thm:maintheorem2}{Theorem~\ref{thm:maintheorem}}
\begin{thm:maintheorem2}
For $i=1,\dots,m$, let $f_i \colon \Sigma X_i \longrightarrow Y_i$ be maps. Let $\Pi = \{P_1,\dots,P_k\}$ be a $k$-partition of $[m]$ for $k \geqslant 3$ and denote by  $P_i = \{i_1,\dots,i_{p_i}\}$ and $Q_i = [m] \setminus P_i = \{j_1,\dots,j_{q_i}\}$. If $X_i$ is a suspension for each $i=1,\dots,m$, then
\begin{equation} \label{eq:maintheorem2}
\sum_{i=1}^k h_w^{\cK_{\Pi}} \left( h_w \left( f_{j_1},\dots,f_{j_{q_i}} \right),f_{i_1},\dots,f_{i_{p_i}} \right) \circ \sigma_i = 0
\end{equation}
in $\left[ \Sigma^{m-2} X_1 \wedge \cdots \wedge X_m, (\underline{Y},\underline{\ast})^{\cK_{\Pi}}\right]$, where
\[
\sigma_i \colon \Sigma^{m-2} X_1 \wedge \cdots \wedge X_m \longrightarrow \Sigma^{p_i} (\Sigma^{q_i-2} X_{j_1} \wedge \cdots \wedge X_{j_{q_i}}) \wedge X_{i_1} \wedge \cdots \wedge X_{i_{p_i}}
\]
is induced by the coordinate permutation
\[
X_1 \times \cdots \times X_m \longrightarrow X_{j_1} \times \cdots \times X_{j_{q_i}} \times X_{i_1} \times \cdots \times X_{i_{p_i}}.
\]
\end{thm:maintheorem2}

%\todo[color=green,inline]{GS: Make sure this matches the actual statement of the theorem in section 5, and also the statement in the introduction. MS: Done}

We establish Theorem~\ref{thm:maintheorem} by generalising the proof of the Jacobi identity \eqref{eq:JacobiIdentity} and relation~\eqref{eq:HardieIdentity2} due to Nakaoka--Toda~\cite{Na-To54} and Hardie~\cite{Hard64}, respectively. We extend their methods in two ways, first by using the combinatorial structure of $\cK_{\Pi}$ to detect the form of the nested higher Whitehead maps appearing in relation~\eqref{eq:maintheorem2}, and second by deriving the claimed relations for general maps $f_i \colon \Sigma X_i \longrightarrow Y_i$.

We begin by outlining the proof. To prove the Jacobi identity, Nakaoka and Toda studied relative Whitehead products, defined in Blakers--Massey~\cite{BlakersMassey53}, appearing in the long exact sequence of relative homotopy groups. We adapt their methods to higher Whitehead maps.

%Let $W = \Sigma^{m-2} X_1 \wedge \cdots \wedge X_m$, $X = FW(Y_1,\dots,Y_m)$, and let $A = (\underline{Y},\underline{\ast})^{\cK}$. 
We will define relative higher Whitehead maps, which are a generalisation of the relative Whitehead product. In the long exact sequence in homotopy for a pair $(X,A)$
\begin{equation} \label{eq:MainLES}
    \begin{tikzcd}
    \cdots \ar[r] & \left[\Sigma W,A\right] \ar[r,"i"] & \left[ \Sigma W,X \right] \ar[r,"j"] & \left[ (CW,W),(X,A) \right] \ar[r,"\partial"] & \left[ W,A \right] \ar[r] & \cdots
    \end{tikzcd}
\end{equation}
%\todo[color=green, inline]{MS: Reviewer asked for clarity on what this long exact sequence is. Added the map on the left, induced by inclusion A into X, per reviewe. George or Jelena please check.}
the map $j$ is the composite of the isomorphism $[\Sigma W,X] \longrightarrow [(CW,W),(X,\ast)]$ with the map induced by the inclusion $(X,\ast) \longrightarrow (X,A)$, and the map $\partial$ is the restriction sending $[f] \in [(CW,W),(X,A)]$ to $[f \vert_W] \in [W,A]$.
We show that relative higher Whitehead maps have the property that 
\begin{equation} \label{eq:RelWHPProperty} 
\partial h_w(g,g_1,\dots,g_n) = h_w(g \vert_{\Sigma X},g_1,\dots,g_n)
\end{equation}
where $\partial$ is the boundary map in~\eqref{eq:MainLES}, and $g$ is a map of pairs. 

Specialising to $W = \Sigma^{m-2} X_1 \wedge \cdots \wedge X_m$, $X = FW(Y_1,\dots,Y_m)$, and $A = (\underline{Y},\underline{\ast})^{\cK_{\Pi}}$, we show that the higher Whitehead map $h_w(f_1,\dots,f_m) \in [\Sigma W,X]$ when composed with the map $j$ decomposes as the sum
\begin{equation} \label{eq:jMapProperty}
j  h_w(f_1,\dots,f_m) = \sum_{i=1}^k \varphi_i  \sigma_i
\end{equation}
where $\{\varphi_i\}$ is a family of relative higher Whitehead maps determined by the combinatorics of the partition $\Pi$, whose image under the boundary map in the long exact sequence \eqref{eq:RelWHPProperty} is
\begin{equation} \label{eq:RelWHPPropertyspecific}
\partial \varphi_i = h_w^{\cK_{\Pi}} \left( h_w \left( f_{j_1},\dots,f_{j_{q_i}} \right),f_{i_1},\dots,f_{i_{p_i}} \right).
\end{equation}
The exactness of sequence \eqref{eq:MainLES} implies that 
\[\sum_{i=1}^k h_w^{\cK_{\Pi}} \left( h_w \left( f_{j_1},\dots,f_{j_{q_i}} \right),f_{i_1},\dots,f_{i_{p_i}} \right) \circ \sigma_i = \partial (j (h_w(f_1,\dots,f_m)))=0\]
 which proves Theorem~\ref{thm:maintheorem}.
 
% \todo[inline,color=green]{Are we using $\simeq$ or are we using $=$, and (related) are $j$ and $\partial$ actual maps or are they just the maps of htpy classes in the LES?)}
%\todo[inline,color=green]{Equation \eqref{eq:RelWHPProperty} needs to be altered to reflect more precisely what we actually do with relative higher Whitehead maps}

%\todo[inline, color=green]{Change dimension "r" to dimension "d"?}

\subsection{The relative higher Whitehead map}

We define the relative higher Whitehead map as a generalisation of the relative Whitehead product, which we begin by recalling. For a sphere $S^d$, let $D^d_+$ and $D^d_-$ denote the upper and lower hemispheres, respectively. Let $\alpha \in \pi_p(X,A)$ and $\beta \in \pi_q(A)$ be represented by $f \colon (D^p,S^{p-1},D^{p-1}_+) \longrightarrow (X,A,\ast)$ and $g \colon (D^q,S^{q-1}) \longrightarrow (A,\ast)$, respectively. The relative Whitehead product $[\alpha,\beta] \in \pi_{p+q-1}(X,A)$ is the homotopy class of the map of pairs
\begin{align*}
(D^{p+q-1},S^{p+q-2}) &= (D^p \times S^{q-1} \cup D^{p-1}_+ \times D^q,D^{p-1}_- \times S^{q-1} \cup S^{p-2} \cup D^q) \\
& \longrightarrow (X \vee A,A \vee A) \longrightarrow (X,A)
\end{align*}
where the first map is induced by the product $f \times g$ and the second map is induced by the fold $\nabla \colon X \vee X \longrightarrow X$.

The relative Whitehead product has the properties of naturality and bilinearity, see~\cite{BlakersMassey53}, analogously with the Whitehead product. It also satisfies the following. Let $\partial \colon \pi_n(X,A) \longrightarrow \pi_{n-1}(A)$ denote the boundary operator in the long exact sequence for the pair $(X,A)$
\begin{equation*}
    \begin{tikzcd}
    \cdots \ar[r] & \pi_{n}(X) \ar[r,"j"] & \pi_{n}(X,A) \ar[r,"\partial"] & \pi_{n-1}(A) \ar[r] & \cdots.
    \end{tikzcd}
\end{equation*} 
Then
\begin{equation} \label{eq:BlakersMasseyRelativeIdentity}
    \partial_{p+q-1} [\alpha,\beta] = -[\partial_p \alpha, \beta].
\end{equation}
This is a key property used by Nakaoka--Toda~\cite{Na-To54} in the proof of the Jacobi identity. 

To analyse relations between higher Whitehead maps, we introduce the relative higher Whitehead map as an analogue of the relative Whitehead product.

%\todo[inline,color=green]{In the following, the maps $f_i$ are not the same as the maps $f_i$ in the statement of the theorem. We should relabel maybe using another letter}

%\begin{definition}Let $W = X \ast X_1 \ast \cdots \ast X_m$. Let $(\cL,\cK)$ be a simplicial pair on $[l] \sqcup [l_1] \sqcup \cdots \sqcup [l_m]$ such that $\partial \Delta \langle \overline{\cL},\cK_1,\dots,\cK_m \rangle \subseteq \cL$ and $\partial \Delta \langle \overline{\cK},\cK_1,\dots,\cK_m \rangle \subseteq \cK$. The \textit{relative higher Whitehead product map} of $f,f_1,\dots,f_m$ over the simplicial pair $(\cL,\cK)$ is the homotopy class of the composite
%\begin{align*}
%h_w^{(\cL,\cK)}(f,f_1,\dots,f_m) \colon CW & \simeq CX \ast X_1 \ast \cdots \ast X_m  \xlongrightarrow{\rho} FW(\Sigma C X, \Sigma X_1, \dots, \Sigma X_m)\\ &\simeq FW(C \Sigma X, \Sigma X_1, \dots, \Sigma X_m) \\
%& \longrightarrow FW((\overline{\underline{Y}}^l,\underline{\ast})^{\overline{\cL}},(\underline{Y}^{l_1},\underline{\ast})^{\cK_1}, \dots , (\underline{Y}^{l_m},\underline{\ast})^{\cK_m}) \\& = (\underline{Y},\underline{\ast})^{\partial \Delta \langle \overline{\cL},\cK_1,\dots,\cK_m \rangle} \hooklongrightarrow (\underline{Y},\underline{\ast})^{\cL}.
%\end{align*}
%\end{definition}
\goodbreak
\begin{definition}
Let $g \colon (C\Sigma X,\Sigma X) \longrightarrow (Z,B)$ be a map of pairs, let $g_i \colon \Sigma X_i \longrightarrow Y_i$ be maps for $i=1,\dots,n$, and let $W = X \ast X_1 \ast \cdots \ast X_n$. The \textit{relative higher Whitehead map} of $g,g_1,\dots,g_n$ is the composite
\begin{align*}
h_w(g,g_1,\dots,g_n) \colon CW & \simeq CX \ast X_1 \ast \cdots \ast X_n \\
& \xlongrightarrow{\rho} FW(\Sigma CX,\Sigma X_1,\dots,\Sigma X_n) \\
& \xlongrightarrow{\simeq} FW(C \Sigma X,\Sigma X_1,\dots,\Sigma X_n) \\
& \longrightarrow FW(Z,Y_1,\dots,Y_n)
\end{align*}
where $\rho \colon CX \ast X_1 \ast \cdots \ast X_n\longrightarrow FW(\Sigma CX,\Sigma X_1,\dots,\Sigma X_n)$ is the quotient map defined in \eqref{eq:hwmapdecomp}, and the last map is induced by $g,g_1,\dots,g_n$.
\end{definition}

The restriction of $h_w(g,g_1,\dots,g_n)$ to $W$,
\begin{align} \label{eq:restrictionofrelhwmap}
    W & = X \ast X_1 \ast \cdots \ast X_n \\
    & \xlongrightarrow{\rho} FW(\Sigma X,\Sigma X_1,\dots,\Sigma X_n) \nonumber\\
    & \longrightarrow(B,Y_1,\dots,Y_n) \nonumber
\end{align}
is the higher Whitehead map $h_w(g \vert_{\Sigma X},g_1,\dots,g_n)$. The homotopy class of the relative higher Whitehead map $h_w(g,g_1,\dots,g_n)$ is therefore an element of the relative homotopy group
\[
\left[ (CW,W),(FW(Z,Y_1,\dots,Y_n),FW(B,Y_1,\dots,Y_n)) \right].
\]

The higher and relative higher Whitehead maps satisfy an analogous relation to~\eqref{eq:BlakersMasseyRelativeIdentity}.

%\todo[inline,color=green]{As with the definition above, the maps $f_i$ in the following proposition are not the same as in the statement of the theorem}

\begin{proposition} \label{prop:RelOrdhw}
The relative higher Whitehead map $h_w(g,g_1,\dots,g_n)$ satisfies
\[
\partial h_w(g,g_1,\dots,g_n) = h_w(g \vert_{\Sigma X},g_1,\dots,g_n)
\]
where
\[
\partial \colon [(CW,W),(FW(Z,Y_1,\dots,Y_n),FW(B,Y_1,\dots,Y_n)] \longrightarrow [(W,FW(B,Y_1,\dots,Y_n)]
\]
is the boundary map in the exact sequence~\eqref{eq:MainLES}.
\end{proposition}

\begin{proof}
Since
\[
\partial h_w(g,g_1,\dots,g_n) = h_w(g,g_1,\dots,g_n) \vert_{W}
\]
the result follows immediately from the definition of the relative higher Whitehead map $h_w(g,g_1,\dots,g_n)$, see~\eqref{eq:restrictionofrelhwmap}.
\end{proof}
%\todo[inline,color=green]{MS: Added a reference back to equation 54, to see how the restriction to W is what we say it is.}
\subsection{The inclusion map $j$}

We now focus on the map $j \colon [\Sigma W, X] \longrightarrow [(CW,W),(X,A)]$. The following lemma enables us to identify the image under $j$ of certain maps with the sum of relative maps.

\begin{lemma} \label{lm:jmap}
Let $W$ be a suspension space. For $i = 1,\dots,k$, consider pairs $(F_i,Z_i)$ with $F_i \subseteq \Sigma W$ satisfying the following conditions. For each $i = 1,\dots, k$, $Z_i$ is closed in $F_i$, pairwise intersections are such that $F_i \cap F_{i'} = Z_i \cap Z_{i'}$ for all $i' \ne i$, and there exists a homotopy equivalence $e_i \colon \Sigma W \longrightarrow \Sigma W$ such that $e_i(F_i) = C_+W$, $e_i(Z_i) = W$, and $e_i((\Sigma W \backslash F_i) \cup Z_i) = C_-W$.

Let $(X,A)$ be a pair and suppose that $f \colon \Sigma W \longrightarrow X$ satisfies $f((\Sigma W \backslash \bigcup_i F_i) \cup (\bigcup_i Z_i)) \subseteq A$. Then
\[
j (f) = \sum_{i=1}^k f_i
\]
where $f_i = f|_{F_i} \in [(F_i,Z_i),(X,A)] \cong [(CW,W),(X,A)]$ for each $i=1,\dots,k$ and $j \colon [\Sigma W,X] \longrightarrow [(CW,W),(X,A)]$ is the map in the long exact sequence~\eqref{eq:MainLES}.
\end{lemma}

To establish Lemma~\ref{lm:jmap}, we first establish the following lemma.

%\todo[color=green,inline]{Note that we are using $j$ for subscript indices as well as the map $j$. George fixed.}

\begin{lemma} \label{lm:jsum}
    Let $Z = C_1 \widetilde{Z} \cup C_2 \widetilde{Z}$ and write $\Sigma Z = C_+ Z \cup C_- Z$. Let $(X,A)$ be a $CW$-pair. Let $g \colon \Sigma Z \longrightarrow X$ be a map such that $g(Z) \subseteq A \subseteq X$. Then $j(g) = g \vert_{C_+ Z} + g \vert_{C_- Z}$.

    %Moreover, if $g(C_+ \widetilde{Z}) \subseteq A$, then $g \vert_{C_+ Z} = g \vert_{C_+ C_1 \widetilde{Z}} + g \vert_{C_+ C_2 \widetilde{Z}}$.
\end{lemma}

\begin{proof}
    The map $j$ is the composite of the isomorphism $\psi \colon [\Sigma Z,X] \longrightarrow [(CZ,Z),(X,\ast)]$ with the map $j'$ induced by the inclusion $(X,\ast) \longrightarrow (X,A)$. The isomorphism $\psi$ is given by $\psi(g) = g \phi$, where $\phi$ is the map of pairs $\phi \colon (CZ,Z) \longrightarrow (\Sigma Z,\ast)$ which collapses the copy of $Z \subseteq CZ$ at the base of the cone to the basepoint. 

    Now write $C Z = C \Sigma \widetilde{Z}$ and define two decompositions $CZ = G_1 \cup G_2 = G_3 \cup G_4$ where
    \begin{align*}
        & G_1 = \left\{ (s,t,z) \in C \Sigma \widetilde{Z} \mid 0 \leqslant t \leqslant \frac{1}{2} \right\}, && G_2 = \left\{ (s,t,z) \in C \Sigma \widetilde{Z} \mid \frac{1}{2} \leqslant t \leqslant 1 \right\} \\
         & G_3 = \left\{ (s,t,z) \in C \Sigma \widetilde{Z} \mid \frac{1}{2} \leqslant s \leqslant 1 \right\}, && G_4 = \left\{ (s,t,z) \in C \Sigma \widetilde{Z} \mid 0 \leqslant s \leqslant \frac{1}{2} \right\}.
    \end{align*}
    
%    \todo[inline, color=green]{MS: Reviewer wants us to define the map $\phi$ properly (rather than stating what the images of the subspaces are). I suggest we circumvent this by rewording to say that "$\phi$ can be chosen to be a map of pairs such that...". Have not made this change myself.GS: done, M or J check please.}
    
By construction, $\phi$ maps the copy of $C \tilde{Z}$ at $t=\frac{1}{2}$ to a copy of $Z = C_1 \tilde{Z} \cup C_2 \tilde{Z}$ such that $\phi (G_1) = C_+ Z$, $\phi (G_2) = C_- Z$, $\phi(G_3) = \Sigma C_1 \widetilde{Z}$ and $\phi(G_4) = \Sigma C_2 \widetilde{Z}$. In particular, $\phi$ sends $G_1 \cap G_2$ to $Z \subseteq \Sigma Z$, $G_1 \cap G_3$ to $C_+ C_1 \widetilde{Z}$ and $G_1 \cap G_4$ to $C_+ C_2 \widetilde{Z}$. Since $g(Z) \subseteq A$, then $\psi(g) = g \phi$ is a map such that $G_1 \cap G_2 \subseteq A$. Since $G_1 \cap G_2$ is contractible, there is a homotopy from $g \phi$ to a map sending $G_1 \cap G_2$ to the basepoint. Therefore $j(g) = j'( \psi g)$ is homotopic to the sum $g \phi \vert_{G_1} + g \phi \vert_{G_2} = g \vert_{C_+ Z} + g \vert_{C_- Z}$.

    %If further $g(C_+ \widetilde{Z}) \subseteq A$, since $G_1 \cap G_3 \cap G_4$ is contractible, there is a homotopy from $g \vert_{C_+Z} = g \phi \vert_{G_1}$ to $g \phi \vert_{G_1 \cap G_3} + g \phi \vert_{G_1 \cap G_4} = g \vert_{C_+ C_1 \widetilde{Z}} + g \vert_{C_+ C_2 \widetilde{Z}}$. 
\end{proof}
\begin{comment}
\begin{lemma} \label{lm:DefRetReps}
Consider a map of pairs $f \colon (X,A) \longrightarrow (Y,B)$, and subspaces $A' \subseteq X' \subseteq X$, such that the inclusion $X' \xlongrightarrow{i} X$ is homotopic to a map $X' \xlongrightarrow{\iota} X$ where $\iota(A') \subseteq A$, and $\iota \colon (X',A') \longrightarrow (X,A)$ is a homotopy equivalence of pairs. Suppose that $f((X \setminus A) \cup A') \subseteq B$. Then $f \circ \iota \simeq f \vert_{X'}$ are homotopic as maps of pairs.
\end{lemma}

\begin{proof}
The fact that $f \circ \iota$ and $f \vert_{X'} = f \circ i$ are homotopic as maps follows immediately from the assumptions. In order to show that they are homotopic as maps of pairs, it remains to show that there exists a homotopy $H_t \colon X' \longrightarrow Y$ such that $H_t(A') \subseteq B$ for all $t \in I$. Such a homotopy is given by $H_t = f \circ H'_t$, where $H'_t \colon X' \longrightarrow X'$ denotes the homotopy of pairs associated to the homotopy equivalence $\iota \colon (X',A') \longrightarrow (X,A)$. Since $H'_t(A') \subseteq A'$ by the fact that $H'_t$ is itself a homotopy equivalence of pairs, and $f((X \setminus A) \cup A') \subseteq B$, it follows that $H_t(A') \subseteq B$ for all $t \in I$.
\end{proof}
\end{comment}
%\todo[inline]{The statement of this lemma is in general, and again the maps $f_i$ are not actually the maps from the original statement of the theorem.}

We now prove Lemma~\ref{lm:jmap}.

\begin{proof}[Proof of Lemma~\ref{lm:jmap}]

Write $\Sigma W = F_k \cup_{Z_k} ((\Sigma W \setminus F_k) \cup Z_k$. Pre-composing with the homotopy equivalence $e_k$, the map $f$ is homotopic to a map $f' = f \circ e_k \colon \Sigma W \longrightarrow X$, where $f'(W) \subseteq A$ since $f(Z_k) \subseteq A$. Then, by Lemma~\ref{lm:jsum}, $j(f) = j(f \circ e_k) = j(f') = f' \vert_{C_-W} + f'\vert_{C_+W} = f \vert_{(\Sigma W \setminus F_k) \cup Z_k} + f \vert_{F_k}$. Since $f((\Sigma W \backslash \bigcup_i F_i) \cup (\bigcup_i Z_i)) \subseteq A$, collapsing $(\Sigma W \backslash \bigcup_i F_i) \cup (\bigcup_i Z_i)$ gives a homotopy of pairs between the map\\ $f\vert_{(\Sigma W \setminus F_k) \cup Z_k} \colon (CW,W) \longrightarrow (X,A)$ and $\sum_{i=1}^{k-1} f\vert_{F_i}$.

\end{proof}

\subsection{The case for the higher Whitehead map}

Let $\Pi$ be the partition of $[m]$ in the statement of Theorem~\ref{thm:maintheorem}. We now specialise to the $j$-image of the higher Whitehead map $h_w(f_1,\dots,f_m)$ and identify it as a sum of relative higher Whitehead maps.

Write $\Pi = \{I_1,\dots,I_k\}$, with $I_i = \{i_1,\dots,i_{p_i}\}$ and $[m] - I_i = \{j_1,\dots,j_{q_i}\}$, let $\psi = f_1 \times \cdots \times f_m$ and $\psi_i = f_{j_1} \times \cdots \times f_{j_{q_i}}$ for $i=1,\dots,k$. Then $\psi$ is a map whose restriction to $X_1 \ast \cdots \ast X_m$ is $h_w(f_1,\dots,f_m)$ and $\psi_i$ is a map whose restriction to $X_{j_1} \ast \cdots \ast X_{j_{q_i}}$ is $h_w \left(f_{j_1},\dots,f_{j_{q_i}} \right)$.

Let $W = \Sigma^{m-2} X_1 \wedge \cdots \wedge X_m$ and $W_{Q_i} = \Sigma^{q_i-2} X_{j_1} \wedge \cdots \wedge X_{j_{q_i}}$. The domain of the relative higher Whitehead map $\varphi_i = h_w(\psi_i,f_{i_1},\dots,f_{i_{n_i }})$ is the pair
\begin{equation} \label{eq:domainofrelhw}
(CW, W) \simeq \left( C W_{Q_i} \ast  X_{i_1} \ast \cdots \ast X_{i_{n_i}}, W_{Q_i} \ast  X_{i_1} \ast \cdots \ast X_{i_{n_i}}\right).
\end{equation}
To apply Lemma~\ref{lm:jmap} to the map $h_w(f_1,\dots,f_m)$, it is sufficient to show the following. Recall the permutation map $\sigma_i$ from the statement of Theorem~\ref{thm:maintheorem}.
\begin{lemma} \label{lm:FiLemma}
    For $i=1,\dots,m$, there exists subspaces $(F_{Q_i}, Z_{Q_i}) \subseteq CX_1 \times \cdots \times CX_m$, with $Z_{Q_i}$ closed in $F_{Q_i}$, and a homotopy equivalence of pairs
    \[
 (CW,W) \longrightarrow \left( C W_{Q_i} \ast  X_{i_1} \ast \cdots \ast X_{i_{n_i}}, W_{Q_i} \ast  X_{i_1} \ast \cdots \ast X_{i_{n_i}}\right) \longrightarrow (F_{Q_i},Z_{Q_i})
\] such that
    \begin{enumerate}[(i)]
    \item $F_{Q_i} \cap F_{Q_{i'}} \subseteq Z_{Q_i} \cap Z_{Q_{i'}}$ for each $i \ne i'$,
    \item The map $\psi$ can be represented by a map such that 
    \[
    \psi \left( \left( X_1 \ast \cdots \ast X_m \backslash \bigcup_i F_{Q_i} \right) \cup \left( \bigcup_i Z_{Q_i} \right) \right) \subseteq (\underline{Y},\underline{\ast})^{\cK_{\Pi}}
    \]
    with each restriction $\psi \vert_{F_{Q_i}}$ being the composite $h_w(\psi_i,f_{i_1},\dots,f_{i_{p_i}}) \circ \sigma_i$.
    \end{enumerate}
\end{lemma}

The establishing of Lemma~\ref{lm:FiLemma} allows us to prove Theorem~\ref{thm:maintheorem}.

\begin{proof}[Proof of Theorem~\ref{thm:maintheorem}]
Let $X = (\underline{Y},\underline{\ast})^{\partial \Delta[1,\dots,m]}$ and $A = (\underline{Y},\underline{\ast})^{\cK_{\Pi}}$. Consider long exact sequence~\eqref{eq:MainLES}
\begin{equation*}
    \begin{tikzcd}
    \cdots \ar[r] & \left[ \Sigma W,X \right] \ar[r,"j"] & \left[ (CW,W),(X,A) \right] \ar[r,"\partial"] & \left[ W,A \right] \ar[r] & \cdots.
    \end{tikzcd}
\end{equation*}
By Proposition~\ref{prop:RelOrdhw}, Lemma~\ref{lm:jmap} and Lemma~\ref{lm:FiLemma},
\begin{align*}
\partial j (\psi \vert_{X_1*\cdots*X_m}) & \simeq \partial \left(\sum_{i=1}^k \psi \vert_{F_{Q_i}} \right) \\ 
& \simeq \sum_{i=1}^k \partial h_w(\psi_i,f_{i_1},\dots,f_{i_{p_i}}) \circ \sigma_i \\
&= \sum_{i=1}^k h_w(\psi_i \vert_{X_{j_1}*\cdots*X_{j_{q_i}}},f_{i_1},\dots,f_{i_{p_i}}) \circ \sigma_i \\
&= \sum_{i=1}^k h_w(h_w(f_{j_1},\dots,f_{j_{q_i}}),f_{i_1},\dots,f_{i_{p_i}}) \circ \sigma_i.
\end{align*}
On the other hand, by exactness, $\partial j$ is the zero homomorphism, completing the proof.
%\todo[color=green,inline]{GS: Did $V^*$ ever get defined? MS: No. Fixed as I have replaced V* with the join. Please George or Jelena checkGS - looks good}
\end{proof}

The remainder of this section is a technical constructive proof of Lemma~\ref{lm:FiLemma}. The method follows that of Hardie~\cite{Hard64}, adapted to work with non-spherical maps, and to detect the richer combinatorial structure coming from identity complexes. 

\subsection{Proof of Lemma~\ref{lm:FiLemma}} \label{sec:FZpairs}

To prove Lemma~\ref{lm:FiLemma}, we explicitly construct pairs $(F_{Q_i},Z_{Q_i})$ together with homotopy equivalences $(CW,W) \longrightarrow (F_{Q_i},Z_{Q_i})$. To do so, we first identify the domain~\eqref{eq:domainofrelhw} of the relative higher Whitehead map $h_w(\psi_i,f_{i_1},\dots,f_{i_{n_i }})$ as a coordinate subspace of $V = CX_1 \times \cdots \times CX_m$ for each $i=1,\dots,k$.

\begin{lemma}
    Let $X_i = \Sigma \tilde{X}_i = C_+ \tilde{X}_i \cup C_-\tilde{X}_i$ be a suspension for $i=1,\dots,m$. Let $V_{Q_i} = CX_{j_1} \times \cdots \times CX_{j_{q_i}}$ and let 
    \[
    U_{Q_i}^- = \bigcup_{n=1}^{q_i} C \Sigma \widetilde{X}_{j_1} \times \cdots \times C_- \widetilde{X}_{j_n} \times \cdots \times C \Sigma \widetilde{X}_{j_{q_i}}
\]
and
\[
    U_{Q_i}^+ = \bigcup_{n=1}^{q_i} C \Sigma \widetilde{X}_{j_1} \times \cdots \times C_+ \widetilde{X}_{j_n} \times \cdots \times C \Sigma \widetilde{X}_{j_{q_i}}.
\]
Then the pair $(CW_{Q_i} \ast X_{i_1} \ast \cdots \ast X_{i_{p_i}}, W_{Q_i} \ast X_{i_1} \ast \cdots \ast X_{i_{p_i}})$
is homotopy equivalent to the pair $(F'_{Q_i},Z'_{Q_i})$ of coordinate subspaces
\begin{align*}
    & \left( V_{Q_i} \times \left( \bigcup_{l=1}^{p_i} CX_{i_1} \times \cdots \times X_l \times \cdots \times CX_{i_{p_i}} \right) \cup U_{Q_i}^- \times \left(CX_{i_1} \times \cdots \times CX_{i_{p_i}} \right) \right., \\
    & \qquad \left. U_{Q_i}^+ \times \left( \bigcup_{l=1}^{p_i} CX_{i_1} \times \cdots \times X_l \times \cdots \times CX_{i_{p_i}} \right) \cup (U_{Q_i}^+ \cap U_{Q_i}^-) \times \left( CX_{i_1} \times \cdots \times CX_{i_{p_i}} \right) \right).
\end{align*}

\end{lemma}

\begin{proof}
    Decompose
\begin{align*}
& (CW_{Q_i} \ast X_{i_1} \ast \cdots \ast X_{i_{p_i}}, W_{Q_i} \ast X_{i_1} \ast \cdots \ast X_{i_{p_i}}) \\
& = \left( C' CW_{Q_i} \times (X_{i_1} \ast \cdots \ast X_{i_{p_i}}) \cup C W_{Q_i} \times C(X_{i_1} \ast \cdots \ast X_{i_{p_i}}), \right. \\
& \left. \qquad \qquad \qquad \qquad \qquad \qquad C' W_{Q_i} \times (X_{i_1} \ast \cdots \ast X_{i_{p_i}}) \cup W_{Q_i} \times C(X_{i_1} \ast \cdots \ast X_{i_{p_i}}) \right).
\end{align*}

The pair $(C(X_{i_1} \ast \cdots \ast X_{i_{p_i}}), X_{i_1} \ast \cdots \ast X_{i_{p_i}})$
is homotopy equivalent to the coordinate subspace
\[
\left( CX_{i_1} \times \cdots \times CX_{i_{p_i}}, \bigcup_{l=1}^{p_i} CX_{i_1} \times \cdots \times X_l \times \cdots \times CX_{i_{p_i}} \right) \subseteq V
\]
and so it remains to identify the spaces $W_{Q_i}, CW_{Q_i}, C'W_{Q_i}$ and $C'CW_{Q_i}$. Both $U_{Q_i}^-$ and $U_{Q_i}^+$ are contractible, since they have the homotopy type of $C_- \widetilde{X}_{j_1} \ast \cdots \ast C_- \widetilde{X}_{j_{q_i}}$ and $C_+ \widetilde{X}_{j_1} \ast \cdots \ast C_+ \widetilde{X}_{j_{q_i}}$, respectively. Moreover, 
%\[
%U_{Q_i}^- \cup U_{Q_i}^+ = \bigcup_{n=1}^{q_i} C \Sigma \widetilde{X}_{j_1} \times \cdots \times X_{j_n} \times \cdots \times C \Sigma \widetilde{X}_{j_{q_i}} =  V_i^*
%\] 
%and
\begin{align*}
    U_{Q_i}^- \cap U_{Q_i}^+ &= \widetilde{X}_{j_1} \times C \Sigma \widetilde{X}_{j_2} \times \cdots \times C \Sigma \widetilde{X}_{j_{q_i}} \\
    & \qquad \qquad \cup C_- \widetilde{X}_{j_1} \times  \left( C_+ \widetilde{X}_{j_2} \times \cdots \times C \Sigma \widetilde{X}_{j_{q_i}} \cup \cdots \cup C \Sigma \widetilde{X}_{j_2} \times \cdots \times C_+ \widetilde{X}_{j_{q_i}} \right) \\
    & \qquad \cup \dots \cup C \Sigma \widetilde{X}_{j_1} \times \cdots \times C \Sigma \widetilde{X}_{j_{q_i-1}} \times \widetilde{X}_{j_{q_i}} \\
    & \qquad \qquad \left( C_+ \widetilde{X}_{j_1} \times \cdots \times C \Sigma \widetilde{X}_{j_{q_i-1}} \cup \cdots \cup C \Sigma \widetilde{X}_{j_1} \times \cdots \times C_+ \widetilde{X}_{j_{q_i-1}} \right) \times C_- \widetilde{X}_{j_{q_i}} \\
    & \simeq \widetilde{X}_{j_1} \ast C_+ \widetilde{X}_{j_2} \ast \cdots \ast C_+ \widetilde{X}_{j_{q_i}} \cup \cdots \cup C_+ \widetilde{X}_{j_1} \ast C_+ \widetilde{X}_{j_2} \ast \cdots \ast \widetilde{X}_{j_{q_i}} \\
    &\simeq \Sigma^{q_i-1} \widetilde{X}_{j_1} \ast \cdots \ast \widetilde{X}_{j_{q_i}} \\
    &\simeq \Sigma^{q_i-2} X_{j_1} \wedge \cdots \wedge X_{j_{q_i}} = W_{Q_i}
\end{align*}
where the homotopy equivalence on the penultimate line is established by induction, starting with the observation that if $B_i$ is a contractible space containing $A_i$ for $i=1,2$ then
\begin{align*}
    (A_1 \ast B_2) \cap (B_1 \ast A_2) &= (CA_1 \times B_2 \cup A_1 \times CB_2) \cap (CB_1 \times A_2 \cup B_1 \times CA_2) \\
    &= CA_1 \times A_2 \cup A_1 \times CA_2 = A_1 \ast A_2
\end{align*}
and therefore $(A_1 \ast B_2) \cup (B_1 \ast A_2) \simeq \Sigma A_1 \ast A_2$.

The claim follows by taking $W_{Q_i} = U_{Q_i}^+ \cap U^-_{Q_i}$, $CW_{Q_i} = U_{Q_i}^+$, $C'W_{Q_i} = U_{Q_i}^-$ and $C'CW_{Q_i} = V_{Q_i}$. 
\end{proof}

The proof continues by applying a homotopy equivalence of pairs to the pairs $(F'_{Q_i},Z'_{Q_i})$, shrinking them inside $V$ until they satisfy the disjointedness condition of Lemma~\ref{lm:FiLemma}(ii). We construct these homotopy equivalences by first constructing homotopy equivalences for the triple $(V_{Q_i}, U_{Q_i}^+, U_{Q_i}^+ \cap U^-_{Q_i})$ in such a way as to create the necessary separation. 

\begin{lemma} \label{lm:GiLemma}
    For each $i = 1,\dots,k$ there exists a homotopy equivalence $V_{Q_i} \longrightarrow G_i$ which restricts to homotopy equivalences between $U_{Q_i}^+$, $U^-_{Q_i}$ and $U^+_{Q_i}\cap U^-_{Q_i}$ and subspaces $L_{Q_i}$, $M_{Q_i}$, and $N_{Q_i}$ of $G_{Q_i}$, respectively, such that for each $i \ne i'$ there exists $i \ne r \neq i'$ such that the projection onto $CX_r$ of $L_{Q_i}$ and the projection of $L_{Q_{i'}}$ onto $CX_r$ are disjoint. 
\end{lemma}

We establish Lemma~\ref{lm:GiLemma} by defining homotopy equivalences between the tuples $(CX_i, X_i, \tilde{X}_i)$, where $X_i = \Sigma \tilde{X}_i$ and subspaces of $CX_i$. We then propagate these equivalences to the tuples $V_{Q_i}$, $U_{Q_i}^+$ and $U_{Q_i}^+ \cap U^-_{Q_i}$. 

Fix $i \in \{1,\dots,m\}$. Since $X_i$ is a suspension, we write $X_i = \Sigma \widetilde{X}_i = C_+ \widetilde{X}_i \cup C_- \widetilde{X}_i$. Consider $C \Sigma \widetilde{X}_i$ as
\[
C \Sigma \widetilde{X}_i = \{(s,t,x) \in I \times I \times \widetilde{X}_i\} / (1,t,x) \sim (s,0,x) \sim (s,1,x) \sim (s,t,*)
\]
so that $\Sigma \widetilde{X}_i$ is the subspace at $s=0$, and $C_+ \widetilde{X}_i$, $C_- \widetilde{X}_i$ and $\widetilde{X}_i$ are the subspaces of $\Sigma \widetilde{X}_i$ with $t \geqslant \frac{1}{2}$, $t \leqslant \frac{1}{2}$, and $t = \frac{1}{2}$, respectively.
\begin{comment}
\begin{align*}
    \widetilde{X}_i &= \left\{(0,(t,x)) \in C \Sigma \widetilde{X}_i \mid t = \frac{1}{2}\right\}\\
C_- \widetilde{X}_i &= \left\{(0,(t,x)) \in C \Sigma \widetilde{X}_i \mid t \leqslant \frac{1}{2}\right\}\\
C_+ \widetilde{X}_i &= \left\{(0,(t,x)) \in C \Sigma \widetilde{X}_i \mid t \geqslant \frac{1}{2}\right\}\\
\Sigma \widetilde{X}_i &= \left\{(0,(t,x)) \in C \Sigma \widetilde{X}_i\right\}.\\
\end{align*} 
\end{comment}

We define two decompositions $C \Sigma \widetilde{X}_i = D^-_i \cup D^+_i = D^1_i \cup D^2_i$, where
\begin{align*}
    D^+_i &= \left\{(s,(t,x)) \in C \Sigma \widetilde{X}_i \mid \frac{1}{2} \leqslant t \right\}  &D^-_i& = \left\{(s,(t,x)) \in C \Sigma \widetilde{X}_i \mid t \leqslant \frac{1}{2} \right\} \\
    D^1_i &= \left\{(s,(t,x)) \in C \Sigma \widetilde{X}_i \mid \frac{3}{4} - t  \leqslant s \right\} &D^2_i& = \left\{(s,(t,x)) \in C \Sigma \widetilde{X}_i \mid s \leqslant \frac{3}{4} - t \right\}.
\end{align*}

In order to create separation of subspaces in $D_i^-$ and $D_i^1$, we further decompose $D_i^1 = E_i^1 \cup C_i^1$ and $D_i^- = E_i^- \cup C_i^-$, where
\begin{align*}
    E^-_i &= \left\{(s,(t,x)) \in D^-_i \mid t \leqslant \frac{1}{8} \right\}
    &E^1_i& = \left\{(s,(t,x)) \in D^1_i \mid t \geqslant \frac{7}{8} \right\} \\  
    C^-_i &= \left\{(s,(t,x)) \in D^-_i \mid t \geqslant \frac{1}{8} \right\}
    &C^1_i& = \left\{(s,(t,x)) \in D^1_i \mid t \leqslant \frac{7}{8} \right\}.
\end{align*}

\begin{figure}[h]
     \centering
     \begin{subfigure}[b]{0.18\textwidth}
        \centering
        \begin{tikzpicture}[scale=0.6]
        \coordinate (1) at (-2,0);
        \coordinate (2) at (2,0);
        \coordinate (3) at (0,-2);
        \coordinate (4) at (0,2);
        \coordinate (5) at (0,0);
        \coordinate (a) at (1,-1);
        \coordinate (b) at (1,1);
        \coordinate (c) at (0,-1.5);
        \coordinate (d) at (0,1.5);
        \coordinate (e) at (-0.75,0);
        \filldraw[fill=red!50] (a) -- (c) -- (e) -- (d) -- (b) -- (4) -- (1) -- (3);
        \filldraw[fill=magenta!50] (a) -- (c) -- (e) -- (d) -- (b) -- (2);
        \draw (1) -- (3) -- (2) -- (4) -- (1) -- (5) -- (3);
        \draw (2) -- (5) -- (4);
        \end{tikzpicture} 
        \caption{$C \Sigma \widetilde{X}_i = D_i^1 \cup D_i^2$}
        \label{fig:y equals x}
     \end{subfigure}
     \hspace{0.03 \textwidth}
     \begin{subfigure}[b]{0.18\textwidth}
        \centering
        \begin{tikzpicture}[scale=0.6]
        \coordinate (1) at (-2,0);
        \coordinate (2) at (2,0);
        \coordinate (3) at (0,-2);
        \coordinate (4) at (0,2);
        \coordinate (5) at (0,0);
        \filldraw[fill=blue!50] (1) -- (3) -- (4);
        \filldraw[fill=cyan!50] (2) -- (3) -- (4);
        \draw (1) -- (3) -- (2) -- (4) -- (1) -- (5) -- (3);
        \draw (2) -- (5) -- (4);
        \end{tikzpicture} 
        \caption{$C \Sigma \widetilde{X}_i = D_i^+ \cup D_i^-$}
     \end{subfigure}
     \hspace{0.03 \textwidth}
     \begin{subfigure}[b]{0.18\textwidth}
        \centering
        \begin{tikzpicture}[scale=0.6]
        \coordinate (1) at (-2,0);
        \coordinate (2) at (2,0);
        \coordinate (3) at (0,-2);
        \coordinate (4) at (0,2);
        \coordinate (5) at (0,0);
        \coordinate (a) at (1,-1);
        \coordinate (b) at (1,1);
        \coordinate (c) at (0,-1.5);
        \coordinate (d) at (0,1.5);
        \coordinate (e) at (-0.75,0);
        \coordinate (f) at (1.5,0.5);
        \coordinate (g) at (1.5,-0.5);
        \filldraw[fill=orange!50] (f) -- (g) -- (2);
        \filldraw[fill=yellow!50] (f) -- (b) -- (d) -- (e) -- (c) -- (a) -- (g);
        \draw (1) -- (3) -- (2) -- (4) -- (1) -- (5) -- (3);
        \draw (2) -- (5) -- (4);
        \end{tikzpicture} 
        \caption{$D_i^1 = E_i^1 \cup C_i^1$}
     \end{subfigure}
     \hspace{0.03 \textwidth}
     \begin{subfigure}[b]{0.18\textwidth}
        \centering
        \begin{tikzpicture}[scale=0.6]
        \coordinate (1) at (-2,0);
        \coordinate (2) at (2,0);
        \coordinate (3) at (0,-2);
        \coordinate (4) at (0,2);
        \coordinate (5) at (0,0);
        \coordinate (a) at (-1.5,0.5);
        \coordinate (b) at (-1.5,-0.5);
        \filldraw[fill=green!50] (a) -- (b) -- (1);
        \filldraw[fill=teal!50] (a) -- (4) -- (3) -- (b);
        \draw (1) -- (3) -- (2) -- (4) -- (1) -- (5) -- (3);
        \draw (2) -- (5) -- (4);
        \end{tikzpicture} 
        \caption{$D_i^- = E_i^- \cup C_i^-$}
     \end{subfigure}
        \caption{Two decompositions of $C \Sigma \widetilde{X}_i$, and further decompositions of $D_i^1$ and $D_i^-$ for $\widetilde{X}_i = S^0$.}
        \label{fig:DecompsOfCSXi}
\end{figure}

%\todo[color=green, inline]{GS: Can I figure out how to draw these in 3D? MS: This is fine for now.}

For $\delta \in \{-,1\}$, we define maps  $\alpha_i^{\delta} \colon C \Sigma \widetilde{X}_i \longrightarrow D_i^{\delta}$ and $\beta_i^{\delta} \colon D_i^{\delta} \longrightarrow E_i^{\delta}$ such that the restrictions of $\alpha_i^{\delta}$ and $\beta_i^{\delta} \alpha_i^{\delta}$ to the subspaces $\Sigma \widetilde{X}_i, C_+ \widetilde{X}_i, C_- \widetilde{X}_i$, and $\widetilde{X}_i$ are homotopy equivalences. We explicitly specify homotopy equivalent copies of $\Sigma \widetilde{X}_i, C_+ \widetilde{X}_i, C_- \widetilde{X}_i$, and $\widetilde{X}_i$ inside $D_i^{\delta}$ and $E_i^{\delta}$. Define

\begin{alignat}{2}
\alpha_i^- &\colon C \Sigma \widetilde{X}_i \longrightarrow D_i^-, \quad \quad &&(s,(t,x)) \mapsto \begin{cases} (2s+2t-2st-1,(\frac{1}{2},x)) & \text{ if } t \geqslant \frac{1}{2} \\ (s,(t,x)) & \text{ otherwise}\end{cases} \nonumber\\ 
\intertext{and}
\alpha_i^1 &\colon C \Sigma \widetilde{X}_i \longrightarrow D_i^1, \quad \quad &&(s,(t,x))  \mapsto \begin{cases} 
 (s,(\frac{3}{4}-s,x))& \text{ if } s \leqslant \frac{3}{4} - t, t \geqslant \frac{1}{2}\\
(1-2t+s,(2t-s-\frac{1}{4},x)) & \text{ if } s \leqslant t - \frac{1}{4},  t \leqslant \frac{1}{2} \\
(\frac{3}{4}-t,(t,x))& \text{ if } t-\frac{1}{4} \leqslant s \leqslant \frac{3}{4} - t\\
(s,(t,x)) & \text{ if } s \geqslant \frac{3}{4} - t.
\end{cases} \nonumber
\end{alignat}
Then $\alpha_i^{\delta}$ takes $\Sigma \widetilde{X}_i, C_- \widetilde{X}_i, C_+ \widetilde{X}_i$ and $\widetilde{X}_i$ to homotopy equivalent spaces $B_i^{\delta}$, $(B_1)^{\delta}_i$, $(B_2)^{\delta}_i$ and $(B_1)^{\delta}_i \cap (B_2)^{\delta}_i$, respectively, where
    \begin{align*}
    (B_1)^-_i &= \left\{(s,(t,x)) \in D^-_i \mid t \leqslant \frac{1}{2}, s=0 \right\}
    & (B_2)^-_i& = \left\{(s,(t,x)) \in D^-_i \mid t = \frac{1}{2} \right\}  \\
    (B_1)^1_i  &= \left\{(s,(t,x)) \in D^1_i \mid s = \frac{3}{4} - t \right\} 
    &(B_2)^1_i &= \left\{(s,(t,x)) \in D^1_i \mid t \geqslant \frac{3}{4}, s=0 \right\}
    \end{align*}
and $B_i^{\delta} = (B_1)_i^{\delta} \cup (B_2)_i^{\delta}$.

Next, define
\begin{align*}
\beta_i^- &\colon D_i^- \longrightarrow E_i^-, \quad \quad &&(s,(t,x)) \mapsto \begin{cases} (s,(\frac{1}{8},x))& \text{ if } t \geqslant \frac{1}{8}\\
(s,(t,x))& \text{ otherwise}
    \end{cases}\\
    \intertext{and}
\beta_i^1 &\colon D_i^1 \longrightarrow E_i^1, \quad \quad &&(s,(t,x)) \mapsto \begin{cases} 
(1-2t+2st,(\frac{7}{8},x))& \text{ if } t \leqslant \frac{1}{2}\\
(s,(\frac{7}{8},x))& \text{ if } \frac{1}{2} \leqslant t \leqslant \frac{7}{8}\\
(s,(t,x)) & \text{ otherwise.}
    \end{cases}\\
\end{align*}
Then the images of $B_i^{\delta}, (B_1)_i^{\delta}, (B_2)_i^{\delta}$ and $(B_1)_i^{\delta} \cap (B_2)_i^{\delta}$ under $\beta_i^{\delta}$ are, respectively, the homotopy equivalent spaces $S_i^{\delta} \cup T_i^{\delta}$, $S_i^{\delta}$, $T_i^{\delta}$ and $R_i^{\delta}$, where
    \begin{align*}
    R^-_i &= \left\{(s,(t,x)) \in D^-_i \mid t = \frac{1}{8}, s = 0 \right\}
    &R^1_i& = \left\{(s,(t,x)) \in D^1_i \mid t = \frac{7}{8}, s = 0 \right\} \\
    S^-_i &= \left\{(s,(t,x)) \in D^-_i \mid t \leqslant \frac{1}{8}, s = 0 \right\}
    &S^1_i& = \left\{(s,(t,x)) \in D^1_i \mid t = \frac{7}{8} \right\} \\
    T^-_i &= \left\{(s,(t,x)) \in D^-_i \mid t = \frac{1}{8} \right\} 
    &T^1_i& = \left\{(s,(t,x)) \in D^1_i \mid t \geqslant \frac{7}{8}, s = 0\right\}.
\end{align*}

\begin{figure}[h]
     \centering
     \begin{subfigure}[b]{0.18\textwidth}
        \centering
        \begin{tikzpicture}[scale=0.6]
        \coordinate (1) at (-2,0);
        \coordinate (2) at (2,0);
        \coordinate (3) at (0,-2);
        \coordinate (4) at (0,2);
        \coordinate (5) at (0,0);
        \coordinate (a) at (1,-1);
        \coordinate (b) at (1,1);
        \coordinate (c) at (0,-1.5);
        \coordinate (d) at (0,1.5);
        \coordinate (e) at (-0.75,0);
        \draw (1) -- (3) -- (2) -- (4) -- (1) -- (5) -- (3);
        \draw (2) -- (5) -- (4);
        \draw[very thick, color=magenta] (a) -- (c) -- (e) -- (d) -- (b);
        \draw[very thick,color=red] (a) -- (2) -- (b);
        \end{tikzpicture} 
        \caption{$ B_i^1\protect\\ =(B_1)_i^1 \cup (B_2)_i^1 $}
     \end{subfigure}
     \hspace{0.03 \textwidth}
     \begin{subfigure}[b]{0.18\textwidth}
        \centering
        \begin{tikzpicture}[scale=0.6]
        \coordinate (1) at (-2,0);
        \coordinate (2) at (2,0);
        \coordinate (3) at (0,-2);
        \coordinate (4) at (0,2);
        \coordinate (5) at (0,0);
        \draw (1) -- (3) -- (2) -- (4) -- (1) -- (5) -- (3);
        \draw (2) -- (5) -- (4);
        \draw [very thick,color=cyan] (3) -- (1) -- (4);
        \draw [very thick,color=blue] (3) -- (5) -- (4);
        \end{tikzpicture} 
        \caption{$B_i^- \protect\\ =(B_1)_i^- \cup (B_2)_i^-$}
        \label{fig:three sin x}
     \end{subfigure}
     \hspace{0.03 \textwidth}
     \begin{subfigure}[b]{0.18\textwidth}
        \centering
        \begin{tikzpicture}[scale=0.6]
        \coordinate (1) at (-2,0);
        \coordinate (2) at (2,0);
        \coordinate (3) at (0,-2);
        \coordinate (4) at (0,2);
        \coordinate (5) at (0,0);
        \coordinate (a) at (-1.5,0.5);
        \coordinate (b) at (-1.5,-0.5);
        \coordinate (c) at (1.5,0.5);
        \coordinate (d) at (1.5,-0.5);
        \draw (1) -- (3) -- (2) -- (4) -- (1) -- (5) -- (3);
        \draw (2) -- (5) -- (4);
        \draw [very thick,color=teal] (a) -- (b);
        \draw [very thick,color=orange] (c) -- (d);
        \draw [very thick,color=green] (a) -- (1) -- (b);
        \draw [very thick,color=yellow] (c) -- (2) -- (d);
        \end{tikzpicture} 
        \caption{$S_i^-$, $T_i^-$, $S_i^1$ and $T_i^1$}
     \end{subfigure}
        \caption{Other subspaces of $D_i^1$ and $D_i^-$ for $\widetilde{X}_i = S^0$.}
        \label{fig:OtherSubspacesofDi}
\end{figure}

Summarising, we have the following.

\begin{lemma} \label{lm:SomePropertiesofPartitions}
There are homotopy equivalences of pairs
\[
(C\Sigma \widetilde{X}_i,\Sigma \widetilde{X}_i) \xlongrightarrow{\alpha_i^{\delta}} (D_i^{\delta},B_i^{\delta}) \xlongrightarrow{\beta_i^{\delta}} (E_i^{\delta},S_i^{\delta} \cup T_i^{\delta}).
\] \qed

\begin{comment}
\todo[color=green,inline]{GS note to self: here's what happens:
\begin{equation*}
    \begin{tikzcd}
        (CX_i, X_i, C_+ \tilde{X}_i, C_- \tilde{X}_i, \tilde{X}_i) \ar{d}{\alpha} \\
        (D_i^{\delta}, B_i^{\delta}, (B_2)_i^{\delta}, (B_1)_i^{\delta}, (B_1)_i^{\delta} \cap (B_2)^{\delta}_i \ar{d}{\beta} \\
        (E_i^{\delta}, S^{\delta}_i \cup T^{\delta}_i, T^{\delta}_i, S^{\delta}_i, R^{\delta}_i)
    \end{tikzcd}
\end{equation*}}
\end{comment}
\end{lemma}

To construct subspaces satisfying Lemma~\ref{lm:GiLemma}, we apply the maps $\alpha_i^{\delta}$ and $\beta^{\delta}_i \alpha_i^{\delta}$ to $(CX_i, X_i, \tilde{X}_i)$ for a choice of $\delta \in \{1, - \}$ depending on $i$ as follows. We use matrix notation to encode the factors in the product. If $k \geqslant 3$ is odd, define the $k \times k$ matrix $H = [\eta(i,j)]$ by
\[
\eta(i,j) = \begin{cases}
    \ast & \text{if } i = j \\ 1 & \text{if } i+j < k+1, \text{ or } i + j = k+1 \text{ and } i < \frac{k+1}{2} \\ - & \text{if } i+j > k+1, \text{ or } i + j = k+1 \text{ and } i > \frac{k+1}{2}.
\end{cases}
\]
If $k \geqslant 4$ is even, define first the $k \times k$ matrix $H' = [\eta'(i,j)]$ by 
\[
\eta'(i,j) = \begin{cases} \ast  & \text{if } i = j \\ 1 & \text{if } i+j < k+1, \text{ or } i+j = k+1 \text{ and } i \leqslant \frac{k}{2}  \\ - & \text{if } i+j > k+1, \text { or } i+j = k+1 \text{ and } i \geqslant \frac{k}{2}+1
\end{cases}
\]
and then define the matrix $H = [\eta(i,j)]$ from $H'$ by swapping the entries $\eta (\frac{k}{2} + 1,\frac{k}{2} - 1) = 1$ and $\eta(\frac{k}{2} + 1,\frac{k}{2}) = -$. 
%\todo[inline,color=green]{Align the ifs with each other}
Figure~\ref{fig:kArrays} shows the matrix $H$, on the left for $k=7$, and on the right for $k=8$.

\begin{figure}[h]
\centering
\begin{subfigure}[b]{0.4\linewidth}
\centering
\[
\left[ \begin{array}{ccccccc}
\ast & 1 & 1 & 1 & 1 & 1 & 1 \\
1 & \ast & 1 & 1 & 1 & 1 & - \\
1 & 1 & \ast & 1 & 1 & - & - \\
1 & 1 & 1 & \ast & - & - & - \\
1 & 1 & - & - & \ast & - & - \\
1 & - & - & - & - & \ast & - \\
- & - & - & - & - & - & \ast 
\end{array} \right]
\]
\caption{$k=7$}
\label{fig:kOddArray}
\end{subfigure}
\hspace{0.1\linewidth}
\begin{subfigure}[b]{0.4\linewidth}
\centering
\[
\left[ \begin{array}{cccccccc}
\ast & 1 & 1 & 1 & 1 & 1 & 1 & 1\\
1 & \ast & 1 & 1 & 1 & 1 & 1 & -\\
1 & 1 & \ast & 1 & 1 & 1 & - & - \\
1 & 1 & 1 & \ast & 1 & - & - & - \\
1 & 1 & - & 1 & \ast & - & - & - \\
1 & 1 & - & - & - & \ast & - & - \\
1 & - & - & - & - & - & \ast & - \\ 
- & - & - & - & - & - & - & \ast 
\end{array} \right]
\]
\caption{$k=8$}
\label{fig:kEvenArray}
\end{subfigure}
\caption{The matrix $H$ for $k=7$ and $k=8$.}
\label{fig:kArrays}
\end{figure}

\begin{comment}
\begin{figure}[h]
\centering
\begin{subfigure}[t]{0.4\linewidth}
\centering
\begin{tikzpicture}
\matrix[matrix of math nodes](m){
\; &&&&&&&&&\\
&\ast & 1 & 1 & 1 & 1 & 1 & 1 & 1 &\\
&1 & \ast & 1 & 1 & 1 & 1 & 1 & - &\\
&1 & 1 & \ast & 1 & 1 & 1 & - & - &\\
&1 & 1 & 1 & \ast & 1 & - & - & - &\\
&1 & 1 & - & 1 & \ast & - & - & - &\\
&1 & 1 & - & - & - & \ast & - & - &\\
&1 & - & - & - & - & - & \ast & - &\\ 
&- & - & - & - & - & - & - & \ast &\\
&&&&&&&&& \;\\
};
\node[fit=(m-5-5)(m-6-6), draw] {};
\node[fit=(m-4-4)(m-7-7), draw] {};
\node[fit=(m-3-3)(m-8-8), draw] {};
\node[fit=(m-2-2)(m-9-9), draw] {};
% \node[fit=(m-1-1)(m-10-10), draw, thick, dotted] {};
\end{tikzpicture}
\caption{$k$ even}
\label{fig:kEvenArray}
\end{subfigure}
\hspace{0.1\linewidth}
\begin{subfigure}[t]{0.4\linewidth}
\centering
\begin{tikzpicture}
\matrix[matrix of math nodes](m){
\; &&&&&&&&\\
&\ast & 1 & 1 & 1 & 1 & 1 & 1& \\
&1 & \ast & 1 & 1 & 1 & 1 & -& \\
&1 & 1 & \ast & 1 & 1 & - & -& \\
&1 & 1 & 1 & \ast & - & - & -& \\
&1 & 1 & - & - & \ast & - & -& \\
&1 & - & - & - & - & \ast & -& \\
&- & - & - & - & - & - & \ast& \\
&&&&&&&& \; \\
};
\node[fit=(m-5-5)(m-5-5), draw] {};
\node[fit=(m-4-4)(m-6-6), draw] {};
\node[fit=(m-3-3)(m-7-7), draw] {};
\node[fit=(m-2-2)(m-8-8), draw] {};
 \node[fit=(m-1-1)(m-9-9), draw, thick, dotted] {};
\end{tikzpicture}
\caption{$k$ odd}
\label{fig:kOddArray}
\end{subfigure}
\caption{Construction of arrays.}
\label{fig:kArrays}
\end{figure}
\end{comment}

We define the space $G_{Q_i}$ to be the image of $V_{Q_i} = \prod_{n=1}^{q_i} CX_{j_n}$ under the product map $\prod_{n=1}^{q_i} \alpha^{\eta(i,\pi(j_n))}_{j_n}$. Explicitly,
\begin{align*}
& G_{Q_i} = \prod_{n=1}^{q_i} D_{j_n}^{\eta(i,\pi(j_n))} = D_{j_1}^{\eta(i,\pi(j_1))} \times \cdots \times D_{j_{q_i}}^{\eta(i,\pi(j_{q_i}))}
\end{align*}
where if $j_i \in P_l$, then $\pi(j_i) = l$. We define $L_{Q_i}$ and $N_{Q_i}$ as the images of $U^+_{Q_i}$ and $U^+_{Q_i} \cap U^-_{Q_i}$, respectively, under the map $\prod_{n=1}^{q_i} \beta^{\eta(i,\pi(j_n))}_{j_n}\alpha^{\eta(i,\pi(j_n))}_{j_n}$. Finally, we define $M_{Q_i}$ as the complement
\[M_{Q_i}= \left(\bigcup_{n=1}^{q_i} D^{\delta}_{j_1} \times \cdots \times B^{\delta}_{j_n} \times \cdots \times D_{j_{q_i}}^{\delta}\right) \setminus L_{Q_i}.\] 

\begin{example} \label{ex:MatrixExample}
Consider the partition $\Pi = \{\{1\},\{2,3\},\{4\}\}$ where we set $P_1 = \{1\}$, $P_2 = \{2,3\}$ and $P_3 = \{4\}$ so that $\pi(1) = 1$, $\pi(2) = \pi(3) = 2$ and $\pi(4) = 3$. Since $\Pi$ is a $3$-partition we choose the $3 \times 3$ matrix $H$ as follows.
\begin{equation*}
\centering
\left[ \begin{array}{cccccccc}
\ast & 1 & 1 \\
1 & \ast & - \\
- & - & \ast
\end{array} \right]
\end{equation*}
This then defines
\[
G_{Q_1} = D_2^{\eta_3(1,\pi(2))} \times D_3^{\eta_3(1,\pi(3))} \times D_4^{\eta_3(1,\pi(4))} = D_2^{\eta_3(1,2)} \times D_3^{\eta_3(1,2)} \times D_4^{\eta_3(1,3)} = D_2^1 \times D_3^1 \times D_4^1
\]
and similarly we obtain
\begin{align*}
    G_{Q_2} &= D_1^{\eta_3(2,1)} \times D_4^{\eta_3(2,3)} = D_1^1 \times D_4^- \\
    G_{Q_3} &= D_1^{\eta_3(3,1)} \times D_2^{\eta_3(3,2)} \times D_3^{\eta_3(3,2)} = D_1^- \times D_2^- \times D_3^-.
\end{align*}
It is also possible to write down explicit expressions for other subspaces. For example,
\[
L_{Q_1} = (T^1_2 \times E^1_3 \times E^1_4) \cup (E^1_2 \times T^1_3 \times E^1_4) \cup (E^1_2 \times E^1_3 \times T^1_4).
\]
\end{example}

%The pair $(G_{Q_i},G_{Q_i}^{\ast})$ is the image of $(V_{Q_i},V_{Q_i}^{\ast})$ under the product of the maps of pairs $\alpha_{j_n}^{\eta(i,\pi(j_n))} \colon (C \Sigma \widetilde{X}_{j_n},\Sigma \widetilde{X}_{j_n}) \longrightarrow (D_{j_n}^{\eta(i,\pi(j_n))}, B_{j_n}^{\eta(i,\pi(j_n))})$. We summarise this observation in the following lemma.
%Since $D_i^{\delta}$ and $B_i^{\delta}$ are deformation retracts of $C X_i$ and $\Sigma X_i$ for all $i$, $G_i$ is a deformation retract of $V_i$. Moreover, we have the following. 

%\begin{lemma} \label{lm:HtpyTypeOfGi}
%For each $i = 1,\dots,k$, there is a homotopy equivalence of pairs $(V_{Q_i},V_{Q_i}^{\ast}) \longrightarrow (G_{Q_i},G_{Q_i}^{\ast})$. \qed
%\end{lemma}

The following property of our construction will later be used to establish properties of the pairwise intersections $F_{Q_i} \cap F_{Q_{i'}}$.

\begin{lemma} \label{lm:GjConsLemma}
Suppose that $1 \leqslant i<j \leqslant k$. Then there exists $r$, different to $i,j$, such that $G_{Q_i}$ contains a factor $D_r^1$ and $G_{Q_j}$ contains a factor $D_r^-$.
\end{lemma}

\begin{proof}
The statement follows from the definition of the matrix $H$ since, given any two rows $i$ and $j$, there exists $r$ different from $i,j$ such that $\eta(i,r) = 1$ and $\eta(j,r) = -$. Then by construction $G_{Q_i}$ has a factor $D_r^1$ and $G_{Q_j}$ has a factor $D_r^-$.
\end{proof}

\begin{example}
To demonstrate Lemma~\ref{lm:GjConsLemma}, consider $G_{Q_1}$, $G_{Q_2}$ and $G_{Q_3}$ from Example~\ref{ex:MatrixExample}. Then $G_{Q_1}$ has a factor $D_4^1$, while $G_{Q_2}$ as a factor $D_4^-$. Similarly, $G_{Q_1}$ has a factor $D_2^1$ while $G_{Q_3}$ has a factor $D_2^-$; and $G_{Q_2}$ has a factor $D_1^1$ while $G_{Q_3}$ has a factor $D_1^-$.
\end{example}

We are now in a position to establish Lemma~\ref{lm:GiLemma}.

\begin{proof}[Proof of Lemma~\ref{lm:GiLemma}]
    As in the case of Lemma~\ref{lm:GjConsLemma}, for each $i \ne i'$, there is $i \ne r \ne i'$ such that $L_{Q_i}$ contains a factor $E^1_r$ and $L_{Q_i'}$ contains a factor $E^-_r$. By construction, $E^1_r \cap E^-_r = \emptyset$, establishing the claim. 
\end{proof}

Let $\rho_{i}$ be the permutation map defined by
\[
\rho_{i} \left(\left(x_{j_1},\dots,x_{j_{q_i}} \right), \left(x_{i_1},\dots,x_{i_{p_i}} \right) \right) = (x_1,\dots,x_m)
\]
and define 
\begin{align*}
F_{Q_i} &= \rho_{i} \left( G_{Q_i} \times \left( \bigcup_{l=1}^{p_i} CX_{i_1} \times \cdots \times X_{i_l} \times \cdots \times CX_{i_{p_i}} \right) \cup L_{Q_i} \times \left( CX_{i_1} \times \cdots \times CX_{i_{p_i}} \right) \right) \\
Z_{Q_i} &= \rho_{i} \left( M_{Q_i} \times \left( \bigcup_{l=1}^{p_i} CX_{i_1} \times \cdots \times X_{i_l} \times \cdots \times CX_{i_{p_i}} \right) \cup N_{Q_i} \times \left( CX_{i_1} \times \cdots \times CX_{i_{p_i}} \right) \right).
\end{align*}

\begin{figure}[h]
    \centering
    \begin{tikzpicture}[scale=2]
    \coordinate (1) at (0,0,0);
    \coordinate (2) at (2,0,0);
    \coordinate (3) at (0,0,2);
    \coordinate (4) at (2,0,2);
    \coordinate (5) at (0,2,0);
    \coordinate (6) at (2,2,0);
    \coordinate (7) at (0,2,2);
    \coordinate (8) at (2,2,2);
    \coordinate (a1) at (0,2,1);
    \coordinate (a2) at (0,1,0);
    \coordinate (a3) at (0,1,1);
    \coordinate (a4) at (2,2,1);
    \coordinate (a5) at (2,1,0);
    \coordinate (a6) at (2,1,1);
    \coordinate (a7) at (0,2,0.25);
    \coordinate (a8) at (0,1.75,0);
    \coordinate (a9) at (2,2,0.25);
    \coordinate (a10) at (2,1.75,0);
    \coordinate (b1) at (1,0,2);
    \coordinate (b2) at (2,1,2);
    \coordinate (b3) at (1,1,2);
    \coordinate (b4) at (1,0,0);
    \coordinate (b5) at (2,1,0);
    \coordinate (b6) at (1,1,0);
    \coordinate (b7) at (1.75,0,2);
    \coordinate (b8) at (2,0.3,2);
    \coordinate (b9) at (1.75,0,0);
    \coordinate (b10) at (2,0.3,0);
    \coordinate (c1) at (0,0,1);
    \coordinate (c2) at (1,0,2);
    \coordinate (c3) at (1,0,1);
    \coordinate (c4) at (0,2,1);
    \coordinate (c5) at (1,2,2);
    \coordinate (c6) at (1,2,1); 
    \coordinate (c7) at (0,0,1.75);
    \coordinate (c8) at (0.25,0,2);
    \coordinate (c9) at (0,2,1.75);
    \coordinate (c10) at (0.25,2,2);
    \filldraw[fill=red!50,opacity=0.5] (a8) -- (a10) -- (6) -- (5);
    \filldraw[fill=red!50,opacity=0.5] (5) -- (a1) -- (a3) -- (a2);
    \filldraw[fill=blue!50,opacity=0.5] (2) -- (b4) -- (b6) -- (b5);
    \filldraw[fill=blue!50,opacity=0.5] (b7) -- (b9) -- (2) -- (4);
    \filldraw[fill=blue!50,opacity=0.5] (b8) -- (b10) -- (2) -- (4);
    \filldraw[fill=green!50,opacity=0.5] (3) -- (c1) -- (c3) -- (c2);
    \filldraw[fill=green!50,opacity=0.5] (c7) -- (c9) -- (7) -- (3);
    \filldraw[fill=green!50,opacity=0.5] (c8) -- (c10) -- (7) -- (3);
    \filldraw[fill=red!50,opacity=0.5] (6) -- (a4) -- (a6) -- (a5);
    \filldraw[fill=green!50,opacity=0.5] (7) -- (c4) -- (c6) -- (c5);
    \filldraw[fill=blue!50,opacity=0.5] (4) -- (b1) -- (b3) -- (b2);
    \filldraw[fill=red!50,opacity=0.5] (a7) -- (a9) -- (6) -- (5);
    \draw (1) -- (2) -- (4) -- (3) -- (1);
    \draw (1) -- (5);
    \draw (2) -- (6);
    \draw (3) -- (7);
    \draw (4) -- (8);
    \draw (5) -- (6) -- (8) -- (7) -- (5);
    %\foreach \point in {1,2,3,4,5,6,7,8,a1,a2,a3,a4,a5,a6,a7,a8,a9,a10,b1,b2,b3,b4,b5,b6,b7,b8,b9,b10,c1,c2,c3,c4,c5,c6,c7,c8,c9,c10}
        %\fill [black] (\point) circle (0.5 pt);
    \end{tikzpicture}
    \caption{The arrangements of the spaces $F_{Q_1}$, $F_{Q_2}$ and $F_{Q_3}$ in $V$ for the partition $\Pi = \{\{1\},\{2\},\{3\}\}$.}
    \label{fig:NiceCube}
\end{figure}

We are finally able to establish Lemma~\ref{lm:FiLemma}.

\begin{proof}[Proof of Lemma~\ref{lm:FiLemma}]
We first prove the first claim of the Lemma. We show that if $i < i'$ then $F_{Q_i} \cap F_{Q_{i'}} \subseteq Z_{Q_i} \cap Z_{Q_{i'}}$. Let
\begin{align*}
    J_{Q_i} &= \rho_i \left(G_{Q_i} \times \bigcup_{l=1}^{p_i} CX_{i_1} \times \cdots \times X_{i_l} \times \cdots \times CX_{i_{p_i}} \right) \\
    K_{Q_i} &= \rho_i \left(L_{Q_i} \times \left( CX_{i_1} \times \cdots \times CX_{i_{p_i}} \right) \right)
\end{align*}
so that $F_{Q_i} = J_{Q_i} \cup K_{Q_i}$. Then
\[
F_{Q_i} \cap F_{Q_{i'}} = (J_{Q_i} \cap J_{Q_{i'}}) \cup (J_{Q_i} \cap K_{Q_{i'}}) \cup (J_{Q_i} \cap K_{Q_{i'}}) \cup (K_{Q_i} \cap K_{Q_{i'}}).
\]
First consider $J_{Q_i} \cap J_{Q_{i'}}$. Since for $l = 1,\ldots,m$, the intersection $D_l^{\delta} \cap X_l \subseteq B_k^{\delta}$, it follows that
\begin{align*}
J_{Q_i} \cap J_{Q_{i'}} & \subseteq \rho_i \left(G_{Q_i}^{\ast} \times \bigcup_{l=1}^{p_i} CX_{i_1} \times \cdots \times X_{i_l} \times \cdots \times CX_{i_{p_i}} \right) \\
& \qquad\cap \rho_{i'} \left(G_{Q_{i'}}^{\ast} \times \bigcup_{l=1}^{p_{i'}} CX_{{i}_1'} \times \cdots \times X_{{i}_l'} \times \cdots \times CX_{{i}_{p_{i'}}'} \right).
\end{align*}
where $G_{Q_i}^{\ast} = L_{Q_i} \cup M_{Q_i}$. By Lemma~\ref{lm:GjConsLemma}, there is $r \ne i,i'$ such that $G_{Q_i}$ contains a factor $D^1_r$ and $G_{Q_{i'}}$ contains a factor $D^-_r$. Moreover, $L_{Q_i}$ contains a factor $E^1_r$ and $L_{Q_{i'}}$ contains a factor $E^-_r$. Since, by construction, $D^1_r \cap E^-_r = \emptyset = D^-_r \cap E^1_r$, we obtain
\[
(J_{Q_i} \cap J_{Q_{i'}}) \cap \rho_i \left(L_{Q_i} \times \bigcup_{l=1}^{p_i} CX_{i_1} \times \cdots \times X_{i_l} \times \cdots \times CX_{i_{p_i}} \right) = \emptyset
\] 
and
\[
(J_{Q_i} \cap J_{Q_{i'}}) \cap \rho_{i'} \left(L_{Q_{i'}} \times \bigcup_{l=1}^{p_{i'}} CX_{{i'}_1} \times \cdots \times X_{{i'}_l} \times \cdots \times CX_{{i'}_{p_{i'}}} \right) = \emptyset.
\]
Therefore, since $G_{Q_i}^* \setminus L_{Q_i} = M_{Q_i}$,
%\todo[inline,color=green]{Check the triple intersection is correct}
\begin{align*}
    J_{Q_i} \cap J_{Q_{i'}} & \subseteq \rho_i \left(M_{Q_i} \times \bigcup_{l=1}^{p_i} CX_{i_1} \times \cdots \times X_{i_l} \times \cdots \times CX_{i_{p_i}} \right) \\
    & \qquad \cap \rho_{i'} \left(M_{Q_{i'}} \times \bigcup_{l=1}^{p_{i'}} CX_{{i}_1'} \times \cdots \times X_{{i}_l'} \times \cdots \times CX_{{i}_{p_{i'}}'} \right) \subseteq Z_{Q_i} \cap Z_{Q_{i'}}.
\end{align*}

Moreover, the conditions $D^1_r \cap E^-_r = \emptyset = D^-_r \cap E^1_r$ imply that $J_{Q_i} \cap K_{Q_i'} = \emptyset = K_{Q_i} \cap J_{Q_{i'}}$. Finally, since $E_r^1 \cap E_r^- = \emptyset$, $K_{Q_i} \cap K_{Q_{i'}} = \emptyset$, the claim is established.

Next we establish the second claim of the Lemma. Since $X_l \subseteq D^+_l \cup D^2_l$, we represent $f_l \colon  \SX_l \longrightarrow Y_l$ by a map of pairs $f_l \colon (C \Sigma \widetilde{X}_l,\Sigma \widetilde{X}_l) \longrightarrow (Y_l,\ast)$ such that $f_l (D_l^+ \cup D_l^2) = \ast$.

Let $\cL = \partial \Delta[1,\dots,m]$ and let $\cL_{i} = \partial \Delta[j_1,\dots,j_{q_i}]$ for $i=1,\dots,k$. Define the map $\psi \colon (V,V^{\ast}) \longrightarrow (\prod_{l=1}^m Y_l,(\underline{Y},\underline{\ast})^{\cL})$ by
\[
\psi(x_1,\dots,x_m) = (f_1(x_1),\dots,f_m(x_m))
\]
and for $i = 1,\dots,k$, define maps $\psi_i \colon (V_{Q_i},V^{\ast}_{Q_i}) \longrightarrow (\prod_{l=1}^{q_i} Y_{j_l},(\underline{Y},\underline{\ast})^{\cL_i})$ by $\psi_i = \psi \vert_{V_{Q_i}}$.

Since $f_l (D^+_l \cup D^2_l) = \ast$ for $l = 1,\dots, q_i$, the map $\psi_i$ takes the complement of $\prod_{l=1}^{q_i} D^1_l \cap D^-_l$ in $V_{Q_i}$ to $(\underline{Y},\underline{*})^{\partial \Delta[j_1,\dots,j_{q_i}]}$.
Equivalently, $\psi_i ((V_{Q_i} \setminus G_{Q_i}) \cup G^*_{Q_i}) \subseteq (\underline{Y},\underline{\ast})^{\cL_i }$. Therefore, $\psi_i \vert_{G_{Q_i}}$ is homotopic as a map of pairs to $\psi_i$.

Now consider 
\begin{align*}
    \psi \rho_i \vert_{\rho_i^{-1} F_{Q_i}} &\colon G_{Q_i} \times \left( \bigcup_{l=1}^{p_i} CX_{i_1} \times \cdots \times X_{i_l} \times \cdots \times CX_{i_{p_i}} \right) \cup L_{Q_i} \times \left( CX_{i_1} \times \cdots \times CX_{i_{p_i}} \right) \\
    & \longrightarrow \psi_i (G_{Q_i}) \times \left(\bigcup_{l=1}^{p_i} \Sigma X_{i_1} \times \cdots \times \ast \times \cdots \times \Sigma X_{i_{p_i}} \right) \cup \ast \times \left(\Sigma X_{i_1} \times \cdots \times \Sigma_{i_{p_i}} \right)  \\
    & \longrightarrow  \prod_{l=1}^{q_i} Y_{j_l} \times (\underline{Y},\underline{\ast})^{\partial \Delta[i_1,\dots,i_{p_i}] } \cup \ast \times \prod_{l=1}^{p_i} Y_{i_l}  \\
    & =  FW \left( \prod_{l=1}^{q_i} Y_{j_l},Y_{i_1},\dots,Y_{i_{p_i}} \right) \\
    & = (\underline{Y},\underline{\ast})^{\partial \Delta \langle \Delta[j_1,\dots,j_{q_i}],i_1,\dots,i_{p_i} \rangle } \longrightarrow (\underline{Y},\underline{\ast})^{\cL}
%    & = (\underline{Y},\underline{\ast})^{\partial \Delta[\hat{j},j_1,\dots,j_{n_j}] \langle \Delta[\hat{j}_1,\dots,\hat{j}_{r_j}] \langle \cK_{\hat{j}_1},\dots,\cK_{\hat{j}_{r_j}} \rangle, \cK_{j_1},\dots,\cK_{j_{n_j}} \rangle} \\
%    & \longrightarrow (\underline{Y},\underline{\ast})^{\cL \langle \cK_1,\dots,\cK_m \rangle}
\end{align*}
and the restriction
\begin{align*}
    \psi \rho_i \vert_{\rho_i^{-1} Z_{Q_i}} &\colon M_{Q_i} \times \left( \bigcup_{l=1}^{p_i} CX_{i_1} \times \cdots \times X_{i_l} \times \cdots \times CX_{i_{p_i}} \right) \cup N_{Q_i} \times \left( CX_{i_1} \times \cdots \times CX_{i_{p_i}} \right) \\
    & \longrightarrow (\underline{Y},\underline{\ast})^{\cL_i} \times (\underline{Y},\underline{\ast})^{\partial \Delta[i_1,\dots,i_{p_i}]} \cup \ast \times \prod_{l=1}^{p_i} Y_{i_l} \\
    & =  FW \left( (\underline{Y},\underline{\ast})^{\cL_i},Y_{i_1},\dots,Y_{i_{p_i}}  \right) \\
    & = (\underline{Y},\underline{\ast})^{\partial \Delta \langle \partial \Delta[j_1,\dots,j_{q_i}],i_1,\dots,i_{p_i} \rangle } \longrightarrow (\underline{Y},\underline{\ast})^{\cK_{\Pi}}.
%    & =  (\underline{Y},\underline{\ast})^{\partial \Delta[\hat{j},j_1,\dots,j_{n_j}] \langle \cL_j \langle \cK_{\hat{j}_1},\dots,\cK_{\hat{j}_{r_j}} \rangle, \cK_{j_1},\dots,\cK_{j_{n_j}} \rangle} \\
%    & \longrightarrow (\underline{Y},\underline{\ast})^{\cK \langle \cK_1,\dots,\cK_m \rangle}.
\end{align*}
Since the complement of $G_{Q_i}$ in $V_{Q_i}$ is mapped into $(\underline{Y},\underline{\ast})^{\cL_i}$ for each $i=1,\dots,k$, a similar calculation also shows that $\psi$ maps the complement of $\bigcup_{i=1}^k F_{Q_i}$ in $X_1 \ast \cdots \ast X_m$ into $(\underline{Y},\underline{\ast})^{\cK_{\Pi}}$.

It follows from the above calculations that $\psi \rho_i \vert_{\rho_i^{-1} F_{Q_i}} \colon (\rho_i^{-1} F_{Q_i}, \rho_i^{-1} Z_{Q_i}) \longrightarrow  ((\underline{Y},\underline{\ast})^{\cL },(\underline{Y},\underline{\ast})^{\cK_{\Pi}})$ is the relative higher Whitehead map $h_w(\psi_i \vert_{G_{Q_i}},f_{i_1},\dots,f_{i_{p_i}})$. Therefore, $\psi \vert_{F_{Q_i}} \colon (F_{Q_i},Z_{Q_i}) \longrightarrow ((\underline{Y},\underline{\ast})^{\cL },(\underline{Y},\underline{\ast})^{\cK_{\Pi}})$ is the composite $h_w(\psi_i \vert_{G_{Q_i}},f_{i_1},\dots,f_{i_{p_i}}) \circ \rho_i^{-1}$. Finally, since $\psi_i \vert_{G_{Q_i}}$ is homotopic to $\psi_i$ and $\rho_i^{-1} = \sigma_i$, the restriction $\psi \vert_{F_{Q_i}}$ is homotopic to $h_w(\psi_i,f_{i_1},\dots,f_{i_{p_i}}) \circ \sigma_i$. 
\end{proof}

\bibliographystyle{abbrv} 
\bibliography{bibliography}

\end{document}